\documentclass[11pt, reqno]{amsart}
\usepackage{amsmath, amsthm, amscd, amsfonts, amssymb, graphicx, color, mathtools}
\usepackage{mathrsfs}
\usepackage[bookmarksnumbered, colorlinks, plainpages]{hyperref}
\usepackage{enumerate}
\usepackage{setspace}
\usepackage{multicol}
\usepackage[margin=1in]{geometry}
\usepackage{comment}
\usepackage{dsfont}
\usepackage{tikz-cd}
\usepackage{algorithmic}

\usepackage[normalem]{ulem}
\newcommand\tsout{\bgroup\markoverwith{\textcolor{red}{\rule[0.5ex]{2pt}{1.4pt}}}\ULon}
\newcommand{\stkout}[1]{\ifmmode\text{\tsout{\ensuremath{#1}}}\else\tsout{#1}\fi}
\allowdisplaybreaks

\theoremstyle{definition}
\newtheorem{theorem}{Theorem}[section]
\newtheorem{lemma}[theorem]{Lemma}
\newtheorem{proposition}[theorem]{Proposition}
\newtheorem{corollary}[theorem]{Corollary}
\newtheorem{definition}[theorem]{Definition}

\newtheorem{remark}[theorem]{Remark}
\numberwithin{equation}{section}

\newcommand{\bb}{\mathbb}

\newcommand{\dt}{\,\mathrm{d}t}

\newcommand{\dr}{\,\mathrm{d}r}
\newcommand{\dx}{\,\mathrm{d}x}
\newcommand{\ds}{\,\mathrm{d}s}

\newcommand{\norm}[2]{\left\|{#1}\right\|_{#2}}
\newcommand{\inpro}[2]{\left\langle#1,#2\right\rangle}

\begin{document}
\setcounter{page}{1}

\title[Beyond the reduced Landau--Lifshitz--Bloch equation]{Beyond the reduced Landau--Lifshitz--Bloch equation: Singular potential, global well-posedness, and micromagnetic limits}

\author[Agus L. Soenjaya]{Agus L. Soenjaya}
\address{Institute of Analysis and Scientific Computing, TU Wien, Wiedner Hauptstrasse 8--10, 1040 Vienna, Austria}
\email{\textcolor[rgb]{0.00,0.00,0.84}{agus.soenjaya@asc.tuwien.ac.at}}

\date{1 October 2026}

\begin{abstract}
We study the full nonlinear mean-field Landau--Lifshitz--Bloch equation in micromagnetics arising from Garanin's spin-polarisation model prior to the usual near-equilibrium reduction. The equation combines a singular inverse-Langevin potential with a mobility whose longitudinal component degenerates at saturation, together with exchange, anisotropy, and the nonlocal whole-space stray field. The natural state space is the closed unit ball $|m|\leq 1$. For any $H^1$ initial data satisfying $|m_0|\leq1$, including saturated data of infinite free energy, we construct global weak solutions and derive instantaneous separation from saturation $|m|=1$. This restores uniform parabolicity and yields the thermodynamic formulation and exact energy identity for positive times. We further establish a quantitative low-temperature limit to the Landau--Lifshitz--Gilbert equation and prove convergence of global trajectories to a single equilibrium. A spectral Curie-type threshold characterises the stability of the paramagnetic state and yields exponential relaxation above the threshold. Finally, for uniformly separated data we prove strong well-posedness, including global strong solutions in one- and two-dimensional domains, and exponential stability of non-degenerate equilibria, with a small-data global theory in three dimensions.
\end{abstract}

\maketitle

\tableofcontents

\section{Introduction}

Classical micromagnetics describes the magnetisation by a vector field of unit length, leading to the Landau--Lifshitz--Gilbert (LLG) equation. This description is appropriate at sufficiently low temperatures,
where longitudinal variations of the magnetisation are negligible, but it becomes inadequate near and above the Curie temperature. In that regime, the magnetisation represents a local thermal average of microscopic spins and its magnitude is no longer fixed. The Landau--Lifshitz--Bloch (LLB) equation introduced by Garanin in~\cite{Gar91, Gar97} removes this constraint and allows the magnetisation to take values in the closed unit ball. It has since been extensively studied in the physics literature~\cite{ChuNowChaGar06, GarChu04, MeoPan20, NieChu16} to model ultrafast heat-assisted magnetisation dynamics relevant in magnetic memory at high temperatures.

Starting from a Fokker--Planck equation for a classical spin, Garanin derives an evolution equation for the averaged spin polarisation $m$ in which the relaxation coefficients depend nonlinearly on $|m|$ through the inverse Langevin function; see~\cite[Eq.~(2.16)]{Gar97}. For a ferromagnet, the magnetic field in this equation is subsequently replaced by a molecular mean field, whose continuum approximation consists of a dominant homogeneous exchange contribution together with the usual spatial exchange, anisotropy, and applied-field contributions.

Motivated by the preceding discussion, we re-examine the nonlinear magnetisation equation originally derived by Garanin~\cite[Eq.~(2.16)]{Gar97}, before the further approximation leading to the reduced macroscopic LLB equation \cite[Eqs.~(4.13)--(4.15)]{Gar97}.
Specifically, after rescaling the physical parameters, we consider a magnetisation $m: [0,\infty)\times\mathscr D\to \bb R^3$, where $|m|\leq 1$, satisfying
\begin{subequations}\label{equ:full nonlinear LLB}
	\begin{align}
		&\partial_t m
		=
		\gamma m\times H_{\mathrm{MFA}}(m)
		+
		\mathbb M(m)H_{\mathrm{MFA}}(m)
		-
		2\theta m
		&&\text{in }(0,\infty)\times\mathscr D,
		\label{equ:intro full nonlinear LLB}
		\\
		&H_{\mathrm{MFA}}(m)
		=
		Jm+\alpha\Delta m+\mathsf A m+H_{\mathrm d}[m]
		&&\text{in }(0,\infty)\times\mathscr D,
		\label{equ: H MFA}
		\\
		&\operatorname{curl}H_{\mathrm d}[m]=0,
		\qquad
		\operatorname{div}
		\bigl(H_{\mathrm d}[m]+m\mathbf 1_{\mathscr D}\bigr)=0
		&&\text{in }(0,\infty)\times\mathbb R^3,
		\label{equ:intro stray field}
		\\
		&\partial_\nu m=0
		&&\text{on } (0,\infty)\times \partial\mathscr D,
		\label{equ:boundary}
		\\
		&m(0)=m_0
		&&\text{in }\mathscr D.
		\label{equ:init}
	\end{align}
\end{subequations}

The constants $\gamma$ and $\alpha>0$ denote the gyromagnetic ratio and the exchange coefficient, respectively, while $\nu$ denotes the outward-pointing normal to $\partial\mathscr{D}$. The term $Jm$ represents the homogeneous molecular exchange, where $J>0$. The matrix $\mathsf A\in\mathbb R^{3\times3}$ contains the local anisotropy contributions and satisfies $\mathsf A=\mathsf A^\top$, while $H_{\mathrm d}[m]$ is the physical stray field. The field $H_{\mathrm{MFA}}(m)$ is the molecular field arising from the mean-field approximation and should be distinguished from the usual micromagnetic field. Accordingly, for emphasis we write
\begin{equation}\label{equ:H MFA mic}
	H_{\mathrm{MFA}}(m) = Jm+H_{\mathrm{mic}}(m), \; \text{ where }\;  H_{\mathrm{mic}}(m) = \alpha\Delta m+\mathsf A m+H_{\mathrm d}[m].
\end{equation}
Here, $H_{\mathrm{mic}}$ is the usual micromagnetic effective field containing the exchange, anisotropy, and magnetostatic contributions~\cite{CarFab01,Pro01}.

The mobility matrix $\mathbb M(m)\in \mathbb R^{3\times3}$ in \eqref{equ:intro full nonlinear LLB} encodes the nonlinear longitudinal and transverse relaxation coefficients in~\cite{Gar97}. Let
\begin{equation}\label{equ:q and L}
	q(s):=\frac{s}{\mathcal L^{-1}(s)}, \quad \text{ where }
	\quad
	\mathcal L(\xi):=\coth\xi-\frac1\xi,
\end{equation}
and $\mathcal L^{-1}$ is the inverse function of $\mathcal{L}$.
We define, for $m\neq0$,
\[
\mathbb M(m)
\coloneq
\bigl(1-q(|m|)\bigr)
\left(I-\frac{m\otimes m}{|m|^2}\right)
+
2q(|m|)
\frac{m\otimes m}{|m|^2}.
\]
Since $q(s)\to 1/3$ as $s\downarrow 0$, the apparent directional singularity is
removable and in fact $\mathbb M(0)\coloneq \frac23 \mathsf I$. On the other hand, $q(s)\to0$ as $s\uparrow1$, so the part of the diffusion acting in the
direction of $m$ degenerates at saturation $|m|=1$.

In Garanin's mean-field description, $m$ is the local average of classical unit spins $s\in \bb{S}^2$. Thus, the natural state space of the LLB model is the closed unit ball $|m|\leq 1$, in contrast to the sphere-valued constraint $|m|=1$ in LLG. Overall, the full nonlinear LLB equation \eqref{equ:full nonlinear LLB} is a \emph{nonlocal quasilinear parabolic system} with \emph{anisotropic degeneracy at saturation} and a \emph{singular inverse-Langevin thermodynamic structure}. A more detailed discussion of \eqref{equ:full nonlinear LLB} and its mathematical structure will be given in Section~\ref{sec:full nonlinear LLB model}.

The LLB equation commonly studied in the mathematical
literature~\cite{HamHam20, Le16, QiuHuaWan26, Soe25b, Soe26} is obtained
only after a further near-equilibrium approximation of
\eqref{equ:full nonlinear LLB}. In regimes where the homogeneous
molecular-exchange field dominates the remaining effective field, the
magnetisation is expanded about the corresponding equilibrium state;
see~\cite[Sec.~IV]{Gar97}. This leads to the reduced equation
\begin{align}
	\partial_t m
	=
	\gamma m\times h_{\mathrm{red}}(m)
	+
	\frac{\beta_\parallel}{|m|^2}
	\bigl(m\cdot h_{\mathrm{red}}(m)\bigr)m
	-
	\frac{\beta_\perp}{|m|^2}
	m\times
	\bigl(
	m\times h_{\mathrm{red}}(m)
	\bigr).
	\label{equ:intro reduced LLB}
\end{align}
A more detailed derivation of this is given in Section~\ref{subsec:relation reduced LLB}. Below the Curie temperature, this formulation is not defined at $m=0$ because of its dependence on $m/|m|$. The weighted formulation used in~\cite{HamHam20} removes the explicit division by $|m|^2$, but degenerates on the zero set and is equivalent to \eqref{equ:intro reduced LLB} only where $m\neq0$. By contrast, the directional singularity of the full nonlinear model is removable and the equation remains well defined at vanishing magnetisation.

In this work, we study the full nonlinear LLB equation \eqref{equ:full nonlinear LLB} before this reduction, which we regard as a more general mean-field starting point whenever the near-equilibrium assumptions underlying \eqref{equ:intro reduced LLB} are not justified. Its analysis is complicated by three interacting features: the longitudinal mobility degenerates as $|m|\uparrow1$, the inverse-Langevin thermodynamic field and free-energy potential are singular at saturation, and the equation is coupled to the nonlocal whole-space stray field. In particular, saturated initial data may have infinite thermodynamic energy, so neither uniform parabolicity nor finite initial energy can be assumed.

Our analysis begins with an equivalent closed-ball formulation, where the singular thermodynamic force and the degenerate mobility cancel exactly. A viscous regularisation and Faedo--Galerkin approximation yield weighted second-order and inverse-distance estimates. Combined with positive-time $W^{2,r}$-regularity and a scalar comparison argument, these give uniform separation from $|m|=1$ for every positive time. Thus solutions may emanate from saturated, possibly infinite-energy,
initial data and enter the finite-energy regime instantaneously. Similar separation mechanisms are central in the analysis of Cahn--Hilliard equations with singular potentials~\cite{AbeWil07,ConGalGatGio25,GalGioGra23}, where separation removes the
singularity and yields higher regularity and long-time analysis. Here, it has the additional role of restoring uniform ellipticity of the principal part. However, the precession term and the nonlocal stray field introduce further difficulties absent from the scalar gradient-flow setting.

Once separation is established, the thermodynamic formulation and exact energy identity become available. We then combine asymptotic compactness with a \L{}ojasiewicz--Simon inequality, in the spirit of~\cite{AbeWil07,ConGalGatGio25}, to prove convergence to a single stationary state. The same structure yields a
spectral Curie-type threshold and exponential relaxation in the paramagnetic regime. We also prove a direct quantitative low-temperature limit from the full nonlinear LLB equation to the Landau--Lifshitz equation, and develop a strong-solution and nonlinear stability theory for uniformly separated data.

The main results of the paper are as follows:
\begin{enumerate}
	\item
	Global existence for arbitrary $m_0\in H^1(\mathscr D)$ with $|m_0|\leq 1$, including saturated initial data; see
	Theorem~\ref{thm:global existence full nonlinear LLB};
	
	\item
	Quantitative low-temperature convergence of the full LLB
	model to the Landau--Lifshitz--Gilbert equation; see
	Theorem~\ref{thm:LLB to LLG};
	
	\item
	Convergence of every global solution to an equilibrium as
	$t\to\infty$; see Theorem~\ref{thm:LLB convergence equilibrium};
	
	\item
	Precise characterisation of the stationary state via the Curie-type threshold $\theta_{\mathrm c}=\Lambda_\ast/3$. Above it, $0$ is the unique
	equilibrium and all solutions converge to it exponentially. Below the threshold, the energy admits nonzero states of lower energy,
	and trajectories attaining negative energy converge to a nonzero
	stationary state; see
	Theorem~\ref{thm:Curie threshold behaviour};
	
	\item
	Local strong well-posedness for uniformly separated initial data $m_0\in H^s(\mathscr{D})$, $s\in (3/2,2]$, and global strong well-posedness for $d\leq2$; see Theorems~\ref{thm:strong separated wellposedness}.
	
	\item 
	Exponential stability of nondegenerate
	equilibria, which implies global strong well-posedness in $d=3$ for sufficiently small initial data; see Theorem~\ref{thm:strong stable equilibrium} and Corollary~\ref{cor:strong small magnetization}.
\end{enumerate}

\section{The full nonlinear mean-field LLB model}
\label{sec:full nonlinear LLB model}

In this section, we summarise the derivation of \eqref{equ:intro full nonlinear LLB} in~\cite{Gar97} and elucidate its mathematical structures.
Let $\mathscr D\subset\mathbb R^3$ be a bounded connected domain with sufficiently smooth boundary. We write $m\colon [0,\infty)\times \mathscr D \to\mathbb R^3$ for the dimensionless spin polarisation. In Garanin's classical-spin description, the microscopic spin $s$ has unit length, $|s|=1$, while $m=\langle s\rangle$ denotes averaging with respect to the local distribution of spin orientations. It follows that the physically admissible range is $|m|\leq 1$.
This differs from the usual LLG convention, in which the magnetisation is constrained to take values in $\mathbb S^2$.

\subsection{The nonlinear spin-polarisation equation}

Let $\mathcal L: [0,\infty) \to [0,1)$ given by
\begin{equation}\label{equ:L xi}
\mathcal L(\xi):=\coth\xi-\frac1\xi,
\qquad \mathcal L(0)\coloneq 0,
\end{equation}
denote the strictly increasing \emph{Langevin function}, whose inverse is denoted by $\mathcal L^{-1}$. We set 
\begin{equation}\label{equ:q s}
	q(s):=\frac{s}{\mathcal L^{-1}(s)}, \qquad s\in (0,1). 
\end{equation}

As argued physically in \cite[Sec.~II]{Gar97}, Garanin introduces a nonequilibrium
reduced field $\xi$ through the relation $|m|=\mathcal L(\xi)$,
and derives the nonlinear LLB equation for the spin polarisation
$m$. In his notation, using the physical relations $\Lambda_N=2\gamma_0\lambda\theta$ and $\xi_0= H/\theta$, where $\theta$ denotes the rescaled temperature, we arrive at \cite[Eq.~(2.16)]{Gar97}:
\begin{align*}
	\partial_t m
	={}&
	\gamma_0 m\times H
	-
	2\gamma_0\lambda\theta m
	+
	\gamma_0\lambda
	\left[
	\bigl(1-q(|m|)\bigr)
	\left(
	\mathsf I-\frac{m\otimes m}{|m|^2}
	\right)
	+
	2q(|m|)
	\frac{m\otimes m}{|m|^2}
	\right] H.
\end{align*}
This is the starting point of our analysis. 

The form of the equation stated above motivates the definition of the \emph{mobility matrix}:
\begin{equation}
	\label{equ:nonlinear LLB mobility}
	\mathbb M(v)
	\coloneq
	\bigl(1-q(|v|)\bigr)
	\left(
	\mathsf I-\frac{v\otimes v}{|v|^2}
	\right)
	+
	2q(|v|)
	\frac{v\otimes v}{|v|^2},
	\qquad
	v\neq0.
\end{equation}
After rescaling time by the positive factor
$\gamma_0\lambda$ and renaming the resulting precession coefficient,
the equation takes the form
\begin{equation}
	\label{equ:Garanin mobility form}
	\partial_t m
	=
	\gamma m\times H
	+
	\mathbb M(m)H
	-
	2\theta m.
\end{equation}
Here, $H$ is a prescribed magnetic field acting on a single-spin ensemble.

\subsection{The continuum molecular mean field}

For a ferromagnet, the model replaces the single-spin magnetic field $H$ in \eqref{equ:Garanin mobility form} by the molecular mean field $H_{\mathrm{MFA}}$. Under the mean-field and weak-anisotropy approximations, this field consists of the homogeneous molecular exchange together with the usual exchange and anisotropy contributions; see~\cite[Sec.~IV]{Gar97}. In the present work, we write
\begin{equation}
	H_{\mathrm{MFA}}(m)
	:=
	\alpha\Delta m+\mathsf Bm+H_{\mathrm d}[m],
	\label{equ:molecular mean field}
\end{equation}
where $\alpha>0$ is the exchange coefficient and
$\mathsf B=\mathsf B^\top\in\mathbb R^{3\times3}$. 

It is useful to regard $\mathsf B=J\mathsf I+\mathsf A$,
where $J>0$ represents the homogeneous molecular exchange and $\mathsf A=\mathsf A^\top$ the local anisotropy. Thus,
\begin{equation*}
H_{\mathrm{MFA}}(m)
=
Jm+H_{\mathrm{mic}}(m), \;\text{ where }\;
H_{\mathrm{mic}}(m)
:=
\alpha\Delta m+\mathsf A m+H_{\mathrm d}[m].
\end{equation*}
In particular, $H_{\mathrm{MFA}}$ should be distinguished from the usual micromagnetic effective field $H_{\mathrm{mic}}$. The term $Jm$ is irrelevant for fixed-length LLG dynamics, since it is parallel to $m$, but it plays an essential role in the longitudinal relaxation of the LLB model.


\subsection{Closed-ball and thermodynamic formulations}

We proceed with a further reformulation of the problem.
First, we note the endpoint behaviour of the coefficient $q$ given in \eqref{equ:q s}. The Taylor expansion of the Langevin function at $0$ gives
\[
\mathcal L^{-1}(s)
=
3s+\frac95s^3+O(s^5),
\qquad
q(s)
=
\frac13-\frac15s^2+O(s^4)
\qquad\text{as }s\to 0,
\]
where $\mathcal{L}$ and $q$ were defined in \eqref{equ:L xi} and \eqref{equ:q s}, respectively.
Therefore, we can extend $q$ continuously by $q(0):=\frac13$.
Furthermore, since $\mathcal L^{-1}(s)\to\infty$ as $s\uparrow1$, we set $q(1):=0$.
It follows from \eqref{equ:nonlinear LLB mobility} that the apparent singularity at $v=0$ is removable, with $\mathbb M(0):=\frac23 \mathsf I$.
Moreover, for $|v|=1$, we have $\mathbb M(v)=\mathsf I -v\otimes v$.
Hence, $\mathbb M$ extends continuously to the closed unit ball
$\overline{B_1(0)}$.

Consequently, the equation
\begin{equation}\label{equ:closed ball LLB}
	\partial_t m
	=
	\gamma m\times H_{\mathrm{MFA}}(m)
	+
	\mathbb M(m)H_{\mathrm{MFA}}(m)
	-
	2\theta m
\end{equation}
is well-defined for all $|m|\leq1$. We shall refer to \eqref{equ:closed ball LLB} as the \emph{closed-ball formulation}. This is the natural formulation for the existence theory, since the initial datum is only assumed to satisfy $|m_0|\leq1$ and may attain $|m_0|=1$ on a set of positive measure.

For $|v|<1$, define
\begin{equation}\label{equ:inverse Langevin field}
	G(v)
	:=
	\frac{\mathcal L^{-1}(|v|)}{|v|}v
	=
	\frac{v}{q(|v|)},
	\qquad
	G(0):=0.
\end{equation}
Since $G(v)$ is parallel to $v$, one has $v\times G(v)=0$ and $\mathbb M(v)G(v)=2v$.
Thus, the closed-ball formulation \eqref{equ:closed ball LLB} can be equivalently written, whenever $|m|<1$, as
\begin{align}
	\label{equ:thermodynamic LLB}
	\partial_t m
	&=
	\gamma m\times h(m)+\mathbb M(m)h(m),
	\\
	h(m)
	&\coloneqq
	H_{\mathrm{MFA}}(m)-\theta G(m),
	\label{equ:thermodynamic effective field}
\end{align}
where $h(m)$ is called the thermodynamic effective field.

These two formulations play complementary roles. On the open unit ball, they are algebraically equivalent, and this equivalence remains valid at the level of weak solutions. The closed-ball formulation \eqref{equ:closed ball LLB} extends continuously to $|m|=1$ and is
therefore used to construct solutions from possibly saturated initial data. Once positive-time separation from saturation has been established, the same weak solution satisfies \eqref{equ:thermodynamic LLB} on every interval $[\tau,T]$, $\tau>0$, where the thermodynamic field is well-defined. This formulation is more natural for the energy and long-time
analysis.

\subsection{Energy structure}

Define the inverse-Langevin potential
\begin{equation}
	\Phi(s)
	:=
	\int_0^s\mathcal L^{-1}(\sigma)\,\mathrm d\sigma,
	\qquad
	0\leq s<1.
	\label{equ:inverse Langevin potential}
\end{equation}
Then, we have a variational derivative identity $\nabla_v\Phi(|v|)=G(v)$. Note that $\Phi(s)\to +\infty$ as $s\uparrow 1$.
Recalling $\theta$ is the rescaled temperature, we define the corresponding finite-temperature energy
\begin{equation}
	\mathcal E_\theta(m)
	:=
	\frac{\alpha}{2}
	\norm{\nabla m}{L^2(\mathscr D)}^2
	-
	\frac12
	\int_{\mathscr D}m\cdot\mathsf Bm\,\mathrm dx
	+
	\frac12
	\norm{H_{\mathrm d}[m]}{L^2(\mathbb R^3)}^2
	+
	\theta
	\int_{\mathscr D}\Phi(|m|)\,\mathrm dx.
	\label{equ:full nonlinear LLB energy}
\end{equation}
The negative variational derivative of the energy functional gives the thermodynamic effective field:
\[
-\frac{\delta\mathcal E_\theta}{\delta m}
=
\alpha\Delta m+\mathsf Bm+H_{\mathrm d}[m]
-\theta G(m)
=
h(m).
\]
Therefore, sufficiently regular solutions of~\eqref{equ:thermodynamic LLB} formally satisfy
\begin{equation}\label{equ:full nonlinear LLB energy identity}
	\frac{\mathrm d}{\mathrm dt}\mathcal E_\theta(m(t))
	+
	\int_{\mathscr D}
	h(m)\cdot\mathbb M(m)h(m)\,\mathrm dx
	=
	0,
\end{equation}
where the precessional term makes no contribution since $(m\times h)\cdot h=0$.

The singularity of $\Phi$ at $|m|=1$ is important for the analysis. If $|m_0|=1$ on a set of positive measure, then $\mathcal E_\theta(m_0)=+\infty$, even though the closed-ball formulation
\eqref{equ:closed ball LLB} remains well-defined. The existence argument will thus be based initially on the closed-ball equation \eqref{equ:closed ball LLB}. We shall prove that the solution enters the open unit ball instantaneously, after which the thermodynamic field, the finite-temperature energy, and the dissipation identity above are well-defined for every positive time.


\subsection{Relation with the standard reduced LLB equations}
\label{subsec:relation reduced LLB}

It is useful to clarify how the full nonlinear mean-field equation
considered here is related to the reduced LLB equations more commonly
used in the mathematical literature. To this end, we recall that $\mathsf B= JI+\mathsf A$, where $J>0$ denotes the homogeneous molecular exchange and
$\mathsf A$ contains the remaining anisotropy contribution. We also recall the expressions for the micromagnetic effective field $H_{\mathrm{mic}}(m)$ and the thermodynamic effective field $h(m)$ from \eqref{equ:H MFA mic} and \eqref{equ:thermodynamic effective field}, respectively.

The reduced LLB equations arise by expanding this thermodynamic field
and the mobility about the appropriate homogeneous equilibrium.
The discussion below is an asymptotic derivation of the usual reduced
models, rather than a rigorous convergence theorem from the full
equation. Rigorous limiting results, when available, are stated and proved
separately later in Section~\ref{sec:low temperature and long time}.

\subsubsection{Reduction below the bulk mean-field Curie temperature}

We first explain how the standard reduced LLB equation below the Curie
temperature is recovered from the full nonlinear model. Recall the thermodynamic field in \eqref{equ:thermodynamic effective field} can be
written as
\begin{equation}
	h(m)
	=
	H_{\mathrm{mic}}(m)
	+
	\left(
	J-\theta\frac{\mathcal L^{-1}(|m|)}{|m|}
	\right)m.
	\label{equ:thermodynamic field bulk decomposition}
\end{equation}

The reference magnitude for the below-Curie reduction is determined by
the corresponding \emph{homogeneous bulk mean-field problem}, obtained
from \eqref{equ:thermodynamic field bulk decomposition} by neglecting
$H_{\mathrm{mic}}$. 
The homogeneous dynamics then simplifies to
\begin{equation}\label{equ:homog dynamics}
\partial_t m
=
2\bigl(Jq(|m|)-\theta\bigr)m.
\end{equation}
A nonzero homogeneous equilibrium of magnitude
$m_{\mathrm e}\in(0,1)$ then satisfies $q(m_{\mathrm e})= \theta/J$.
Since $q$ decreases from $q(0)=1/3$ to $q(1)=0$, such an equilibrium exists
precisely when $\theta/J\in (0, 1/3)$. Thus, $\theta_{\mathrm c}^{\mathrm{bulk}}
\coloneq J/3$ is the Curie threshold of the \emph{homogeneous bulk mean-field problem}. This threshold arises here in the asymptotic reduction to the
usual model and should be distinguished from the rigorous spectral
threshold $\Lambda_\ast/3$ introduced later in
\eqref{equ:critical spectral value} for the full micromagnetic PDE. In the homogeneous isotropic setting where
$\Lambda_\ast=J$, the two thresholds coincide.

We now return to the full thermodynamic field \eqref{equ:thermodynamic field bulk decomposition}, and consider states for which $\theta<J/3$ and $|m|$ remains close to $m_{\mathrm e}\in (0,1)$. Introduce the longitudinal susceptibility
\begin{equation}
	\chi_\parallel
	\coloneq
	\left[
	\theta(\mathcal L^{-1})'(m_{\mathrm e})-J
	\right]^{-1}>0.
	\label{equ:longitudinal susceptibility}
\end{equation}
Consider the scalar longitudinal coefficient
\[
f(r)
\coloneq
J-\theta\frac{\mathcal L^{-1}(r)}{r}.
\]
Since $f(m_{\mathrm e})=0$, a Taylor expansion about
$r=m_{\mathrm e}$, together with
$q(m_{\mathrm e})=\theta/J$, gives
\[
f(r)
=
\frac{1}{2\chi_\parallel}
\left(
1-\frac{r^2}{m_{\mathrm e}^2}
\right)
+
O\bigl((r-m_{\mathrm e})^2\bigr).
\]
Consequently, for $|m|$ close to $m_{\mathrm e}$,
\begin{equation}
	h(m)
	=
	h_{\mathrm{red}}(m)
	+
	O\bigl((|m|-m_{\mathrm e})^2\bigr)m,
	\label{equ:below Curie field expansion}
\end{equation}
where
\begin{equation}\label{equ:h red}
h_{\mathrm{red}}(m)
\coloneq
H_{\mathrm{mic}}(m)
+
\frac{1}{2\chi_\parallel}
\left(
1-\frac{|m|^2}{m_{\mathrm e}^2}
\right)m.
\end{equation}
This is the first approximation made in the full nonlinear problem, leading to the usual effective field for this temperature regime considered in the mathematical literature~\cite{HamHam20, Soe25b,Soe25a, SoeTra23}.

The second approximation consists in freezing the longitudinal and
transverse mobilities at the equilibrium magnitude $m_{\mathrm e}$. Recall the mobility matrix $\mathbb M(m)$ defined in \eqref{equ:nonlinear LLB mobility}.
Since $q(m_{\mathrm e})=\theta/J$, this approximation at $m_{\mathrm e}$ gives
\[
\mathbb M(m)
=
\beta_\perp \left(I-\frac{m\otimes m}{|m|^2}\right)
+
\beta_\parallel \frac{m\otimes m}{|m|^2}
+
O\bigl(\bigl| |m|-m_{\mathrm e} \bigr|\bigr),
\]
where $\beta_\perp \coloneq 1-\theta/J$ and $\beta_\parallel \coloneq 2\theta/J$.

Substituting these two near-equilibrium approximations into the full
thermodynamic equation \eqref{equ:thermodynamic LLB}
and neglecting the higher-order terms yields
\begin{align}
	\partial_t m
	=
	\gamma m\times h_{\mathrm{red}}(m)
	+
	\frac{\beta_\parallel}{|m|^2}
	\bigl(m\cdot h_{\mathrm{red}}(m)\bigr)m
	-
	\frac{\beta_\perp}{|m|^2}
	m\times
	\bigl(
	m\times h_{\mathrm{red}}(m)
	\bigr).
	\label{equ:standard reduced below Curie LLB}
\end{align}
This is the standard reduced LLB equation below the bulk mean-field
Curie temperature~\cite{HamHam20, Soe25b}.

\subsubsection{Reduction above the bulk mean-field Curie temperature.}

For $\theta> J/3$, it is clear from \eqref{equ:homog dynamics} that the homogeneous equilibrium is $m=0$. We recall the expansions
\[
q(r)
=
\frac13-\frac15r^2+O(r^4),
\qquad
\mathbb M(m)
=
\frac23I+O(|m|^2),
\]
together with the exact identity
\[
\mathbb M(m)
\bigl(
Jm-\theta G(m)
\bigr)
=
\bigl(
2Jq(|m|)-2\theta
\bigr)m.
\]
Substituting these into the equation \eqref{equ:thermodynamic LLB}, and noting that $m\times (Jm-\theta G(m))=0$, we obtain
\begin{align}
	\partial_t m
	={}&
	\gamma m\times H_{\mathrm{mic}}(m)
	+
	\frac23H_{\mathrm{mic}}(m)
	-
	2\left(
	\theta-\frac J3
	\right)m
	-
	\frac{2J}{5}|m|^2m
	+
	\mathcal R(m),
	\label{equ:above Curie expansion}
\end{align}
where, formally near $m=0$,
\[
|\mathcal R(m)|
\leq
C\left(
|m|^2|H_{\mathrm{mic}}(m)|
+
|m|^5
\right).
\]
If anisotropy and stray-field contributions are omitted for simplicity,
then \eqref{equ:above Curie expansion} reduces, to leading order, to
\begin{equation}
	\partial_t m
	=
	\widetilde\gamma
	m\times\Delta m
	+
	\kappa_1\Delta m
	-
	\kappa_2
	\bigl(
	1+\mu_{\mathrm C}|m|^2
	\bigr)m,
	\label{equ:standard above Curie LLB}
\end{equation}
with
\[
\widetilde\gamma=\alpha\gamma,
\qquad
\kappa_1
=
\frac{2\alpha}{3},
\qquad
\kappa_2
=
2\left(
\theta-\frac J3
\right),
\qquad
\mu_{\mathrm C}
=
\frac{J}{
	5\left(\theta-J/3\right)
}.
\]
Up to the normalisation of the material coefficients,
\eqref{equ:standard above Curie LLB} is the simplified LLB equation
commonly analysed in the regime above the Curie temperature~\cite{BrzGolLe20, Le16,LiGuoLiuLiu21, QiuHuaWan26,Soe26}.

\begin{remark}[Low-temperature consistency]
	Consider the reduced below-Curie equation
	\eqref{equ:standard reduced below Curie LLB}. This equation was
	itself obtained by a near-equilibrium approximation of the full
	nonlinear model, so the following low-temperature argument is only
	formal. Since $q(m_{\mathrm e})=\theta/J$ and $q$ is strictly
	decreasing with $q(1)=0$, we have $m_{\mathrm e}\to1$ as
	$\theta/J\to0$. Moreover,
	\[
	\beta_\parallel
	=
	\frac{2\theta}{J}
	\to0,
	\qquad
	\beta_\perp
	=
	1-\frac{\theta}{J}
	\to1.
	\]
	Noting \eqref{equ:h red}, the longitudinal part of
	\eqref{equ:standard reduced below Curie LLB} is
	\[
	\frac{\beta_\parallel}{|m|^2}
	\bigl(m\cdot H_{\mathrm{mic}}(m)\bigr)m
	+
	\frac{\beta_\parallel}{2\chi_\parallel}
	\left(
	1-\frac{|m|^2}{m_{\mathrm e}^2}
	\right)m.
	\]
	Using
	$\mathcal L^{-1}(s)=s/q(s)$ and
	$q(m_{\mathrm e})=\theta/J$, while noting the definition of $\chi_\parallel$ in
	\eqref{equ:longitudinal susceptibility}, we find
	\[
	\frac{\beta_\parallel}{2\chi_\parallel}
	=
	-Jm_{\mathrm e}q'(m_{\mathrm e})
	\to J.
	\]
	Hence, formally for near-equilibrium families satisfying
	$|m|-m_{\mathrm e}\to0$, both longitudinal terms vanish as $\theta\to0$. Since the radial part of $h_{\mathrm{red}}(m)$ does not contribute to the cross-product terms, the remaining equation is
	\[
	\partial_t m
	=
	\gamma m\times H_{\mathrm{mic}}(m)
	-
	m\times\bigl(m\times H_{\mathrm{mic}}(m)\bigr),
	\qquad |m|=1,
	\]
	which is the Landau--Lifshitz--Gilbert (LLG) equation.
	Thus, the standard reduced model is formally consistent with the
	low-temperature Landau--Lifshitz limit. A direct rigorous passage from
	the full nonlinear LLB equation to the LLG equation is
	proved independently in Section~\ref{subsec:low temperature limit}.
\end{remark}


\section{Global existence and instantaneous separation}
\label{sec:global existence}

In this section, we establish the existence of a global weak solution to
\eqref{equ:thermodynamic LLB}--\eqref{equ:thermodynamic effective field}.
Since $G(m)$ is not defined at $|m|=1$, the construction is first carried out using the equivalent closed-ball formulation \eqref{equ:closed ball LLB}. We then prove that the solution is strictly separated from the saturation boundary $|m|=1$ for every positive time. This recovers~\eqref{equ:thermodynamic LLB} and its energy structure on every interval bounded away from $t=0$.

\subsection{Weak solutions and the main theorem}

The notion of weak solution employed in this paper is as follows.

\begin{definition}\label{def:global weak solution}
	Let $m_0\in H^1(\mathscr D)$ such that $|m_0|\leq 1$ a.e. in $\mathscr D$, and let $T>0$.
	A function $m$ is called a \emph{global weak solution} if
	\begin{enumerate}[(i)]
		\item $m
		\in
		H^1\bigl(0,T;L^2(\mathscr D)\bigr)
		\cap
		L^\infty\bigl(0,T;H^1(\mathscr D)\bigr)$;
		\item $|m|\leq 1$ a.e. in $(0,T)\times\mathscr D$;
		\item $m(0)=m_0$ in $L^2(\mathscr{D})$;
		\item for every $\tau \in (0,T)$, $m\in L^2\bigl(\tau,T;H^2(\mathscr D)\bigr)$ and $h(m)\in L^2\bigl(\tau,T;L^2(\mathscr D)\bigr)$; and
		\item the following equations hold for almost every $t\in (\tau,T)$:
		\begin{align}
			\label{equ:weak evolution equation}
			\inpro{\partial_t m(t)}{\varphi}
			&=
			\gamma \inpro{m(t)\times h\bigl(m(t)\bigr)}{\varphi}
			+
			\inpro{\mathbb M\bigl(m(t)\bigr) h\bigl(m(t)\bigr)}{\varphi},
			\\
			\label{equ:weak effective field}
			\inpro{h\bigl(m(t)\bigr)}{\psi}
			&=
			- \alpha \inpro{\nabla m(t)}{\nabla\psi}
			+
			\inpro{\mathsf Bm(t) + H_{\mathrm d}[m(t)]}{\psi}
			-
			\theta \inpro{G\bigl(m(t)\bigr)}{\psi},
		\end{align}
		for every $\varphi\in L^2(\mathscr D)$ and $\psi\in H^1(\mathscr D)$.
	\end{enumerate}
\end{definition}

We have the following main theorem of this section.

\begin{theorem}
	\label{thm:global existence full nonlinear LLB}
	Let $m_0\in H^1(\mathscr D)$ such that $|m_0|\leq 1$ a.e. in $\mathscr D$.
	Then there exists a global weak solution in the sense of
	Definition~\ref{def:global weak solution}. This solution enjoys the following additional properties:
	\begin{enumerate}[(i)]
		\item for every $p\geq1$, the solution satisfies
		\begin{equation}
			\int_{\mathscr D}
			\frac{1}
			{\bigl(1-|m(t,x)|^2\bigr)^p}\dx
			\leq
			C_p\bigl(1+t^{-p}\bigr)
			\qquad
			\text{for every }t>0;
			\label{equ:inverse distance final}
		\end{equation}
		\item for every $\tau>0$, there exists $\rho_\tau\in(0,1)$ such that
		\begin{equation}
			\operatorname*{ess\,sup}_{x\in\mathscr D}
			|m(t,x)|
			\leq
			\rho_\tau
			\qquad
			\text{for every }t\in[\tau,\infty);
			\label{equ:positive time separation final}
		\end{equation}
		\item for every $T>0$, $m \in C\bigl([0,T];H^1(\mathscr D)\bigr)$. In particular, \begin{equation}\label{equ:m0 converge H1}
		m(t)\to m_0 \qquad \text{strongly in }H^1(\mathscr D) \quad\text{as }t\downarrow0; 
		\end{equation}
		\item the energy identity
		\begin{equation}
			\mathcal E_\theta(m(t))
			+
			\int_s^t
			\int_{\mathscr D}
			h(m)\cdot\mathbb M(m)h(m)\,\mathrm dx\,\mathrm dr
			=
			\mathcal E_\theta(m(s))
			\label{equ:positive time energy identity}
		\end{equation}
		holds for every $0<s<t$. Furthermore, if $\bigl|\{x\in\mathscr D:|m_0(x)|=1\}\bigr|>0$, then
		\begin{equation}
			\mathcal E_\theta(m(t))
			\to +\infty
			\qquad
			\text{as }t\downarrow0.
			\label{equ:initial energy blowup}
		\end{equation}
	\end{enumerate}
\end{theorem}

\begin{remark}
	The solution is global and weak in the analytical sense, but no uniqueness is
	asserted. Nevertheless, on every interval $[\tau,T]$, $\tau>0$, it
	has the strong parabolic regularity
	\[
	m\in
	H^1(\tau,T;L^2(\mathscr D))
	\cap
	L^2(\tau,T;H^2(\mathscr D)).
	\]
\end{remark}

\subsection{Viscous Faedo--Galerkin approximation}

This section introduces a viscous regularisation of the problem and its Faedo--Galerkin approximation.
Let
\begin{equation}
	\pi(v)
	:=
	\begin{cases}
		v,&|v|\leq1,
		\\[1mm]
		\dfrac{v}{|v|},&|v|>1,
	\end{cases}
	\label{equ:radial projection existence}
\end{equation}
be the metric projection onto $\overline{B_1(0)}$, and set
\[
\widehat{\mathbb M}(v)
:=
\mathbb M \bigl(\pi(v)\bigr).
\]
Both maps $\pi$ and $\widehat{\mathbb M}$ are globally Lipschitz. Notice also that
\begin{equation}
	\widehat{\mathbb M}(v)v
	=
	\begin{cases}
		2q(|v|)v,&|v|\leq1,\\
		0,&|v|>1.
	\end{cases}
	\label{equ:extended mobility radial identity}
\end{equation}

For $\varepsilon\in(0,1]$, we consider a regularisation of \eqref{equ:closed ball LLB}, namely
\begin{align}
	\partial_t m_\varepsilon
	={}&
	\gamma\pi(m_\varepsilon)
	\times H_{\mathrm{MFA}}(m_\varepsilon)
	+
	\widehat{\mathbb M}(m_\varepsilon)
	H_{\mathrm{MFA}}(m_\varepsilon)
	-
	2\theta m_\varepsilon
	+
	\varepsilon\Delta m_\varepsilon,
	\label{equ:viscous closed ball problem}
\end{align}
with $\partial_\nu m_\varepsilon=0$ on $\partial\mathscr{D}$ and $m_\varepsilon(0)= m_0$.

Let $\{e_j\}_{j\geq1}$ be an $L^2(\mathscr D)$-orthonormal basis of
Neumann Laplacian eigenfunctions:
\[
-\Delta e_j=\lambda_j e_j \quad \text{in }\mathscr{D},
\qquad
\partial_\nu e_j=0 \quad \text{on }\partial\mathscr{D},
\]
and set $V_n:=\operatorname{span}\{e_1,\ldots,e_n\}$.
and let $P_n$ denote the $L^2$-orthogonal projection onto $V_n$. We seek
$m_{\varepsilon,n}\colon[0,\infty)\to V_n$ such that
\begin{align}
	\inpro{\partial_t m_{\varepsilon,n}}{\varphi_n}
	&=
	\gamma
	\inpro{
	\pi(m_{\varepsilon,n})
	\times H_{\mathrm{MFA}}(m_{\varepsilon,n})}{\varphi_n}
	+
	\inpro{
	\widehat{\mathbb M}(m_{\varepsilon,n})
	H_{\mathrm{MFA}}(m_{\varepsilon,n})}
	{\varphi_n}
	\nonumber\\
	&\quad
	-2\theta \inpro{m_{\varepsilon,n}}{\varphi_n}
	-
	\varepsilon
	\inpro{\nabla m_{\varepsilon,n}}{\nabla\varphi_n}
	\label{equ:viscous Galerkin}
\end{align}
for all $\varphi_n\in V_n$, with $m_{\varepsilon,n}(0)=P_n m_0$. Equivalently, we solve
\begin{align}
	\label{equ:projected viscous Galerkin}
	\partial_t m_{\varepsilon,n}
	=
	\gamma P_n \bigl(\pi(m_{\varepsilon,n})
	\times H_{\mathrm{MFA}}(m_{\varepsilon,n}) \bigr) 
	+
	P_n \left(\widehat{\mathbb M}(m_{\varepsilon,n})
	H_{\mathrm{MFA}}(m_{\varepsilon,n})\right) 
	-
	2\theta m_{\varepsilon,n}
	+
	\varepsilon \Delta m_{\varepsilon,n}.
\end{align}
Here, $H_{\mathrm{MFA}}$ is given by \eqref{equ:H MFA mic}.

We begin deriving the necessary uniform estimates in the following proposition.

\begin{proposition}
	\label{prop:viscous Galerkin estimates}
	For every $\varepsilon\in(0,1]$ and every $n\in\mathbb N$, problem~\eqref{equ:viscous Galerkin} has a unique global solution. For every $T>0$,
	\begin{align}
		&\sup_{0\leq t\leq T}
		\norm{m_{\varepsilon,n}(t)}{H^1(\mathscr D)}^2
		+
		\int_0^T
		\int_{\mathscr D}
		\Delta m_{\varepsilon,n}\cdot
		\widehat{\mathbb M}(m_{\varepsilon,n})
		\Delta m_{\varepsilon,n}\,\mathrm dx\,\mathrm dt
		\nonumber\\
		&\qquad
		+
		\varepsilon
		\int_0^T
		\norm{\Delta m_{\varepsilon,n}}{L^2(\mathscr D)}^2\,\mathrm dt
		+
		\int_0^T
		\norm{\partial_t m_{\varepsilon,n}}{L^2(\mathscr D)}^2\,\mathrm dt
		\leq
		C_T,
		\label{equ:viscous Galerkin uniform}
	\end{align}
	where $C_T$ is independent of $n$ and $\varepsilon$.
\end{proposition}

\begin{proof}
	Since all coefficients in the finite-dimensional system are locally Lipschitz, the Cauchy--Lipschitz theorem gives a unique local solution.
	
	Taking $\varphi_n=m_{\varepsilon,n}$ in~\eqref{equ:viscous Galerkin}, noting that $\pi(v)$ is parallel to $v$, gives
	\begin{align}
		\frac12\frac{\mathrm d}{\mathrm dt}
		\norm{m_{\varepsilon,n}}{L^2}^2
		+
		2\theta\norm{m_{\varepsilon,n}}{L^2}^2
		+
		\varepsilon \norm{\nabla m_{\varepsilon,n}}{L^2}^2
		=
		\inpro{H_{\mathrm{MFA}}(m_{\varepsilon,n})}{\widehat{\mathbb M}(m_{\varepsilon,n}) m_{\varepsilon,n}}.
		\label{equ:Galerkin L2 identity}
	\end{align}
	The map $v\mapsto \widehat{\mathbb M}(v)v$ is globally Lipschitz and vanishes for $|v|>1$ by \eqref{equ:extended mobility radial identity}. By the Sobolev chain rule, $\bigl|\nabla \widehat{\mathbb M}(m_{\varepsilon,n})
	m_{\varepsilon,n}\bigr| \leq C |\nabla m_{\varepsilon,n}|$ a.e. in $\mathscr{D}$.  Hence, after integrating by parts,
	\[
	\left|
	\bigl(
	\Delta m_{\varepsilon,n},
	\widehat{\mathbb M}(m_{\varepsilon,n})
	m_{\varepsilon,n}
	\bigr)_{\mathscr D}
	\right|
	\leq
	C
	\norm{\nabla m_{\varepsilon,n}}{L^2}^2.
	\]
	Furthermore, by \eqref{equ:stray field W2r},
	\[
	\norm{
		\mathsf Bm_{\varepsilon,n}
		+
		H_{\mathrm d}[m_{\varepsilon,n}]
	}{L^2}
	\leq
	C\norm{m_{\varepsilon,n}}{L^2},
	\]
	Therefore, continuing from \eqref{equ:Galerkin L2 identity}, we obtain
	\begin{equation}
		\frac{\mathrm d}{\mathrm dt}
		\norm{m_{\varepsilon,n}}{L^2}^2
		+
		4\theta\norm{m_{\varepsilon,n}}{L^2}^2
		+
		2\varepsilon \norm{\nabla m_{\varepsilon,n}}{L^2}^2
		\leq
		C
		\norm{m_{\varepsilon,n}}{H^1}^2.
		\label{equ:Galerkin L2 estimate}
	\end{equation}
	
	We next test \eqref{equ:viscous Galerkin} with
	$-\Delta m_{\varepsilon,n}$. In this case, the highest-order precessional term
	vanishes.
	Using
	\eqref{equ:mobility square estimate} and
	\eqref{equ:precession mobility control}, we obtain
	\begin{align}
		&\frac12\frac{\mathrm d}{\mathrm dt}
		\norm{\nabla m_{\varepsilon,n}}{L^2}^2
		+
		\frac{\alpha}{2}
		\int_{\mathscr D}
		\Delta m_{\varepsilon,n}\cdot
		\widehat{\mathbb M}(m_{\varepsilon,n})
		\Delta m_{\varepsilon,n}\,\mathrm dx
		+
		2\theta
		\norm{\nabla m_{\varepsilon,n}}{L^2}^2
		+
		\varepsilon
		\norm{\Delta m_{\varepsilon,n}}{L^2}^2
		\nonumber\\
		&\leq
		C
		\norm{
			\mathsf Bm_{\varepsilon,n}
			+
			H_{\mathrm d}[m_{\varepsilon,n}]
		}{L^2}^2
		\nonumber\\
		&\leq
		C\norm{m_{\varepsilon,n}}{L^2}^2.
		\label{equ:Galerkin H1 estimate}
	\end{align}
	Combining
	\eqref{equ:Galerkin L2 estimate} and
	\eqref{equ:Galerkin H1 estimate}, we obtain a bound for the first three terms in \eqref{equ:viscous Galerkin uniform} by the Gronwall lemma. In particular,
	the finite-dimensional solution cannot blow up and is therefore
	global.
	
	Finally, the Galerkin-projected equation \eqref{equ:projected viscous Galerkin}, together with \eqref{equ:mobility square estimate} and
	\eqref{equ:precession mobility control}, gives
	\begin{align}\label{equ:Galerkin time der estimate}
		\norm{\partial_t m_{\varepsilon,n}}{L^2}^2
		\leq
		C
		\int_{\mathscr D}
		\Delta m_{\varepsilon,n}\cdot
		\widehat{\mathbb M}(m_{\varepsilon,n})
		\Delta m_{\varepsilon,n}\,\mathrm dx
		+
		C\norm{m_{\varepsilon,n}}{L^2}^2
		+
		C\varepsilon^2
		\norm{\Delta m_{\varepsilon,n}}{L^2}^2.
	\end{align}
	Integration in time and using the bound for the first three terms in \eqref{equ:viscous Galerkin uniform} derived earlier then completes the proof.
\end{proof}

For fixed $\varepsilon>0$, the preceding estimates imply, after extraction
of a subsequence,
\[
\begin{aligned}
	m_{\varepsilon,n}
	&\overset{*}{\rightharpoonup} m_\varepsilon
	&&\text{in } L^\infty\bigl(0,T;H^1(\mathscr D)\bigr),
	\\
	m_{\varepsilon,n}
	&\rightharpoonup m_\varepsilon
	&&\text{in } L^2\bigl(0,T;H^2(\mathscr D)\bigr),
	\\
	\partial_t m_{\varepsilon,n}
	&\rightharpoonup \partial_t m_\varepsilon
	&&\text{in } L^2\bigl(0,T;L^2(\mathscr D)\bigr).
\end{aligned}
\]
Aubin--Lions compactness theorem gives
\begin{equation}
	m_{\varepsilon,n}
	\to
	m_\varepsilon
	\qquad
	\text{strongly in }L^2\bigl(0,T;H^1(\mathscr D)\bigr)
	\cap C\bigl([0,T];L^2(\mathscr D)\bigr).
	\label{equ:Galerkin strong convergence}
\end{equation}
The boundedness and continuity of $\pi$ and
$\widehat{\mathbb M}$, together with the linear continuity of the stray
field, permit passage to the limit of
\eqref{equ:viscous Galerkin} by a standard argument. Hence, for fixed $\varepsilon>0$,
\begin{equation}\label{equ:regularity viscous Galerkin}
m_\varepsilon\in H^1\bigl(0,T;L^2(\mathscr{D})\bigr) \cap L^\infty\bigl(0,T;H^1(\mathscr{D})\bigr) \cap L^2\bigl(0,T;H^2(\mathscr{D})\bigr)
\end{equation}
is a strong solution of~\eqref{equ:viscous closed ball problem}, or equivalently \eqref{equ:viscous closed ball problem} holds in $L^2(\mathscr{D})$ for a.e. $t>0$.

\subsection{Further $\varepsilon$-uniform estimates}

We now derive further estimates which are uniform in $\varepsilon$. The following lemma shows invariance of the physical state space.

\begin{lemma}
	\label{lem:viscous invariant ball}
	For every $\varepsilon\in(0,1]$,
	\begin{equation}
		|m_\varepsilon(t,x)|
		\leq 1
		\qquad
		\text{for a.e. }x\in\mathscr D
		\text{ and every }t\geq 0.
		\label{equ:bound ae viscous}
	\end{equation}
\end{lemma}

\begin{proof}
	Define
	\[
	\widehat q(s)
	:=
	\begin{cases}
		q(s),&0\leq s\leq1,\\
		0,&s>1.
	\end{cases}
	\]
	Using \eqref{equ:extended mobility radial identity} and taking the
	scalar product of \eqref{equ:viscous closed ball problem} with
	$2m_\varepsilon$, we obtain
	\begin{align}
		\partial_t|m_\varepsilon|^2
		={}&
		\bigl(
		2\alpha\widehat q(|m_\varepsilon|)
		+\varepsilon
		\bigr)
		\Delta |m_\varepsilon|^2
		-
		2
		\bigl(
		2\alpha\widehat q(|m_\varepsilon|)
		+\varepsilon
		\bigr)
		|\nabla m_\varepsilon|^2
		\nonumber\\
		&+
		4\widehat q(|m_\varepsilon|)
		m_\varepsilon\cdot
		\bigl(
		\mathsf Bm_\varepsilon
		+
		H_{\mathrm d}[m_\varepsilon]
		\bigr)
		-
		4\theta |m_\varepsilon|^2.
		\label{equ:viscous scalar equation}
	\end{align}
	Here, we have used the fact that $2m_\varepsilon\cdot\Delta m_\varepsilon
	=
	\Delta|m_\varepsilon|^2
	-
	2|\nabla m_\varepsilon|^2$.
	
	Now, set
	\[
	w_\varepsilon
	:=
	\bigl(|m_\varepsilon|^2-1\bigr)_+.
	\]
	We first justify that $w_\varepsilon$ is an admissible test function.
	For fixed $\varepsilon>0$, the regularity \eqref{equ:regularity viscous Galerkin} implies
	$m_\varepsilon
	\in
	L^2(0,T;L^\infty(\mathscr D))$ since $H^2(\mathscr D)\hookrightarrow L^\infty(\mathscr D)$.
	Consequently,
	\begin{align*}
	\partial_t|m_\varepsilon|^2
	&=
	2m_\varepsilon\cdot\partial_t m_\varepsilon
	\in
	L^1(0,T;L^2(\mathscr D)),
	\\
	\nabla|m_\varepsilon|^2
	&=
	2m_\varepsilon\cdot\nabla m_\varepsilon
	\in
	L^2(0,T;L^2(\mathscr D)).
	\end{align*}
	Hence,
	\[
	|m_\varepsilon|^2-1
	\in
	W^{1,1}(0,T;L^2(\mathscr D))
	\cap
	L^2(0,T;H^1(\mathscr D)).
	\]
	Since the map $s\mapsto s_+$ is Lipschitz, the Sobolev chain rule gives
	\[
	w_\varepsilon
	\in
	W^{1,1}(0,T;L^2(\mathscr D))
	\cap
	L^2(0,T;H^1(\mathscr D)),
	\]
	with
	\begin{align}
		\partial_t w_\varepsilon
		&=
		\mathbf 1_{\{|m_\varepsilon|^2>1\}}
		\partial_t|m_\varepsilon|^2,
		\qquad
		\nabla w_\varepsilon
		=
		\mathbf 1_{\{|m_\varepsilon|^2>1\}}
		\nabla|m_\varepsilon|^2
		\quad \text{a.e.}
		\label{equ:positive part space chain rule}
	\end{align}
	In particular,
	$t\mapsto\norm{w_\varepsilon(t)}{L^2}^2$ is absolutely continuous
	and
	\begin{equation}
		\frac12\frac{\mathrm d}{\mathrm dt}
		\norm{w_\varepsilon}{L^2}^2
		=
		\int_{\mathscr D}
		\partial_t|m_\varepsilon|^2
		w_\varepsilon\,\mathrm dx
		\qquad\text{for a.e. }t.
		\label{equ:positive part energy chain rule}
	\end{equation}
	Indeed, this is the chain rule applied to
	$s\mapsto \frac12(s_+)^2$, whose derivative is $s_+$.
	
	We may therefore test \eqref{equ:viscous scalar equation} by $w_\varepsilon$. On the set where $w_\varepsilon>0$, one has $|m_\varepsilon|>1$, and thus $\widehat q(|m_\varepsilon|)=0$.
	Therefore, all terms multiplied by $\widehat q(|m_\varepsilon|)$ vanish after multiplication by $w_\varepsilon$. Moreover, using
	\eqref{equ:positive part space chain rule}, we have
	\[
	-\varepsilon
	\int_{\mathscr D}
	\nabla|m_\varepsilon|^2
	\cdot\nabla w_\varepsilon\,\mathrm dx
	=
	-\varepsilon
	\norm{\nabla w_\varepsilon}{L^2}^2,
	\]
	where there is no boundary contribution since $\partial_\nu |m_\varepsilon|^2
	=
	2m_\varepsilon\cdot\partial_\nu m_\varepsilon
	=
	0$ on $\partial\mathscr D$.
	Hence, for almost every $t>0$, we obtain
	\begin{align}
		\frac12\frac{\mathrm d}{\mathrm dt}
		\norm{w_\varepsilon}{L^2(\mathscr D)}^2
		+
		\varepsilon
		\norm{\nabla w_\varepsilon}{L^2(\mathscr D)}^2
		+
		2\varepsilon
		\int_{\mathscr D}
		|\nabla m_\varepsilon|^2
		w_\varepsilon\,\mathrm dx
		+
		4\theta
		\int_{\mathscr D}
		|m_\varepsilon|^2
		w_\varepsilon\,\mathrm dx
		=0.
		\label{equ:positive part invariant ball estimate}
	\end{align}
	
	Since $|m_0|\leq1$ almost everywhere,
	\[
	w_\varepsilon(0)
	=
	\bigl(|m_0|^2-1\bigr)_+
	=0.
	\]
	Integrating \eqref{equ:positive part invariant ball estimate} over
	$(0,t)$ and dropping all non-negative terms on the left-hand side yields $\norm{w_\varepsilon(t)}{L^2(\mathscr D)}^2=0$ for every $t \geq 0$, which implies $|m_\varepsilon (t,x)|\leq 1$ for a.e. $(t,x)$.
	Moreover, since $m_\varepsilon\in C\bigl([0,T];L^2(\mathscr D)\bigr)$, the bound
	extends from a.e. $t$ to every $t\geq 0$ by approximation with a
	sequence of times for which \eqref{equ:bound ae viscous} holds and strong
	$L^2$-convergence.
	This proves the assertion.
\end{proof}

After Lemma~\ref{lem:viscous invariant ball} is established, the artificial extensions are
inactive, so that
\begin{equation}\label{equ:proj inactive}
\pi(m_\varepsilon)=m_\varepsilon,
\qquad
\widehat{\mathbb M}(m_\varepsilon)
=
\mathbb M(m_\varepsilon)
\qquad\text{for a.e. }x\in\mathscr D
\text{ and every }t\geq 0.
\end{equation}
The following uniform estimates will be used repeatedly in subsequent argument.

\begin{lemma}
	\label{lem:uniform viscous estimates}
	Let $m_\varepsilon$ be a solution of
	\eqref{equ:viscous closed ball problem}, with
	$\varepsilon\in(0,1]$. There exists a constant $C>0$ such that
	\begin{align}
		&\sup_{t\geq0}
		\norm{m_\varepsilon(t)}{H^1(\mathscr D)}^2
		+
		\sup_{t\geq0}
		\int_t^{t+1}
		\int_{\mathscr D}
		\Delta m_\varepsilon\cdot
		\mathbb M(m_\varepsilon)
		\Delta m_\varepsilon
		\,\mathrm dx\,\mathrm ds
		\nonumber\\
		&\qquad
		+
		\varepsilon
		\sup_{t\geq0}
		\int_t^{t+1}
		\norm{\Delta m_\varepsilon(s)}{L^2(\mathscr D)}^2
		\,\mathrm ds
		+
		\sup_{t\geq0}
		\int_t^{t+1}
		\norm{\partial_t m_\varepsilon(s)}{L^2(\mathscr D)}^2
		\,\mathrm ds
		\leq C.
		\label{equ:global viscous estimates}
	\end{align}
	The constant $C>0$ depends only on
	$\mathscr D$, $\alpha$, $\theta$, $\gamma$, $\mathsf B$, and $m_0$,
	but is independent of $\varepsilon$, $t$, and the final time $T$.
\end{lemma}

\begin{proof}
	By Lemma~\ref{lem:viscous invariant ball}, we have $\norm{m_\varepsilon(t)}{L^2(\mathscr D)}^2 \leq |\mathscr D|$ for all $t\geq 0$. Furthermore, we have the regularity \eqref{equ:regularity viscous Galerkin} so that the viscous equation \eqref{equ:viscous closed ball problem} holds in $L^2(\mathscr{D})$ for a.e. $t>0$. Thus, by the same argument leading to \eqref{equ:Galerkin H1 estimate}, but applied directly to \eqref{equ:viscous closed ball problem}, we obtain for a.e. $t>0$,
	\begin{align}
		\frac12\frac{\mathrm d}{\mathrm dt}
		\norm{\nabla m_\varepsilon}{L^2(\mathscr D)}^2
		+
		\frac{\alpha}{2}
		\int_{\mathscr D}
		\Delta m_\varepsilon\cdot
		\mathbb M(m_\varepsilon)
		\Delta m_\varepsilon
		\,\mathrm dx
		+
		2\theta
		\norm{\nabla m_\varepsilon}{L^2(\mathscr D)}^2
		+
		\varepsilon
		\norm{\Delta m_\varepsilon}{L^2(\mathscr D)}^2
		\leq C,
		\label{equ:uniform viscous differential inequality}
	\end{align}
	where $C$ is independent of $\varepsilon$.
	Discarding the non-negative second-order dissipation terms yields
	\[
	\frac{\mathrm d}{\mathrm dt}
	\norm{\nabla m_\varepsilon}{L^2(\mathscr D)}^2
	+
	4\theta
	\norm{\nabla m_\varepsilon}{L^2(\mathscr D)}^2
	\leq C.
	\]
	By the Gronwall lemma, we obtain
	\begin{equation}\label{equ:gronwall nabla me}
	\norm{\nabla m_\varepsilon(t)}{L^2(\mathscr D)}^2
	\leq
	e^{-4\theta t}
	\norm{\nabla m_0}{L^2(\mathscr D)}^2
	+
	C\bigl(1-e^{-4\theta t}\bigr),
	\end{equation}
	and hence
	\begin{equation}\label{equ:sup m eps H1}
	\sup_{t\geq0}
	\norm{m_\varepsilon(t)}{H^1(\mathscr D)}^2
	\leq C.
	\end{equation}
	Next, integrating \eqref{equ:uniform viscous differential inequality} over
	$(t,t+1)$ gives
	\begin{align*}
		&\frac12
		\norm{\nabla m_\varepsilon(t+1)}{L^2}^2
		-
		\frac12
		\norm{\nabla m_\varepsilon(t)}{L^2}^2
		+
		\frac{\alpha}{2}
		\int_t^{t+1}
		\int_{\mathscr D}
		\Delta m_\varepsilon\cdot
		\mathbb M(m_\varepsilon)
		\Delta m_\varepsilon
		\dx\ds
		\\
		&\quad+
		2\theta
		\int_t^{t+1}
		\norm{\nabla m_\varepsilon(s)}{L^2}^2
		\ds
		+
		\varepsilon
		\int_t^{t+1}
		\norm{\Delta m_\varepsilon(s)}{L^2}^2
		\ds
		\leq C.
	\end{align*}
	Using \eqref{equ:sup m eps H1} and dropping the non-negative terms gives
	\begin{align}
		&\sup_{t\geq0}
		\int_t^{t+1}
		\int_{\mathscr D}
		\Delta m_\varepsilon\cdot
		\mathbb M(m_\varepsilon)
		\Delta m_\varepsilon
		\dx \ds
		+
		\varepsilon \,
		\sup_{t\geq0}
		\int_t^{t+1}
		\norm{\Delta m_\varepsilon(s)}{L^2(\mathscr D)}^2
		\ds
		\leq C.
		\label{equ:uniform viscous intermediate dissipation}
	\end{align}
	
	By the same argument leading to 
	\eqref{equ:Galerkin time der estimate}, but using \eqref{equ:viscous closed ball problem} and noting that $|m_\varepsilon|\leq 1$, we obtain
	\begin{align*}
		\norm{\partial_t m_\varepsilon}{L^2(\mathscr D)}^2
		\leq{}&
		C
		\int_{\mathscr D}
		\Delta m_\varepsilon\cdot
		\mathbb M(m_\varepsilon)
		\Delta m_\varepsilon
		\,\mathrm dx
		+
		C
		+
		C\varepsilon^2
		\norm{\Delta m_\varepsilon}{L^2(\mathscr D)}^2.
	\end{align*}
	Integrating this inequality over $(t,t+1)$ and using
	\eqref{equ:uniform viscous intermediate dissipation}, together with
	$\varepsilon\leq1$, yields
	\[
	\sup_{t\geq0}
	\int_t^{t+1}
	\norm{\partial_t m_\varepsilon(s)}{L^2(\mathscr D)}^2
	\,\mathrm ds
	\leq C.
	\]
	Combining this with the preceding estimates \eqref{equ:sup m eps H1} and \eqref{equ:uniform viscous intermediate dissipation} proves \eqref{equ:global viscous estimates}.
\end{proof}

\subsection{Instantaneous inverse-distance and regularity estimates}

Next, we derive an entropy-type estimate controlling the approach of $m_\varepsilon$ to saturation. Since
\[
d_\varepsilon(s)
\coloneq
2\alpha q(\sqrt s)+\varepsilon
\simeq
\varepsilon+1-s,
\]
negative powers of $d_\varepsilon(|m_\varepsilon|^2)$ control the inverse regularised gap to $|m_\varepsilon|=1$. Testing the scalar equation for $|m_\varepsilon|^2$ by $d_\varepsilon(|m_\varepsilon|^2)^{-p-1}$ yields the following $\varepsilon$-uniform estimate, which will later be used to prove strict separation.

\begin{proposition}
	\label{prop:inverse distance viscous}
	For every $p\geq1$, there exists $C_p>0$, independent of
	$\varepsilon\in(0,1]$, such that
	\begin{equation}
		\int_{\mathscr D}
		\frac{\mathrm dx}
		{\bigl(
			\varepsilon+1-|m_\varepsilon(t,x)|^2
			\bigr)^p}
		\leq
		C_p\bigl(1+t^{-p}\bigr)
		\qquad
		\text{for a.e. }t>0.
		\label{equ:viscous inverse distance}
	\end{equation}
\end{proposition}

\begin{proof}
	Let $q$ be defined by \eqref{equ:q and L} and $p\geq 1$. For $s\in [0,1]$, set
	\[
	d_\varepsilon(s)
	:=
	2\alpha q(\sqrt s)+\varepsilon,
	\qquad
	\Psi_{p,\varepsilon}(s)
	:=
	\int_0^s
	d_\varepsilon(r)^{-p-1} \dr.
	\]
	By \eqref{equ:q saturation equivalence}, there exist constants $c,C>0$, independent of $\varepsilon\in(0,1]$, such that
	\begin{equation}
		c\bigl(\varepsilon+1-s\bigr)
		\leq
		d_\varepsilon(s)
		\leq
		C\bigl(\varepsilon+1-s\bigr)
		\qquad
		\text{for all }s\in[0,1].
		\label{equ:d epsilon equivalence}
	\end{equation}
	A direct integration gives
	\begin{equation}
		c\,d_\varepsilon(s)^{-p}-C
		\leq
		\Psi_{p,\varepsilon}(s)
		\leq
		C\,d_\varepsilon(s)^{-p}
		\qquad
		\text{for all }s\in[0,1],
		\label{equ:Psi comparison}
	\end{equation}
	with constants independent of $\varepsilon$.
	
	By Lemma~\ref{lem:viscous invariant ball},
	$|m_\varepsilon|\leq1$, and hence
	$\widehat q(|m_\varepsilon|)=q(|m_\varepsilon|)$ in
	\eqref{equ:viscous scalar equation}. We now multiply that equation by
	$d_\varepsilon(|m_\varepsilon|^2)^{-p-1} = \Psi_{p,\varepsilon}'(|m_\varepsilon|^2)$
	and integrate over $\mathscr D$.
	First, note that this nonlinear testing is legitimate. Indeed, for fixed $\varepsilon>0$, the regularity \eqref{equ:regularity viscous Galerkin} and \eqref{equ:bound ae viscous} imply
	\[
	|m_\varepsilon|^2
	\in
	H^1(0,T;L^2(\mathscr D))
	\cap
	L^\infty(0,T;H^1(\mathscr D)).
	\]
	Moreover, $d_\varepsilon\geq\varepsilon>0$, and
	$s\mapsto q(\sqrt{s})$ is continuously differentiable on $[0,1]$.
	Hence $\Psi_{p,\varepsilon}'$ is Lipschitz on $[0,1]$ for fixed
	$\varepsilon$. The temporal Sobolev chain rule thus gives
	\[
	\frac{\mathrm d}{\mathrm dt}
	\int_{\mathscr D}
	\Psi_{p,\varepsilon}(|m_\varepsilon|^2)\,\mathrm dx
	=
	\int_{\mathscr D}
	d_\varepsilon(|m_\varepsilon|^2)^{-p-1}
	\partial_t|m_\varepsilon|^2\,\mathrm dx
	\qquad \text{for a.e. } t>0.
	\]
	Consequently, using \eqref{equ:viscous scalar equation}, we obtain
	\begin{align}
		&\frac{\mathrm d}{\mathrm dt}
		\int_{\mathscr D}
		\Psi_{p,\varepsilon}(|m_\varepsilon|^2)\,\mathrm dx
		\nonumber\\
		&=
		\int_{\mathscr D}
		d_\varepsilon(|m_\varepsilon|^2)^{-p}
		\Delta |m_\varepsilon|^2\,\mathrm dx
		-
		2
		\int_{\mathscr D}
		d_\varepsilon(|m_\varepsilon|^2)^{-p}
		|\nabla m_\varepsilon|^2\,\mathrm dx
		\nonumber\\
		&\qquad
		+
		4
		\int_{\mathscr D}
		q(|m_\varepsilon|)
		d_\varepsilon(|m_\varepsilon|^2)^{-p-1}
		m_\varepsilon\cdot\mathsf Bm_\varepsilon\,\mathrm dx
		+
		4
		\int_{\mathscr D}
		q(|m_\varepsilon|)
		d_\varepsilon(|m_\varepsilon|^2)^{-p-1}
		m_\varepsilon\cdot H_{\mathrm d}[m_\varepsilon]\,\mathrm dx
		\nonumber\\
		&\qquad
		-
		4\theta
		\int_{\mathscr D}
		|m_\varepsilon|^2
		d_\varepsilon(|m_\varepsilon|^2)^{-p-1}\,\mathrm dx
		\nonumber\\
		&\quad=
		-
		\int_{\mathscr D}
		\left.
		\frac{\mathrm d}{\mathrm ds}
		\bigl[d_\varepsilon(s)^{-p}\bigr]
		\right|_{s=|m_\varepsilon|^2}
		\bigl|\nabla |m_\varepsilon|^2\bigr|^2\,\mathrm dx
		-
		2
		\int_{\mathscr D}
		d_\varepsilon(|m_\varepsilon|^2)^{-p}
		|\nabla m_\varepsilon|^2\,\mathrm dx
		\nonumber\\
		&\qquad
		+
		\underbrace{
			4
			\int_{\mathscr D}
			q(|m_\varepsilon|)
			d_\varepsilon(|m_\varepsilon|^2)^{-p-1}
			m_\varepsilon\cdot\mathsf Bm_\varepsilon\,\mathrm dx
		}_{\text{anisotropy term}}
		+
		\underbrace{
			4
			\int_{\mathscr D}
			q(|m_\varepsilon|)
			d_\varepsilon(|m_\varepsilon|^2)^{-p-1}
			m_\varepsilon\cdot H_{\mathrm d}[m_\varepsilon]\,\mathrm dx
		}_{\text{stray-field term}}
		\nonumber\\
		&\qquad
		-
		4\theta
		\int_{\mathscr D}
		|m_\varepsilon|^2
		d_\varepsilon(|m_\varepsilon|^2)^{-p-1}\,\mathrm dx.
		\label{equ:inverse distance weighted identity}
	\end{align}
	Here in the last equality, the first term on the right-hand side is obtained using integration by parts and the fact that $\partial_\nu|m_\varepsilon|^2= 0$ on $\partial\mathscr D$.
	Since $q$ is decreasing, $d_\varepsilon$ is decreasing, and hence
	\[
	\frac{\mathrm d}{\mathrm ds}
	\bigl[d_\varepsilon(s)^{-p}\bigr]
	=
	-pd_\varepsilon(s)^{-p-1}d_\varepsilon'(s)
	\geq0.
	\]
	Thus the first two terms and the last term on the right-hand side of
	\eqref{equ:inverse distance weighted identity} are nonpositive.
	It remains only to estimate the anisotropy and stray-field terms.
	
	For the anisotropy term, using the fact that $\mathsf B$ is symmetric and the definition of $d_\varepsilon$, we have
	\[
	m_\varepsilon\cdot\mathsf Bm_\varepsilon
	\leq
	\norm{\mathsf B}{}|m_\varepsilon|^2,
	\qquad
	q(|m_\varepsilon|)
	\leq
	\frac{1}{2\alpha}
	d_\varepsilon(|m_\varepsilon|^2).
	\]
	Consequently, for every $\eta>0$,
	\begin{align}
		4 \int_{\mathscr D}
		q(|m_\varepsilon|)
		d_\varepsilon(|m_\varepsilon|^2)^{-p-1}
		m_\varepsilon\cdot\mathsf Bm_\varepsilon\,\mathrm dx
		&\leq
		C
		\int_{\mathscr D}
		|m_\varepsilon|^2
		d_\varepsilon(|m_\varepsilon|^2)^{-p}\,\mathrm dx
		\nonumber\\
		&\leq
		\eta
		\int_{\mathscr D}
		|m_\varepsilon|^2
		d_\varepsilon(|m_\varepsilon|^2)^{-p-1}\,\mathrm dx
		+
		C_{\eta,p}.
		\label{equ:inverse distance anisotropy estimate}
	\end{align}
	Indeed, on the set $\{|m_\varepsilon|^2\leq1/2\}$, we have
	\[
	d_\varepsilon(|m_\varepsilon|^2)
	\geq
	2\alpha q(1/\sqrt2)>0,
	\]
	so the corresponding integrand is uniformly bounded. On
	$\{|m_\varepsilon|^2>1/2\}$, thanks to the elementary inequality $x^p\leq\eta x^{p+1}+C_{\eta,p}$ for any $x\geq 0$, applied to
	$x=d_\varepsilon(|m_\varepsilon|^2)^{-1}$, we obtain
	\eqref{equ:inverse distance anisotropy estimate}.
	
	For the stray-field term in \eqref{equ:inverse distance weighted identity}, we similarly infer that
	\begin{align}
		&4
		\int_{\mathscr D}
		q(|m_\varepsilon|)
		d_\varepsilon(|m_\varepsilon|^2)^{-p-1}
		m_\varepsilon\cdot H_{\mathrm d}[m_\varepsilon]\dx
		\nonumber\\
		&\qquad
		\leq
		C
		\int_{\mathscr D}
		d_\varepsilon(|m_\varepsilon|^2)^{-p}
		\bigl|H_{\mathrm d}[m_\varepsilon]\bigr| \dx
		\nonumber\\
		&\qquad
		\leq
		\eta
		\int_{\mathscr D}
		|m_\varepsilon|^2
		d_\varepsilon(|m_\varepsilon|^2)^{-p-1}\dx
		+
		C_{\eta,p}
		\left(
		1+
		\norm{H_{\mathrm d}[m_\varepsilon]}
		{L^{p+1}(\mathbb R^3)}^{p+1}
		\right)
		\nonumber\\
		&\qquad\leq
		\eta
		\int_{\mathscr D}
		|m_\varepsilon|^2
		d_\varepsilon(|m_\varepsilon|^2)^{-p-1}\,\mathrm dx
		+
		C_{\eta,p}.
		\label{equ:inverse distance stray estimate}
	\end{align}
	To justify the penultimate inequality, we see that on the set $\{|m_\varepsilon|^2\leq1/2\}$, the factor $d_\varepsilon(|m_\varepsilon|^2)^{-p}$ is uniformly bounded, so
	H\"older's inequality controls this part by $C\norm{H_{\mathrm d}[m_\varepsilon]}{L^{p+1}}$.
	On $\{|m_\varepsilon|^2>1/2\}$, Young's inequality with conjugate
	exponents $(p+1)/p$ and $p+1$ gives
	\[
	d_\varepsilon(|m_\varepsilon|^2)^{-p}
	\bigl|H_{\mathrm d}[m_\varepsilon]\bigr|
	\leq
	\widetilde{\eta}
	d_\varepsilon(|m_\varepsilon|^2)^{-p-1}
	+
	C_{\widetilde{\eta},p}
	|H_{\mathrm d}[m_\varepsilon]|^{p+1},
	\]
	and, since $|m_\varepsilon|^2>1/2$, the first term is bounded
	by $2\widetilde{\eta} |m_\varepsilon|^2 d_\varepsilon(|m_\varepsilon|^2)^{-p-1}$.
	Finally, Lemma~\ref{lem:stray field mapping} and
	$|m_\varepsilon|\leq1$ imply
	\[
	\norm{H_{\mathrm d}[m_\varepsilon]}
	{L^{p+1}(\mathbb R^3)}^{p+1}
	\leq
	C_p
	\norm{m_\varepsilon}{L^{p+1}(\mathscr D)}^{p+1}
	\leq C_p.
	\]
	Thus, \eqref{equ:inverse distance stray estimate} follows after
	renaming $\eta$.

	Substituting
	\eqref{equ:inverse distance anisotropy estimate} and
	\eqref{equ:inverse distance stray estimate} into
	\eqref{equ:inverse distance weighted identity}, and choosing
	$\eta>0$ sufficiently small compared with $\theta$, we obtain
	\begin{align}
		&\frac{\mathrm d}{\mathrm dt}
		\int_{\mathscr D}
		\Psi_{p,\varepsilon}(|m_\varepsilon|^2) \dx
		+
		\int_{\mathscr D}
		\left.
		\frac{\mathrm d}{\mathrm ds}
		\bigl[d_\varepsilon(s)^{-p}\bigr]
		\right|_{s=|m_\varepsilon|^2}
		\bigl|\nabla|m_\varepsilon|^2\bigr|^2 \dx
		\nonumber\\
		&\quad+
		2
		\int_{\mathscr D}
		d_\varepsilon(|m_\varepsilon|^2)^{-p}
		|\nabla m_\varepsilon|^2\,\mathrm dx
		+
		2\theta
		\int_{\mathscr D}
		|m_\varepsilon|^2
		d_\varepsilon(|m_\varepsilon|^2)^{-p-1} \dx
		\leq
		C_p.
		\label{equ:inverse distance differential estimate first}
	\end{align}
	We now convert the last term on the left-hand side into the full singular weight. On
	$\{|m_\varepsilon|^2>1/2\}$,
	\[
	d_\varepsilon(|m_\varepsilon|^2)^{-p-1}
	\leq
	2|m_\varepsilon|^2
	d_\varepsilon(|m_\varepsilon|^2)^{-p-1},
	\]
	whereas on $\{|m_\varepsilon|^2\leq1/2\}$, we have
	\[
	d_\varepsilon(|m_\varepsilon|^2)^{-p-1}
	\leq
	\bigl(2\alpha q(1/\sqrt2)\bigr)^{-p-1}.
	\]
	Hence,
	\[
	\int_{\mathscr D}
	d_\varepsilon(|m_\varepsilon|^2)^{-p-1}\dx
	\leq
	C_p
	+
	2
	\int_{\mathscr D}
	|m_\varepsilon|^2
	d_\varepsilon(|m_\varepsilon|^2)^{-p-1}\dx.
	\]
	Therefore, dropping the nonnegative gradient terms from
	\eqref{equ:inverse distance differential estimate first}, we obtain
	\begin{equation}
		\frac{\mathrm d}{\mathrm dt}
		\int_{\mathscr D}
		\Psi_{p,\varepsilon}(|m_\varepsilon|^2)\dx
		+
		c_p
		\int_{\mathscr D}
		d_\varepsilon(|m_\varepsilon|^2)^{-p-1}\dx
		\leq
		C_p.
		\label{equ:inverse distance differential estimate}
	\end{equation}
	
	It remains to express the above inequality as a scalar differential inequality.
	From \eqref{equ:Psi comparison},
	\[
	\int_{\mathscr D}
	\Psi_{p,\varepsilon}(|m_\varepsilon|^2) \dx
	\leq
	C
	\int_{\mathscr D}
	d_\varepsilon(|m_\varepsilon|^2)^{-p} \dx
	\leq
	C |\mathscr D|^{\frac{1}{p+1}}
	\left(
	\int_{\mathscr D}
	d_\varepsilon(|m_\varepsilon|^2)^{-p-1} \dx
	\right)^{\frac{p}{p+1}},
	\]
	where in the last step we used H\"older's inequality.
	Consequently, we have
	\begin{equation}
		\int_{\mathscr D}
		d_\varepsilon(|m_\varepsilon|^2)^{-p-1} \dx
		\geq
		c_p
		\left(
		\int_{\mathscr D}
		\Psi_{p,\varepsilon}(|m_\varepsilon|^2) \dx
		\right)^{1+ \frac1p}.
		\label{equ:inverse distance holder}
	\end{equation}
	Thus, by writing
	\[
	E_{\varepsilon,p}(t)
	\eqqcolon
	\int_{\mathscr D}
	\Psi_{p,\varepsilon}(|m_\varepsilon(t)|^2)\dx,
	\]
	and inserting \eqref{equ:inverse distance holder} into
	\eqref{equ:inverse distance differential estimate}, we arrive at
	\begin{equation}
		E_{\varepsilon,p}'(t)
		+
		c_p \left[E_{\varepsilon,p}(t)\right]^{1+ \frac1p}
		\leq
		C_p.
		\label{equ:inverse distance ODE}
	\end{equation}
	By a form of the Gronwall lemma with algebraic dissipation~\cite[Lemma~5.1]{Tem97}, we obtain
	\begin{equation}\label{equ:bound E eps p unif gron}
	E_{\varepsilon,p}(t)
	\leq
	C_p\bigl(1+t^{-p}\bigr),
	\qquad t>0.
	\end{equation}
	Finally, \eqref{equ:Psi comparison} and
	\eqref{equ:d epsilon equivalence} imply
	\[
	\bigl(
	\varepsilon+1-|m_\varepsilon|^2
	\bigr)^{-p}
	\leq
	C_p
	\left(
	1+
	\Psi_{p,\varepsilon}(|m_\varepsilon|^2)
	\right).
	\]
	Integrating over $\mathscr D$ and applying \eqref{equ:bound E eps p unif gron} then yields \eqref{equ:viscous inverse distance}, as claimed.
\end{proof}

The inverse-distance estimate obtained above controls the approach to saturation only in an integral sense. We next complement it with a positive-time spatial regularity estimate. In particular, choosing $r>3/2$ will later give the spatial continuity needed to turn the integral control of $1-|m_\varepsilon|^2$ into a uniform pointwise separation from saturation.

\begin{lemma}
	\label{lem:positive time W2r viscous}
	Let $\tau>0$ and $r\in (1,2)$. For every time interval $I\subset[\tau,\infty)$ with $|I|\leq1$,
	\begin{equation}
		\int_I
		\norm{m_\varepsilon(t)}{W^{2,r}(\mathscr D)}^r
		\dt
		\leq
		C_{\tau,r},
		\label{equ:viscous W2r estimate}
	\end{equation}
	where $C_{\tau,r}$ is independent of $\varepsilon$.
\end{lemma}

\begin{proof}
	By \eqref{equ:mobility radial coercivity}, for a.e. $(t,x)\in (0,T)\times \mathscr{D}$,
	\[
	\alpha
	\Delta m_\varepsilon\cdot
	\mathbb M(m_\varepsilon)
	\Delta m_\varepsilon
	+
	\varepsilon|\Delta m_\varepsilon|^2
	\geq
	c
	\bigl(
	q(|m_\varepsilon|)+\varepsilon
	\bigr)
	|\Delta m_\varepsilon|^2,
	\]
	where $c>0$ is independent of $\varepsilon$. Hence, by H\"older's
	inequality with exponents $2/r$ and $2/(2-r)$,
	\begin{align*}
		\int_{I\times\mathscr D}
		|\Delta m_\varepsilon|^r \dx\dt
		&=
		\int_{I\times\mathscr D}
		\Bigl[
		\bigl(q(|m_\varepsilon|)+\varepsilon\bigr)^{\frac12}
		|\Delta m_\varepsilon|
		\Bigr]^r
		\bigl(q(|m_\varepsilon|)+\varepsilon\bigr)^{-\frac{r}{2}}
		\dx\dt
		\\
		&\leq
		C
		\left(
		\int_{I\times\mathscr D}
		\bigl(q(|m_\varepsilon|)+\varepsilon\bigr)
		|\Delta m_\varepsilon|^2
		\dx\dt
		\right)^{\frac{r}{2}}
		\left(
		\int_{I\times\mathscr D}
		\bigl(q(|m_\varepsilon|)+\varepsilon\bigr)^{-\frac{r}{2-r}}
		\dx\dt
		\right)^{1-\frac{r}{2}}
		\\
		&\leq C_{\tau,r}.
	\end{align*}
	Indeed, since $|I|\leq1$, the first factor is uniformly bounded by
	\eqref{equ:global viscous estimates}. For the second factor,
	\eqref{equ:q saturation equivalence} gives
	$q(|m_\varepsilon|)+\varepsilon \simeq \varepsilon+1-|m_\varepsilon|^2$. On the other hand, Proposition~\ref{prop:inverse distance viscous} with
	$p=r/(2-r)$, yields for $t\in I\subset[\tau,\infty)$,
	\[
	\int_{\mathscr D}
	\bigl(q(|m_\varepsilon|)+\varepsilon\bigr)^{-r/(2-r)}
	\dx
	\leq
	C_{\tau,r}.
	\]
	Thus, since $|I|\leq1$, we obtain
	\[
	\int_I
	\norm{\Delta m_\varepsilon(t)}{L^r(\mathscr D)}^r
	\dt
	\leq
	C_{\tau,r}.
	\]
	The standard $W^{2,r}$-regularity estimate for the Neumann Laplacian gives
	\[
	\norm{m_\varepsilon(t)}{W^{2,r}(\mathscr D)}
	\leq
	C
	\left(
	\norm{\Delta m_\varepsilon(t)}{L^r(\mathscr D)}
	+
	\norm{m_\varepsilon(t)}{L^r(\mathscr D)}
	\right).
	\]
	Since $|m_\varepsilon| \leq 1$ a.e. by
	Lemma~\ref{lem:viscous invariant ball}, integration over $I$ proves
	\eqref{equ:viscous W2r estimate}.
\end{proof}


\subsection{Uniform positive-time separation}

We now upgrade the preceding positive-time estimates to a uniform pointwise separation from saturation. Given $\tau>0$, the inverse-distance bound \eqref{equ:viscous inverse distance} and the $W^{2,r}$-regularity \eqref{equ:viscous W2r estimate} first yield a time $s_\tau\in(0,\tau)$ at which the solution is uniformly separated from $|m|=1$. A scalar comparison argument, initiated at $s_\tau$, then propagates this separation forward in time and therefore yields a
uniform gap on $[\tau,\infty)$.

\begin{proposition}
	\label{prop:uniform viscous separation}
	For every $\tau>0$, there exist $\rho_\tau\in(0,1)$ and
	$\varepsilon_\tau\in(0,1]$ such that
	\begin{equation}
		\operatorname*{ess\,sup}_{x\in \mathscr{D}}
		|m_\varepsilon (t,x)|
		\leq
		\rho_\tau
		\qquad
		\text{for every } t\in[\tau,\infty)
		\label{equ:viscous uniform separation}
	\end{equation}
	for every $0<\varepsilon\leq\varepsilon_\tau$.
\end{proposition}

\begin{proof}
	We first obtain a uniform pointwise gap at a suitable time preceding
	each $t\geq\tau$, and then propagate this gap to $t$ by a scalar
	comparison argument.
	
	Fix $r\in\left(\frac32,2\right)$ and $\mu:=2-\frac3r >0$, and choose $p>\frac{3}{\mu}$. By the Sobolev--Morrey embedding, we have $W^{2,r}(\mathscr D)
	\hookrightarrow
	C^{0,\mu}(\overline{\mathscr D})$.
	For arbitrary $t\geq\tau$, define the closed interval
	\[
	I_t
	:=
	\bigl[\max\{\tau/4,t-1\},t\bigr].
	\]
	Then $I_t\subset[\tau/4, \infty)$ and $\min\left\{1,\frac{3\tau}{4}\right\} \leq |I_t|\leq 1$.
	
	By Lemma~\ref{lem:positive time W2r viscous}, applied with $\tau/4$, we have
	\[
	\int_{I_t}
	\norm{m_\varepsilon(s)}{W^{2,r}(\mathscr D)}^r\ds
	\leq C_{\tau/4,r}.
	\]
	Set $\ell_\tau \coloneq \min\left\{1,\frac{3\tau}{4}\right\}$,
	so that $|I_t|\geq\ell_\tau>0$. By Chebyshev's inequality,
	\begin{align*}
		\left|
		\left\{
		s\in I_t:
		\norm{m_\varepsilon(s)}{W^{2,r}(\mathscr D)}^r
		>
		\frac{2C_{\tau/4,r}}{\ell_\tau}
		\right\}
		\right|
		&\leq
		\frac{\ell_\tau}{2C_{\tau/4,r}}
		\int_{I_t}
		\norm{m_\varepsilon(s)}{W^{2,r}(\mathscr D)}^r\ds
		\leq
		\frac{\ell_\tau}{2}.
	\end{align*}
	Hence, the set
	\[
	\left\{
	s\in I_t:
	\norm{m_\varepsilon(s)}{W^{2,r}(\mathscr D)}
	\leq
	\left(
	\frac{2C_{\tau/4,r}}{\ell_\tau}
	\right)^{\frac1r}
	\right\}
	\]
	has measure at least $\ell_\tau/2$.
	On the other hand, Proposition~\ref{prop:inverse distance viscous} implies, for almost every $s\in I_t$,
	\[
	\int_{\mathscr D}
	\bigl(
	\varepsilon+1-|m_\varepsilon(s)|^2
	\bigr)^{-p}\dx
	\leq
	C_p
	\left[
	1+\left(\frac4\tau\right)^p
	\right],
	\]
	since $I_t\subset[\tau/4,\infty)$. The exceptional set where this
	estimate may fail has measure zero. We may therefore choose
	$s_{\varepsilon,t}\in I_t$ such that
	\begin{equation}\label{equ:inv dist m s eps}
	\norm{m_\varepsilon(s_{\varepsilon,t})}{W^{2,r}(\mathscr D)}
	\leq
	\left(
	\frac{2C_{\tau/4,r}}{\ell_\tau}
	\right)^{\frac1r}
	\quad \text{and} \quad 
	\int_{\mathscr D}
	\bigl(
	\varepsilon+1-|m_\varepsilon(s_{\varepsilon,t})|^2
	\bigr)^{-p}\dx
	\leq
	C_p
	\left[
	1+\left(\frac4\tau\right)^p
	\right].
	\end{equation}

	By the choice of $s_{\varepsilon,t}$ and the Sobolev--Morrey embedding,
	\[
	\norm{m_\varepsilon(s_{\varepsilon,t})}
	{C^{0,\mu}(\overline{\mathscr D})}
	\leq
	C
	\left(
	\frac{2C_{\tau/4,r}}{\ell_\tau}
	\right)^{1/r}
	\eqqcolon K_\tau,
	\]
	where $K_\tau$ is independent of $\varepsilon$ and $t$. Since
	$m_\varepsilon(s_{\varepsilon,t})$ is continuous in $x$ and
	$|m_\varepsilon(s_{\varepsilon,t})|\leq1$ a.e., the latter inequality
	holds in fact everywhere on $\overline{\mathscr D}$. Therefore, for any $x,y\in \overline{\mathscr D}$,
	\begin{align}\label{equ:Holder eps 1 m eps}
		\left|
		\bigl(
		\varepsilon+1-|m_\varepsilon(s_{\varepsilon,t},x)|^2
		\bigr)
		-
		\bigl(
		\varepsilon+1-|m_\varepsilon(s_{\varepsilon,t},y)|^2
		\bigr)
		\right|
		&\leq
		2
		\left|
		m_\varepsilon(s_{\varepsilon,t},x)
		-
		m_\varepsilon(s_{\varepsilon,t},y)
		\right|
		\nonumber\\
		&\leq
		2K_\tau |x-y|^\mu .
	\end{align}
	Moreover, $0< \varepsilon+1-|m_\varepsilon(s_{\varepsilon,t},x)|^2 \leq 2$.
	Thus, the $C^{0,\mu}(\overline{\mathscr D})$-norm of
	$x\mapsto\varepsilon+1-|m_\varepsilon(s_{\varepsilon,t},x)|^2$
	is bounded independently of $\varepsilon$ and $t$.
	
	Now, set
	\begin{equation}\label{equ:delta eps t}
	\delta_{\varepsilon,t}
	\coloneq
	\min_{x\in\overline{\mathscr D}}
	\left\{
	\varepsilon+1
	-
	|m_\varepsilon(s_{\varepsilon,t},x)|^2
	\right\}>0,
	\end{equation}
	and let $x_{\varepsilon,t}\in\overline{\mathscr D}$ be a point at which
	this minimum is attained.  Since $\mathscr D$ has sufficiently smooth boundary, it satisfies a uniform interior cone condition: there exist a height
	$R_{\mathscr D}>0$ and an opening angle $\vartheta_{\mathscr D}\in(0,\pi/2)$, such that, for every $x\in\overline{\mathscr D}$ and every $R\in (0, R_{\mathscr D}]$,
	the set $\mathscr D\cap B_R(x)$ contains a cone of height $R$ and opening
	angle $\vartheta_{\mathscr D}$. Consequently, there exists
	$c_{\mathscr D}>0$, depending only on $\vartheta_{\mathscr D}$, such that
	\begin{equation}
		|\mathscr D\cap B_R(x)|
		\geq
		c_{\mathscr D}R^3
		\qquad
		\text{for every } x\in\overline{\mathscr D}
		\text{ and } R\in (0, R_{\mathscr D}].
		\label{equ:uniform cone volume}
	\end{equation}
	
	From the H\"older estimate \eqref{equ:Holder eps 1 m eps},
	\[
	\varepsilon+1
	-
	|m_\varepsilon(s_{\varepsilon,t},x)|^2
	\leq
	\delta_{\varepsilon,t}
	+
	2K_\tau |x-x_{\varepsilon,t}|^\mu
	\qquad
	\text{for every }x\in\overline{\mathscr D}.
	\]
	Suppose first that $\delta_{\varepsilon,t} \leq 2K_\tau R_{\mathscr D}^{\mu}$.
	We then choose the radius $R \coloneq \bigl( \delta_{\varepsilon,t}/(2K_\tau)
	\bigr)^{\frac{1}{\mu}}$.
	By the preceding assumption, $0<R\leq R_{\mathscr D}$, and thus, for every
	$x\in\mathscr D\cap B_R(x_{\varepsilon,t})$,
	\begin{align*}
		\varepsilon+1
		-
		|m_\varepsilon(s_{\varepsilon,t},x)|^2
		&\leq
		\delta_{\varepsilon,t}
		+
		2K_\tau R^\mu
		=
		2\delta_{\varepsilon,t}.
	\end{align*}
	Using \eqref{equ:uniform cone volume} and the inverse-distance estimate \eqref{equ:inv dist m s eps}, we obtain
	\begin{align*}
		C_p
		\left(
		1+\left(\frac4\tau\right)^p
		\right)
		&\geq
		\int_{\mathscr D}
		\bigl(
		\varepsilon+1
		-
		|m_\varepsilon(s_{\varepsilon,t})|^2
		\bigr)^{-p}\dx
		\\
		&\geq
		\int_{\mathscr D\cap B_R(x_{\varepsilon,t})}
		\bigl(
		\varepsilon+1
		-
		|m_\varepsilon(s_{\varepsilon,t})|^2
		\bigr)^{-p}\dx
		\\
		&\geq
		(2\delta_{\varepsilon,t})^{-p}
		|\mathscr D\cap B_R(x_{\varepsilon,t})|
		\\
		&\geq
		c_{\mathscr D}
		2^{-p}
		R^3
		\delta_{\varepsilon,t}^{-p}
		\\
		&=
		c_{\mathscr D}
		2^{-p}
		(2K_\tau)^{-\frac{3}{\mu}}
		\delta_{\varepsilon,t}^{\frac{3}{\mu}-p}.
	\end{align*}
	Since $p-3/\mu>0$, we infer that
	\[
	\delta_{\varepsilon,t}
	\geq
	\left[
	\frac{
		c_{\mathscr D}2^{-p}(2K_\tau)^{-\frac{3}{\mu}}
	}{
		C_p\bigl(1+(4/\tau)^p \bigr)
	}
	\right]^{\frac{1}{p-\frac{3}{\mu}}}.
	\]
	If instead $\delta_{\varepsilon,t} > 2K_\tau R_{\mathscr D}^{\mu}$, then a positive lower bound is immediate. Hence, in either case, we have $\delta_{\varepsilon,t} \geq \delta_\tau>0$, where $\delta_\tau$ is independent of $\varepsilon$ and $t$.
	Finally, recall that from \eqref{equ:delta eps t},
	\[
	\varepsilon+1
	-
	|m_\varepsilon(s_{\varepsilon,t},x)|^2
	\geq
	\delta_\tau
	\qquad
	\forall x\in\overline{\mathscr D}.
	\]
	Hence, whenever $\varepsilon\leq\delta_\tau/2$,
	\begin{equation}
		|m_\varepsilon(s_{\varepsilon,t},x)|^2
		\leq
		1+\varepsilon-\delta_\tau
		\leq
		1-\frac{\delta_\tau}{2},
		\qquad
		\forall x\in\overline{\mathscr D}.
		\label{equ:good time pointwise separation}
	\end{equation}

	We now propagate \eqref{equ:good time pointwise separation} from
	$s_{\varepsilon,t}$ to $t$. Define
	\begin{equation}\label{equ:b epsilon}
	b_\varepsilon(s)
	\coloneq
	\int_0^s
	\frac{1}{2\alpha q(\sqrt \sigma)+\varepsilon}\,\mathrm{d} \sigma,
	\qquad
	s\in [0,1].
	\end{equation}
	Since $q$ is decreasing and non-negative, we have
	\[
	b_\varepsilon'(s)
	=
	\frac{1}{2\alpha q(\sqrt s)+\varepsilon}>0,
	\qquad
	b_\varepsilon''(s)\geq 0,
	\]
	so $b_\varepsilon$ is strictly increasing and convex on $[0,1]$.
	
	Using \eqref{equ:viscous scalar equation} and the chain rule in temporal variable, we
	obtain, for a.e. $s>0$,
	\begin{align}
		&\partial_s b_\varepsilon(|m_\varepsilon|^2)
		-
		\Delta|m_\varepsilon|^2
		\nonumber\\
		&=
		-2|\nabla m_\varepsilon|^2
		+
		\frac{4q(|m_\varepsilon|)}
		{2\alpha q(|m_\varepsilon|)+\varepsilon}
		m_\varepsilon\cdot
		\bigl(
		\mathsf Bm_\varepsilon
		+
		H_{\mathrm d}[m_\varepsilon]
		\bigr)
		-
		\frac{4\theta|m_\varepsilon|^2}
		{2\alpha q(|m_\varepsilon|)+\varepsilon}
		\nonumber\\
		&\leq
		\frac2\alpha
		\left(
		\norm{\mathsf B}{}
		+
		\norm{H_{\mathrm d}[m_\varepsilon](s)}
		{L^\infty(\mathscr D)}
		\right),
		\label{equ:entropy comparison inequality}
	\end{align}
	where we used the fact that $|m_\varepsilon|\leq1$ and
	\[
	\frac{q(|m_\varepsilon|)}
	{2\alpha q(|m_\varepsilon|)+\varepsilon}
	\leq
	\frac1{2\alpha},
	\]
	while the first and last terms in the equality are non-positive. 
	The purpose of \eqref{equ:entropy comparison inequality} is to propagate the
	pointwise gap obtained at $s_{\varepsilon,t}$. We first show that the
	right-hand side has a uniformly bounded integral over
	$[s_{\varepsilon,t},t]$, then construct a scalar comparison function
	$z_\varepsilon$ starting from $1-\delta_\tau/2$ whose evolution accounts
	for this forcing. A logarithmic lower bound for $b_\varepsilon$ near $1$, uniform for sufficiently small $\varepsilon$, then prevents the comparison function from reaching the saturation boundary. A comparison argument then give
	\[
	|m_\varepsilon(s)|^2
	\leq
	z_\varepsilon(s)
	\leq
	1-\eta_\tau
	\qquad
	\text{on }[s_{\varepsilon,t},t].
	\]
	
	We now proceed with our aim.
	By Lemmas~\ref{lem:stray field mapping} and~\ref{lem:positive time W2r viscous}, the embedding $W^{2,r}(\mathscr{D}) \hookrightarrow L^\infty(\mathscr{D})$, and the fact that
	$[s_{\varepsilon,t},t]\subset I_t$ has length at most one, we have
	\begin{align}
		\int_{s_{\varepsilon,t}}^t
		\norm{H_{\mathrm d}[m_\varepsilon](s)}
		{L^\infty(\mathscr D)}\ds
		&\leq
		C
		\int_{s_{\varepsilon,t}}^t
		\norm{m_\varepsilon(s)}{W^{2,r}(\mathscr D)}\ds
		\nonumber\\
		&\leq
		C
		|t-s_{\varepsilon,t}|^{1-\frac1r}
		\left(
		\int_{s_{\varepsilon,t}}^t
		\norm{m_\varepsilon(s)}{W^{2,r}(\mathscr D)}^r\ds
		\right)^{1/r}
		\nonumber\\
		&\leq
		\widetilde{C_\tau}.
		\label{equ:stray field comparison interval}
	\end{align}
	Consequently,
	\begin{equation}
		\frac2\alpha
		\int_{s_{\varepsilon,t}}^t
		\left(
		\norm{\mathsf B}{}
		+
		\norm{H_{\mathrm d}[m_\varepsilon](s)}
		{L^\infty(\mathscr D)}
		\right)\ds
		\leq
		\widetilde{C_\tau},
		\label{equ:comparison forcing bound}
	\end{equation}
	where $\widetilde{C_\tau}$ is a highlighted constant independent of $\varepsilon$ and $t$.
	
	We next choose a uniform upper barrier below $1$. By
	\eqref{equ:q saturation equivalence},
	\[
	2\alpha q(\sqrt s)+\varepsilon
	\leq
	C(1-s+\varepsilon),
	\qquad \forall s\in [0,1].
	\]
	Hence, for $\eta\in (0,\delta_\tau/2)$, we obtain from \eqref{equ:b epsilon} that
	\begin{align}
		b_\varepsilon(1-\eta)
		-
		b_\varepsilon
		\left(1-\frac{\delta_\tau}{2}\right)
		&=
		\int_{1-\delta_\tau/2}^{1-\eta}
		\frac{1}{2\alpha q(\sqrt s)+\varepsilon}\ds
		\nonumber\\
		&\geq
		c
		\int_{1-\delta_\tau/2}^{1-\eta}
		\frac{1}{1-s+\varepsilon}\ds
		=
		c
		\log
		\left(
		\frac{\delta_\tau/2+\varepsilon}
		{\eta+\varepsilon}
		\right).
		\label{equ:b epsilon logarithmic gap}
	\end{align}
	Choose $\eta_\tau\in(0,\delta_\tau/2)$ sufficiently small such that
	\begin{equation}\label{equ:tilde Ctau}
	c\log\left(\frac{\delta_\tau}{4\eta_\tau}\right)
	>
	\widetilde{C_\tau},
	\end{equation}
	and then choose $\varepsilon_\tau$ such that
	\[
	0<\varepsilon_\tau
	\leq
	\min\left\{
	1,\frac{\delta_\tau}{2},\eta_\tau
	\right\}.
	\]
	For every $\varepsilon\in (0, \varepsilon_\tau]$, we then have
	\[
	\frac{\delta_\tau/2+\varepsilon}
	{\eta_\tau+\varepsilon}
	\geq
	\frac{\delta_\tau/2}{2\eta_\tau}
	=
	\frac{\delta_\tau}{4\eta_\tau},
	\]
	thus \eqref{equ:b epsilon logarithmic gap} and \eqref{equ:tilde Ctau} gives
	\begin{equation}
		b_\varepsilon(1-\eta_\tau)
		-
		b_\varepsilon
		\left(1-\frac{\delta_\tau}{2}\right)
		>
		\widetilde{C_\tau}.
		\label{equ:b epsilon dominates forcing}
	\end{equation}
	Now, for $s\in[s_{\varepsilon,t},t]$, set
	\begin{equation}\label{equ:Y epsilon s}
	Y_\varepsilon(s)
	\coloneq
	b_\varepsilon\left(1-\frac{\delta_\tau}{2}\right)
	+
	\frac2\alpha
	\int_{s_{\varepsilon,t}}^s
	\left(
	\norm{\mathsf B}{}
	+
	\norm{H_{\mathrm d}[m_\varepsilon](\sigma)}
	{L^\infty(\mathscr D)}
	\right)\mathrm{d}\sigma.
	\end{equation}
	For every $s\in[s_{\varepsilon,t},t]$, by
	\eqref{equ:comparison forcing bound},
	\eqref{equ:b epsilon dominates forcing}, and the fact that $\eta_\tau\in(0,\delta_\tau/2)$,
	\[
	b_\varepsilon\left(1-\frac{\delta_\tau}{2}\right)
	\leq
	Y_\varepsilon(s)
	<
	b_\varepsilon(1-\eta_\tau).
	\]
	Since $b_\varepsilon$ is continuous and strictly increasing,
	the intermediate value theorem and strict monotonicity imply that there
	exists a unique $z_\varepsilon(s)\in \left[ 1-\frac{\delta_\tau}{2}, 1-\eta_\tau\right)$
	such that $b_\varepsilon(z_\varepsilon(s))= Y_\varepsilon(s)$.
	Equivalently, $z_\varepsilon(s)= b_\varepsilon^{-1}(Y_\varepsilon(s))$.
	
	Moreover, the map $s\mapsto Y_\varepsilon(s)$ is absolutely continuous.
	Since $b_\varepsilon^{-1}$ is continuously differentiable on the relevant
	interval, $z_\varepsilon$ is absolutely continuous and
	\[
	\partial_s b_\varepsilon(z_\varepsilon(s))
	=
	\partial_s Y_\varepsilon (s)
	=
	\frac2\alpha
	\left(
	\norm{\mathsf B}{}
	+
	\norm{H_{\mathrm d}[m_\varepsilon](s)}
	{L^\infty(\mathscr D)}
	\right)
	\]
	for a.e. $s\in[s_{\varepsilon,t},t]$.
	Subtracting this identity from
	\eqref{equ:entropy comparison inequality} gives
	\begin{equation*}
		\partial_s
		\left[
		b_\varepsilon(|m_\varepsilon|^2)
		-
		b_\varepsilon(z_\varepsilon)
		\right]
		-
		\Delta|m_\varepsilon|^2
		\leq 0
	\end{equation*}
	in the distributional sense on
	$(s_{\varepsilon,t},t)\times\mathscr D$; that is,
	\begin{equation}\label{equ:comparison difference inequality}
	-\int_{s_{\varepsilon,t}}^t
	\inpro{
	b_\varepsilon(|m_\varepsilon|^2)
	-
	b_\varepsilon(z_\varepsilon)}{
	\partial_s\varphi} \ds
	+
	\int_{s_{\varepsilon,t}}^t
	\inpro{
	\nabla|m_\varepsilon|^2}{
	\nabla\varphi}\ds
	\leq 0
	\end{equation}
	for every nonnegative
	$\varphi\in C_c^\infty\bigl((s_{\varepsilon,t},t);H^1(\mathscr D)\bigr)$.
	
	Set $F_\varepsilon \coloneq b_\varepsilon(|m_\varepsilon|^2) -
	b_\varepsilon(z_\varepsilon)$.
	For fixed $\varepsilon>0$, the boundedness of $b_\varepsilon'$ on $[0,1]$, together with the regularity of $|m_\varepsilon|^2$ and the absolute continuity of $z_\varepsilon$, gives
	\[
	F_\varepsilon
	\in
	W^{1,r}\bigl(s_{\varepsilon,t},t;L^2(\mathscr D)\bigr)
	\cap
	L^2\bigl(s_{\varepsilon,t},t;H^1(\mathscr D)\bigr).
	\]
	Hence, $(F_\varepsilon)_+$ has the same regularity and may be used, after
	multiplication by a standard nonnegative time cut-off, as a test function
	in \eqref{equ:comparison difference inequality}. Using the Sobolev chain
	rule, the monotonicity of $b_\varepsilon$, and the fact that
	$z_\varepsilon$ is independent of $x$, we obtain for every
	$s\in(s_{\varepsilon,t},t)$
	\begin{align}\label{equ:F eps plus ineq}
		\frac12
		\norm{(F_\varepsilon(s))_+}{L^2(\mathscr D)}^2
		-
		\frac12
		\norm{(F_\varepsilon(s_{\varepsilon,t}))_+}{L^2(\mathscr D)}^2
		+
		\int_{s_{\varepsilon,t}}^s
		\int_{\{
			F_\varepsilon>0
			\}}
		b_\varepsilon'(|m_\varepsilon|^2)
		\bigl|\nabla|m_\varepsilon|^2\bigr|^2
		\dx\ds
		\leq 0,
	\end{align}
	where we used the fact that $\nabla(F_\varepsilon)_+ =
	\mathbf 1_{\{F_\varepsilon>0\}}
	b_\varepsilon'(|m_\varepsilon|^2)
	\nabla|m_\varepsilon|^2$ and that $z_\varepsilon$ is independent of $x$. Moreover, since
	$b_\varepsilon$ is strictly increasing,
	$\{F_\varepsilon>0\}$ is equivalent with $\{|m_\varepsilon|^2 > z_\varepsilon\}$. However, by \eqref{equ:good time pointwise separation}, using the injectivity of $b_\varepsilon$ and the fact that $b_\varepsilon(z_\varepsilon(s_{\varepsilon,t}))= Y_\varepsilon(s_{\varepsilon,t})= b_\varepsilon(1-\delta_\tau/2)$, we infer
	\[
	|m_\varepsilon(s_{\varepsilon,t},x)|^2
	\leq
	1-\frac{\delta_\tau}{2}
	=
	z_\varepsilon(s_{\varepsilon,t})
	\qquad
	\text{for every }x\in\overline{\mathscr D},
	\]
	and hence $(F_\varepsilon(s_{\varepsilon,t}))_+=0$.
	We thus conclude from \eqref{equ:F eps plus ineq} that
	\[
	\norm{(F_\varepsilon(s))_+}{L^2(\mathscr D)}^2=0
	\qquad
	\text{for every }s\in[s_{\varepsilon,t},t].
	\]
	Therefore, $b_\varepsilon \bigl(|m_\varepsilon(s,x)|^2\bigr) \leq b_\varepsilon\bigl(z_\varepsilon(s)\bigr)$ for a.e. $x\in\mathscr D$.
	Since $b_\varepsilon$ is strictly increasing, we obtain
	\[
	|m_\varepsilon(s,x)|^2
	\leq
	z_\varepsilon(s)
	<
	1-\eta_\tau
	\]
	for a.e. $x\in\mathscr D$ and every
	$s\in[s_{\varepsilon,t},t]$. Since $t\geq\tau$ was arbitrary,
	\[
	\operatorname*{ess\,sup}_{x\in\mathscr D}
	|m_\varepsilon(t,x)|
	\leq
	\rho_\tau
	\coloneq
	\sqrt{1-\eta_\tau}
	<1
	\qquad
	\text{for every }t\geq\tau.
	\]
	This completes the proof of \eqref{equ:viscous uniform separation}.
\end{proof}


\subsection{Passage to the limit}

Finally, we are in a position to prove Theorem~\ref{thm:global existence full nonlinear LLB}.

\begin{proof}[Proof of Theorem~\ref{thm:global existence full nonlinear LLB}]
	We first consolidate the compactness properties of $m_\varepsilon$ and
	identify the limit equation, then prove properties {\rm (i)}--{\rm
		(iv)}.
	
	Fix $T>0$. By \eqref{equ:global viscous estimates}, we have
	\[
	\sup_{t\in[0,T]}
	\norm{m_\varepsilon(t)}{H^1(\mathscr D)}^2
	\leq C,
	\]
	while, covering $[0,T]$ by finitely many intervals of length at most
	one,
	\[
	\int_0^T
	\norm{\partial_t m_\varepsilon(s)}{L^2(\mathscr D)}^2\ds
	+
	\varepsilon
	\int_0^T
	\norm{\Delta m_\varepsilon(s)}{L^2(\mathscr D)}^2\ds
	\leq C_T,
	\]
	where $C_T$ is independent of $\varepsilon$. Hence, after passing to
	a subsequence,
	\[
	\begin{aligned}
		m_\varepsilon
		&\overset{*}{\rightharpoonup} m
		&&\text{in }L^\infty\bigl(0,T;H^1(\mathscr D)\bigr),
		\\
		\partial_t m_\varepsilon
		&\rightharpoonup \partial_t m
		&&\text{in }L^2\bigl(0,T;L^2(\mathscr D)\bigr).
	\end{aligned}
	\]
	Moreover, for every $s,t\in[0,T]$,
	\[
	\norm{m_\varepsilon(t)-m_\varepsilon(s)}{L^2(\mathscr D)}
	\leq
	|t-s|^{\frac12}
	\norm{\partial_t m_\varepsilon}{L^2(0,T;L^2(\mathscr D))}
	\leq
	C_T|t-s|^{\frac12}.
	\]
	Since $H^1(\mathscr D)\Subset L^2(\mathscr D)$, the
	Arzel\`a--Ascoli theorem yields
	\begin{equation}
		m_\varepsilon
		\to
		m
		\qquad
		\text{strongly in }C\bigl([0,T];L^2(\mathscr D)\bigr).
		\label{equ:epsilon strong L2}
	\end{equation}
	Since $T>0$ was arbitrary, by a standard diagonal argument over
	$T\in\mathbb N$, we may choose the above subsequence independently of
	$T$ and obtain a single limit $m$ on $[0,\infty)$. In what follows, we
	work with this subsequence without further relabelling.
	In particular, $m_\varepsilon\to m$ strongly in $L^2\bigl((0,T)\times\mathscr D\bigr)$, which implies, after passing to a further subsequence,
	\begin{equation}
		m_\varepsilon(t,x)
		\to
		m(t,x)
		\qquad
		\text{for a.e. }(t,x)\in(0,T)\times\mathscr D.
		\label{equ:epsilon ae convergence}
	\end{equation}
	Since $m_\varepsilon(0)=m_0$, by \eqref{equ:epsilon strong L2} we have
	$m(0)=m_0$ in $L^2(\mathscr D)$. Furthermore, since
	$|m_\varepsilon(t)|\leq1$ a.e. in $\mathscr D$ for every $t\geq0$,
	the $1$-Lipschitz continuity of $v\mapsto(|v|-1)_+$ gives
	\[
	\norm{(|m(t)|-1)_+}{L^2(\mathscr D)}
	\leq
	\norm{m(t)-m_\varepsilon(t)}{L^2(\mathscr D)}
	\to0.
	\]
	Thus,
	\begin{equation}
		|m(t,x)|\leq1
		\qquad
		\text{for a.e. }x\in\mathscr D
		\text{ and every }t\in[0,T].
		\label{equ:limit unit ball}
	\end{equation}
	
	We next obtain the stronger compactness needed to identify the
	thermodynamic equation. Fix $\tau\in(0,T)$. By
	Proposition~\ref{prop:uniform viscous separation}, for all
	sufficiently small $\varepsilon$,
	\begin{equation}
		|m_\varepsilon(t,x)|
		\leq
		\rho_\tau<1
		\qquad
		\text{for a.e. }x\in\mathscr D
		\text{ and every }t\in[\tau,T].
		\label{equ:limit unit ball me}
	\end{equation}
	Since $q$ is decreasing,
	$q(|m_\varepsilon|)\geq q(\rho_\tau)>0$. Therefore, by
	\eqref{equ:mobility radial coercivity} and
	\eqref{equ:global viscous estimates},
	\begin{align}\label{equ:Lap bound by weight}
		\int_\tau^T
		\norm{\Delta m_\varepsilon(s)}{L^2(\mathscr D)}^2\ds
		&\leq
		\frac{1}{2q(\rho_\tau)}
		\int_\tau^T
		\int_{\mathscr D}
		\Delta m_\varepsilon\cdot
		\mathbb M(m_\varepsilon)
		\Delta m_\varepsilon
		\dx\ds
		\leq C_{\tau,T}.
	\end{align}
	The Neumann elliptic estimate and \eqref{equ:limit unit ball me} imply
	\begin{equation}\label{equ:me H2 bound}
	\norm{m_\varepsilon}
	{L^2(\tau,T;H^2(\mathscr D))}
	\leq C_{\tau,T}.
	\end{equation}
	Since the limit $m$ has already been fixed by
	\eqref{equ:epsilon strong L2}, every weak accumulation point in
	$L^2(\tau,T;H^2(\mathscr D))$ necessarily coincides with $m$.
	Hence,
	\[
	m_\varepsilon
	\rightharpoonup
	m
	\qquad
	\text{in }L^2\bigl(\tau,T;H^2(\mathscr D)\bigr).
	\]
	Moreover, \eqref{equ:me H2 bound}, the uniform estimate for
	$\partial_t m_\varepsilon$, and the Aubin--Lions lemma gives relative
	compactness in $L^2(\tau,T;H^1(\mathscr D))$. Again, the limit is
	uniquely identified by \eqref{equ:epsilon strong L2}, so that
	\begin{equation}
		m_\varepsilon
		\to
		m
		\qquad
		\text{strongly in }L^2\bigl(\tau,T;H^1(\mathscr D)\bigr).
		\label{equ:epsilon strong H1}
	\end{equation}
	
	We now collect the remaining convergences needed to pass to the
	nonlinear equation. From \eqref{equ:epsilon strong L2} and the
	$L^2$-boundedness of the stray-field operator,
	\[
	H_{\mathrm d}[m_\varepsilon]
	\to
	H_{\mathrm d}[m]
	\qquad
	\text{strongly in }C\bigl([0,T];L^2(\mathbb R^3)\bigr).
	\]
	Furthermore, by \eqref{equ:global viscous estimates},
	\begin{align*}
		\norm{\varepsilon\Delta m_\varepsilon}
		{L^2(0,T;L^2(\mathscr D))}^2
		&=
		\varepsilon^2
		\int_0^T
		\norm{\Delta m_\varepsilon}{L^2(\mathscr D)}^2\dt
		\leq
		\varepsilon
		\left[
		\varepsilon
		\int_0^T
		\norm{\Delta m_\varepsilon}{L^2(\mathscr D)}^2\dt
		\right]
		\to 0.
	\end{align*}
	By \eqref{equ:limit unit ball me} and
	\eqref{equ:epsilon ae convergence}, we also have
	$|m|\leq\rho_\tau$ a.e. on $(\tau,T)\times\mathscr D$.
	Thus $G(m_\varepsilon)$, $G(m)$, $h(m_\varepsilon)$, and $h(m)$
	are well defined there. Since $G$ is Lipschitz on
	$\{|v|\leq\rho_\tau\}$,
	\[
	G(m_\varepsilon)
	\to
	G(m)
	\qquad
	\text{strongly in }L^2(\tau,T;L^2(\mathscr D)).
	\]
	Combining these convergences with the weak convergence of
	$\Delta m_\varepsilon$ gives
	\begin{equation}
		h(m_\varepsilon)
		\rightharpoonup
		h(m)
		\qquad
		\text{in }L^2(\tau,T;L^2(\mathscr D)).
		\label{equ:h epsilon weak convergence}
	\end{equation}
	Let $\varphi\in L^2((\tau,T)\times\mathscr D;\mathbb R^3)$. By
	\eqref{equ:epsilon ae convergence}, \eqref{equ:limit unit ball me}, and
	the dominated convergence theorem, we have
	\[
	\begin{aligned}
		\varphi\times m_\varepsilon
		&\to \varphi\times m
		&&\text{strongly in }L^2\bigl((\tau,T)\times\mathscr D\bigr),
		\\
		\mathbb M(m_\varepsilon)\varphi
		&\to \mathbb M(m)\varphi
		&&\text{strongly in }L^2\bigl((\tau,T)\times\mathscr D\bigr),
	\end{aligned}
	\]
	noting that $\mathbb M$ is continuous and bounded on the closed
	unit ball.
	Utilising the symmetry of $\mathbb M(m_\varepsilon)$, we then have
	\begin{align*}
		\int_\tau^T
		\inpro{\partial_t m_\varepsilon}{\varphi}\ds
		&=
		\gamma
		\int_\tau^T
		\inpro{h(m_\varepsilon)}
		{\varphi\times m_\varepsilon}\ds
		+
		\int_\tau^T
		\inpro{h(m_\varepsilon)}
		{\mathbb M(m_\varepsilon)\varphi}\ds
		+
		\int_\tau^T
		\inpro{\varepsilon\Delta m_\varepsilon}{\varphi}\ds
		\\
		&\longrightarrow
		\gamma
		\int_\tau^T
		\inpro{h(m)}{\varphi\times m}\ds
		+
		\int_\tau^T
		\inpro{h(m)}{\mathbb M(m)\varphi}\ds.
	\end{align*}
	Hence, we infer that
	\begin{equation}
		\partial_t m
		=
		\gamma m\times h(m)
		+
		\mathbb M(m)h(m)
		\qquad
		\text{in }L^2(\tau,T;L^2(\mathscr D)).
		\label{equ:limit thermodynamic equation}
	\end{equation}
	To recover the weak formulation of the effective field, let
	$\psi\in H^1(\mathscr D)$ and
	$\zeta\in C_c^\infty(\tau,T)$. Since
	$\partial_\nu m_\varepsilon=0$, we have
	\begin{align*}
		\int_\tau^T
		\zeta(s)
		\inpro{h(m_\varepsilon)}{\psi}\ds
		={}&
		-\alpha
		\int_\tau^T
		\zeta(s)
		\inpro{\nabla m_\varepsilon}{\nabla\psi}\ds
		+
		\int_\tau^T
		\zeta(s)
		\inpro{
			\mathsf Bm_\varepsilon
			+
			H_{\mathrm d}[m_\varepsilon]
			-
			\theta G(m_\varepsilon)
		}{\psi}\ds.
	\end{align*}
	Passing to the limit gives
	\begin{equation}\label{equ:limit thermodynamic equation hm}
	\inpro{h(m)}{\psi}
	=
	-\alpha\inpro{\nabla m}{\nabla\psi}
	+
	\inpro{
		\mathsf Bm
		+
		H_{\mathrm d}[m]
		-
		\theta G(m)
	}{\psi},
	\end{equation}
	for every $\psi\in H^1(\mathscr D)$ and a.e.
	$s\in(\tau,T)$. Thus
	\eqref{equ:weak evolution equation} and
	\eqref{equ:weak effective field} hold on $(\tau,T)$, for every $\tau\in (0,T)$. This shows the existence of a global weak solution in the sense of Definition~\ref{def:global weak solution}.
	
	We next prove properties {\rm (i)}--{\rm (ii)}, beginning with the
	uniform separation. Fix $\tau>0$. For all sufficiently small
	$\varepsilon$, Proposition~\ref{prop:uniform viscous separation}
	gives
	\[
	|m_\varepsilon(t,x)|
	\leq
	\rho_\tau
	\qquad
	\text{for a.e. }x\in\mathscr D
	\text{ and every }t\in[\tau,T].
	\]
	Hence, using again the Lipschitz continuity of the positive-part
	map and \eqref{equ:epsilon strong L2},
	\[
	\norm{
		\bigl(|m(t)|-\rho_\tau\bigr)_+
	}{L^2(\mathscr D)}
	\leq
	\norm{m(t)-m_\varepsilon(t)}{L^2(\mathscr D)}
	\to0
	\]
	for every $t\in[\tau,T]$. Therefore,
	\[
	\operatorname*{ess\,sup}_{x\in\mathscr D}
	|m(t,x)|
	\leq
	\rho_\tau
	\qquad
	\text{for every }t\in[\tau,T].
	\]
	Since $T>\tau$ is arbitrary, this proves
	\eqref{equ:positive time separation final}.
	
	We now pass the inverse-distance estimate to the limit. By
	\eqref{equ:epsilon ae convergence} and Fubini's theorem, there is a
	null set of times outside which
	\[
	m_\varepsilon(t,x)\to m(t,x)
	\qquad
	\text{for a.e. }x\in\mathscr D.
	\]
	Since the selected sequence $(\varepsilon_m)_{m\in\bb{N}}$ is countable, enlarging
	this null set by the union of the exceptional time sets in
	Proposition~\ref{prop:inverse distance viscous}, we may assume that
	both this convergence and \eqref{equ:viscous inverse distance} hold
	for every $t$ outside one null set. Fatou's lemma therefore gives,
	for every $p\geq1$,
	\begin{align}
		\int_{\mathscr D}
		\frac{1}{(1-|m(t)|^2)^p}\dx
		&\leq
		\liminf_{\varepsilon\to0}
		\int_{\mathscr D}
		\frac{1}{
			(\varepsilon+1-|m_\varepsilon(t)|^2)^p
		}\dx
		\nonumber\\
		&\leq
		C_p(1+t^{-p})
		\qquad
		\text{for a.e. }t\in(0,T).
		\label{equ:inverse distance final ae}
	\end{align}
	To extend this estimate to every positive time, fix $\tau\in (0,T)$. By the separation result \eqref{equ:positive time separation final} already proved, the map
	$v\mapsto(1-|v|^2)^{-p}$
	is Lipschitz on $\{|v|\leq\rho_\tau\}$. Hence, for
	$s,t\in[\tau,T]$,
	\begin{align*}
		\left|
		\int_{\mathscr D}
		\frac{1}{(1-|m(t)|^2)^p}\dx
		-
		\int_{\mathscr D}
		\frac{1}{(1-|m(s)|^2)^p}\dx
		\right|
		&\leq
		C_{\tau,p}
		\norm{m(t)-m(s)}{L^1(\mathscr D)}.
	\end{align*}
	Thus, the left-hand side of
	\eqref{equ:inverse distance final ae} is continuous in time on
	$[\tau,T]$. Approximating an arbitrary $t\in[\tau,T]$ by times for
	which \eqref{equ:inverse distance final ae} holds gives
	\eqref{equ:inverse distance final} for every $t\in[\tau,T]$.
	Since $\tau>0$ is arbitrary, property {\rm (i)} follows.
	
	Next, we verify property (iii), including \eqref{equ:m0 converge H1}. The differential estimate \eqref{equ:gronwall nabla me} obtained in the proof of Lemma~\ref{lem:uniform viscous estimates}
	gives
	\[
	\norm{\nabla m_\varepsilon(t)}{L^2(\mathscr D)}^2
	\leq
	e^{-4\theta t}
	\norm{\nabla m_0}{L^2(\mathscr D)}^2
	+
	C\bigl(1-e^{-4\theta t}\bigr),
	\qquad t\geq0,
	\]
	where $C$ is independent of $\varepsilon$ and $t$. For every fixed
	$t\in[0,T]$, the strong $L^2$-convergence
	\eqref{equ:epsilon strong L2} and the uniform $H^1$-bound imply $m_\varepsilon(t)
	\rightharpoonup m(t)$ in $H^1(\mathscr{D})$.
	Thus, by weak lower semicontinuity,
	\[
	\norm{\nabla m(t)}{L^2(\mathscr D)}^2
	\leq
	e^{-4\theta t}
	\norm{\nabla m_0}{L^2(\mathscr D)}^2
	+
	C\bigl(1-e^{-4\theta t}\bigr).
	\]
	Hence, we obtain
	\[
	\limsup_{t\downarrow0}
	\norm{\nabla m(t)}{L^2(\mathscr D)}
	\leq
	\norm{\nabla m_0}{L^2(\mathscr D)}.
	\]
	On the other hand,
	\eqref{equ:epsilon strong L2} gives
	$m(t)\to m_0$ strongly in $L^2(\mathscr D)$ as $t\downarrow0$.
	The uniform $H^1$-bound then implies $m(t) \rightharpoonup m_0$ in $H^1(\mathscr D)$
	as $t\downarrow 0$. Therefore,
	\[
	\norm{\nabla m_0}{L^2(\mathscr D)}
	\leq
	\liminf_{t\downarrow0}
	\norm{\nabla m(t)}{L^2(\mathscr D)}.
	\]
	Consequently, $m(t) \to m_0$ strongly in $H^1(\mathscr D)$ as $t \downarrow 0$.
	Note that we have already shown that $m \in H^1(\tau,T;L^2(\mathscr D)) \cap L^2(\tau,T;H^2(\mathscr D))$. The standard Lions--Magenes interpolation result gives $m\in C([\tau,T];H^1(\mathscr D))$. Combining this with the
	continuity at $t=0$ yields $m\in C([0,T];H^1(\mathscr D))$, as required.
	
	Finally, we prove property {\rm (iv)}. Fix $0<s<t$. By \eqref{equ:positive time separation final}, applied with
	$\tau=s$, there exists $\rho_s\in(0,1)$ such that $|m(r,x)|\leq \rho_s$ for a.e. $x\in\mathscr D$ and every $r\in[s,t]$.
	Hence $m$ takes values a.e. in the compact set $\overline{B_{\rho_s}(0)}\subset B_1(0)$. Since $G\in C^1(B_1(0);\mathbb R^3)$, it is Lipschitz on $\overline{B_{\rho_s}(0)}$. By the regularity $m \in H^1(s,t;L^2(\mathscr D)) \cap L^2(s,t;H^2(\mathscr D))$, we can apply the chain rule, symmetry of $\mathsf B$, and the stray-field identities to obtain, for a.e. $r\in(s,t)$,
	\begin{align*}
		\frac{\mathrm d}{\dr}\mathcal E_\theta(m(r))
		&=
		-\alpha
		\inpro{\Delta m(r)}{\partial_t m(r)}
		-
		\inpro{\mathsf Bm(r)}{\partial_t m(r)}
		-
		\inpro{H_{\mathrm d}[m(r)]}{\partial_t m(r)}
		+
		\theta
		\inpro{G(m(r))}{\partial_t m(r)}
		\\
		&=
		-\inpro{h(m(r))}{\partial_t m(r)}.
	\end{align*}
	Taking the $L^2(\mathscr D)$ inner product of
	\eqref{equ:thermodynamic LLB} with $h(m)$ gives
	\begin{align*}
		\inpro{h(m(r))}{\partial_t m(r)}
		&=
		\gamma
		\inpro{h(m(r))}{m(r)\times h(m(r))}
		+
		\inpro{h(m(r))}{\mathbb M(m(r))h(m(r))}
		\\
		&=
		\inpro{\mathbb M(m(r))h(m(r))}{h(m(r))}.
	\end{align*}
	Hence,
	\[
	\frac{\mathrm d}{\dr}\mathcal E_\theta(m(r))
	+
	\int_{\mathscr D}
	h(m(r))\cdot\mathbb M(m(r))h(m(r))\dx
	=
	0
	\]
	for a.e. $r\in(s,t)$. Since
	$r\mapsto\mathcal E_\theta(m(r))$ is absolutely continuous on
	$[s,t]$, integration in time yields
	\[
	\mathcal E_\theta(m(t))
	+
	\int_s^t
	\int_{\mathscr D}
	h(m(r))\cdot\mathbb M(m(r))h(m(r))\dx\dr
	=
	\mathcal E_\theta(m(s)),
	\]
	which proves \eqref{equ:positive time energy identity}.
	
	Now, we show \eqref{equ:initial energy blowup}.
	Suppose that $\bigl| \{x\in\mathscr D:|m_0(x)|=1\} \bigr|>0$.
	If \eqref{equ:initial energy blowup} were false, there would exist
	$t_j\downarrow0$ such that
	$\mathcal E_\theta(m(t_j))$ remains bounded. Since
	$m(t_j)\to m_0$ strongly in $H^1(\mathscr D)$, after passing to a
	subsequence, we have $m(t_j,x)\to m_0(x)$ for a.e. $x\in\mathscr D$. Let $\Phi$ be the function defined in \eqref{equ:inverse Langevin potential}.
	Noting that $\Phi\geq0$ and $\Phi(s)\to+\infty$ as $s\uparrow1$, we obtain, by Fatou's
	lemma,
	\begin{equation}\label{equ:lim inf int Phi}
	\liminf_{j\to\infty}
	\int_{\mathscr D}
	\Phi(|m(t_j)|)\dx
	=
	+\infty.
	\end{equation}
	The exchange term in $\mathcal E_\theta(m(t_j))$ converges to a finite limit by the strong $H^1$-convergence, while
	the anisotropy and stray-field terms converge by the strong
	$L^2$-convergence and the $L^2$-continuity of
	$H_{\mathrm d}$. Therefore, $\mathcal E_\theta(m(t_j)) \to +\infty$,
	in contradiction with~\eqref{equ:lim inf int Phi}. This shows
	\eqref{equ:initial energy blowup}. 
	
	The proof is now complete.
\end{proof}


\section{Low-temperature limit and long-time behaviour}
\label{sec:low temperature and long time}

Having established the global weak theory and the positive-time separation from saturation, we now investigate two  asymptotic regimes of the full nonlinear LLB dynamics. The first concerns the low-temperature limit $\theta\to 0$, in which the longitudinal relaxation is expected to disappear and the dynamics approach the
classical Landau--Lifshitz--Gilbert equation on the unit sphere.

The second concerns the asymptotic behaviour for fixed temperature as $t\to\infty$. Positive-time separation makes the mobility uniformly coercive and allows us to exploit the thermodynamic formulation and the free-energy dissipation. Combined with a \L{}ojasiewicz--Simon inequality, this yields convergence of every global solution to a single stationary state. We then identify a spectral Curie-type threshold associated with the Hessian of the free energy and show that, above this threshold, the zero equilibrium is globally exponentially attracting.

Throughout this section, global weak solutions are understood to satisfy
the additional properties established in Theorem~\ref{thm:global existence full nonlinear LLB}.

\subsection{Convergence to the Landau--Lifshitz--Gilbert equation}
\label{subsec:low temperature limit}

For $v\in H^2(\mathscr D)$ satisfying $\partial_\nu v=0$, recall from \eqref{equ:H MFA mic} that $H_{\mathrm{MFA}}$ is a self-adjoint operator given by
\[
H_{\mathrm{MFA}}(v)
=
\alpha\Delta v+\mathsf Bv+H_{\mathrm d}[v],
\]
where $\mathsf B=JI+\mathsf A$.
Since $\mathbb M(n)=I-n\otimes n$ whenever $|n|=1$, the corresponding
normalised Landau--Lifshitz equation is
\begin{equation}
	\partial_t n
	=
	\gamma n\times H_{\mathrm{MFA}}(n)
	+
	\mathbb M(n)H_{\mathrm{MFA}}(n),
	\qquad
	|n|=1.
	\label{equ:limiting LLG}
\end{equation}
Equivalently,
\[
\partial_t n
=
\gamma n\times H_{\mathrm{MFA}}(n)
-
n\times
\bigl(
n\times H_{\mathrm{MFA}}(n)
\bigr).
\]
If $\mathsf B=JI+\mathsf A$, the isotropic contribution $Jn$ disappears
from \eqref{equ:limiting LLG}, while the anisotropic contribution
$\mathsf A n$ remains. This is the usual form of the Landau--Lifshitz--Gilbert equation~\cite{AloSoy92, Pro01}.

We recall that the extension of $\mathbb M$ to $\overline{B_1(0)}$ belongs to
$C^{1,1}(\overline{B_1(0)};\mathbb R^{3\times3})$; see Lemma~\ref{lem:properties M}. Fix $\kappa\coloneq 1+\norm{\mathsf B}{}$, and define, for $v\in H^1(\mathscr D)$,
\begin{equation}
	\mathscr R(v)
	\coloneq
	\frac{\kappa}{2}\norm{v}{L^2(\mathscr D)}^2
	+
	\frac{\alpha}{2}\norm{\nabla v}{L^2(\mathscr D)}^2
	-
	\frac12\inpro{\mathsf Bv}{v}
	+
	\frac12
	\norm{H_{\mathrm d}[v]}{L^2(\mathbb R^3)}^2.
	\label{equ:relative micromagnetic energy}
\end{equation}
Since $\kappa-\norm{\mathsf B}{}\geq1$, there exist constants
$c,C>0$ such that
\begin{equation}
	c\norm{v}{H^1(\mathscr D)}^2
	\leq
	\mathscr R(v)
	\leq
	C\norm{v}{H^1(\mathscr D)}^2.
	\label{equ:relative energy equivalence}
\end{equation}

The following theorem shows that, in effect, the Landau--Lifshitz equation \eqref{equ:limiting LLG} is an approximation of the full nonlinear LLB equation at low temperatures ($\theta\to 0^+$). We also provide a quantitative approximation rate.

\begin{theorem}[Low-temperature weak--strong convergence]
	\label{thm:LLB to LLG}
	Let $T>0$, and let $n$ be a strong solution of the Landau--Lifshitz equation
	\eqref{equ:limiting LLG} on $[0,T]$ such that
	\begin{equation}\label{equ:LLG strong regularity}
		\begin{split}
		n
		&\in
		C([0,T];H^1(\mathscr D))
		\cap
		H^1(0,T;L^2(\mathscr D))
		\cap
		L^2(0,T;H^2(\mathscr D)),
		\\
		n
		&\in
		L^\infty(0,T;W^{1,\infty}(\mathscr D)),
		\qquad
		H_{\mathrm{MFA}}(n)
		\in
		L^2(0,T;W^{1,\infty}(\mathscr D)),
		\end{split}
	\end{equation}
	with initial datum $n(0)=n_0$ satisfying $|n_0|=1$ a.e. in $\mathscr D$.
	
	For $\theta\in(0,1]$, let $m_\theta$ be a global weak solution of the
	full nonlinear LLB equation with initial datum $m_\theta(0)= m_{\theta,0}$, given by Theorem~\ref{thm:global existence full nonlinear LLB}. Set
	$e_\theta \coloneq m_\theta-n$.
	Then, there exists $C_T>0$, independent of $\theta$, such that for
	every $t\in[0,T]$,
	\begin{align}
		\mathscr R(e_\theta(t))
		+
		\frac12
		\int_0^t
		\inpro{
			\mathbb M(m_\theta)
			H_{\mathrm{MFA}}(e_\theta)}
		{H_{\mathrm{MFA}}(e_\theta)}
		\ds
		\leq
		C_T
		\left(
		\mathscr R(e_\theta(0))
		+
		\theta^2
		\right).
		\label{equ:low temperature relative estimate}
	\end{align}
	In particular, we have
	\begin{equation}
		\sup_{t\in[0,T]}
		\norm{m_\theta(t)-n(t)}{H^1(\mathscr D)}
		\leq
		C_T
		\left(
		\norm{m_{\theta}(0)-n(0)}{H^1(\mathscr D)}
		+
		\theta
		\right).
		\label{equ:low temperature H1 estimate}
	\end{equation}
	Furthermore,
	\begin{equation} 
		\sup_{t\in[0,T]} \norm{1-|m_\theta(t)|^2}{L^2(\mathscr D)} 
		\leq 
		C_T \left( \norm{m_{\theta}(0)-n(0)}{H^1(\mathscr D)} + \theta \right). 
		\label{equ:low temperature length defect} 
	\end{equation}
\end{theorem}

\begin{proof}
	Le $e_\theta = m_\theta-n$.
	For every $\delta>0$, the solution $m_\theta$ on $[\delta,T]$ is strictly separated from saturation and has the positive-time
	regularity given by Definition~\ref{def:global weak solution}. Hence, all calculations below are legitimate on $[\delta,T]$.
	
	Since $H_{\mathrm{MFA}}$ is linear, subtracting
	\eqref{equ:limiting LLG} from the closed-ball formulation of the
	LLB equation gives
	\begin{align}
		\partial_t e_\theta
		={}&
		\gamma m_\theta\times
		H_{\mathrm{MFA}}(e_\theta)
		+
		\mathbb M(m_\theta)
		H_{\mathrm{MFA}}(e_\theta)
		\nonumber\\
		&+
		\gamma e_\theta\times H_{\mathrm{MFA}}(n)
		+
		\bigl(
		\mathbb M(m_\theta)-\mathbb M(n)
		\bigr)
		H_{\mathrm{MFA}}(n)
		-
		2\theta m_\theta.
		\label{equ:low temperature difference equation}
	\end{align}
	Here, we used the fact that $\mathbb M(n)=I-n\otimes n$.
	
	We first record some preliminary estimates that will be used in the relative-energy identity. If $v\in W^{1,\infty}(\mathscr D)$, integration by parts, \eqref{equ:relative energy equivalence}, and \eqref{equ:stray field Lp} give
	\begin{align}
		\bigl|
		\inpro{H_{\mathrm{MFA}}(e_\theta)}
		{e_\theta\times v}
		\bigr|
		&\leq
		\alpha
		\sum_{i=1}^3
		\bigl|
		\inpro{\partial_i e_\theta}
		{e_\theta\times\partial_i v}
		\bigr|
		+
		C\norm{v}{L^\infty(\mathscr D)}
		\norm{e_\theta}{L^2(\mathscr D)}^2
		\nonumber\\
		&\leq
		C
		\bigl(
		1+\norm{v}{W^{1,\infty}(\mathscr D)}
		\bigr)
		\mathscr R(e_\theta).
		\label{equ:cross relative estimate}
	\end{align}
	
	Since $\mathbb M\in C^{1,1}(\overline{B_1(0)};\mathbb R^{3\times3})$,
	the Sobolev chain rule gives, for $i=1,2,3$,
	\begin{align*}
		\partial_i
		\bigl(
		\mathbb M(m_\theta)-\mathbb M(n)
		\bigr)
		&=
		D\mathbb M(m_\theta)\partial_i m_\theta
		-
		D\mathbb M(n)\partial_i n
		\\
		&=
		D\mathbb M(m_\theta)\partial_i e_\theta
		+
		\bigl(
		D\mathbb M(m_\theta)-D\mathbb M(n)
		\bigr)\partial_i n.
	\end{align*}
	Using the boundedness of $D\mathbb M$ and the Lipschitz continuity of
	$D\mathbb M$ on $\overline{B_1(0)}$, we obtain
	\begin{align*}
		\norm{
			\mathbb M(m_\theta)-\mathbb M(n)
		}{L^2(\mathscr D)}
		&\leq
		C\norm{e_\theta}{L^2(\mathscr D)},
		\\
		\norm{
			\nabla\bigl(
			\mathbb M(m_\theta)-\mathbb M(n)
			\bigr)
		}{L^2(\mathscr D)}
		&\leq
		C
		\bigl(
		\norm{\nabla e_\theta}{L^2(\mathscr D)}
		+
		\norm{\nabla n}{L^\infty(\mathscr D)}
		\norm{e_\theta}{L^2(\mathscr D)}
		\bigr).
	\end{align*}
	Consequently, we have
	\begin{align}
		&
		\norm{
			\bigl(
			\mathbb M(m_\theta)-\mathbb M(n)
			\bigr)
			H_{\mathrm{MFA}}(n)
		}{H^1(\mathscr D)}
		\leq
		C
		\bigl(
		1+\norm{\nabla n}{L^\infty(\mathscr D)}
		\bigr)
		\norm{H_{\mathrm{MFA}}(n)}
		{W^{1,\infty}(\mathscr D)}
		\norm{e_\theta}{H^1(\mathscr D)}.
		\label{equ:mobility difference H1 estimate}
	\end{align}
	
	By the symmetry of $\mathsf B$ and of the stray-field operator,
	\[
	\frac{\mathrm d}{\dt}\mathscr R(e_\theta)
	=
	\inpro{
		\kappa e_\theta
		-
		H_{\mathrm{MFA}}(e_\theta)}
	{\partial_t e_\theta}
	\]
	for a.e. $t\in(\delta,T)$. Substituting
	\eqref{equ:low temperature difference equation}, we obtain
	\begin{equation}
	\begin{split}
		\frac{\mathrm d}{\dt}\mathscr R(e_\theta)
		={}&
		\underbrace{
			\gamma
			\inpro{
				\kappa e_\theta-H_{\mathrm{MFA}}(e_\theta)}
			{m_\theta\times H_{\mathrm{MFA}}(e_\theta)}
		}_{\text{principal precession}}
		\\
		&+
		\underbrace{
			\inpro{
				\kappa e_\theta-H_{\mathrm{MFA}}(e_\theta)}
			{\mathbb M(m_\theta)H_{\mathrm{MFA}}(e_\theta)}
		}_{\text{principal damping}}
		\\
		&+
		\underbrace{
			\gamma
			\inpro{
				\kappa e_\theta-H_{\mathrm{MFA}}(e_\theta)}
			{e_\theta\times H_{\mathrm{MFA}}(n)}
		}_{\text{lower-order precession}}
		\\
		&+
		\underbrace{
			\inpro{
				\kappa e_\theta-H_{\mathrm{MFA}}(e_\theta)}
			{
				\bigl(
				\mathbb M(m_\theta)-\mathbb M(n)
				\bigr)
				H_{\mathrm{MFA}}(n)}
		}_{\text{mobility difference}}
		\\
		&-
		\underbrace{
			2\theta
			\inpro{
				\kappa e_\theta-H_{\mathrm{MFA}}(e_\theta)}
			{m_\theta}
		}_{\text{thermal forcing}}.
		\label{equ:relative energy decomposition}
	\end{split}
	\end{equation}
	We estimate these five terms in the order in which they appear.
	
	For the principal precession term, we note that  $e_\theta\times m_\theta= e_\theta\times n$.
	Hence, \eqref{equ:cross relative estimate}, applied with $v=n$,
	gives
	\begin{align}
		\left|
		\gamma
		\inpro{
			\kappa e_\theta-H_{\mathrm{MFA}}(e_\theta)}
		{m_\theta\times H_{\mathrm{MFA}}(e_\theta)}
		\right|
		&=
		|\gamma|\kappa
		\left|
		\inpro{H_{\mathrm{MFA}}(e_\theta)}
		{e_\theta\times n}
		\right|
		\leq
		C\mathscr R(e_\theta).
		\label{equ:relative principal precession}
	\end{align}
	
	We estimate the principal damping term next.
	Note that since $0\leq\mathbb M(m_\theta)\leq I$, we have $\mathbb M(m_\theta)^2 \leq \mathbb M(m_\theta)$. By \eqref{equ:mobility square estimate}, we infer that
	\begin{align}
		&
		\inpro{
			\kappa e_\theta-H_{\mathrm{MFA}}(e_\theta)}
		{\mathbb M(m_\theta)H_{\mathrm{MFA}}(e_\theta)}
		\nonumber\\
		&\quad
		=
		\kappa
		\inpro{e_\theta}
		{\mathbb M(m_\theta)H_{\mathrm{MFA}}(e_\theta)}
		-
		\inpro{
			\mathbb M(m_\theta)
			H_{\mathrm{MFA}}(e_\theta)}
		{H_{\mathrm{MFA}}(e_\theta)}
		\nonumber\\
		&\quad
		\leq
		\kappa
		\norm{e_\theta}{L^2(\mathscr D)}
		\bigl|\inpro{
			\mathbb M(m_\theta)
			H_{\mathrm{MFA}}(e_\theta)}
		{H_{\mathrm{MFA}}(e_\theta)}\bigr|^{1/2}
		-
		\inpro{
			\mathbb M(m_\theta)
			H_{\mathrm{MFA}}(e_\theta)}
		{H_{\mathrm{MFA}}(e_\theta)}
		\nonumber\\
		&\quad
		\leq
		-\frac12 \inpro{
			\mathbb M(m_\theta)
			H_{\mathrm{MFA}}(e_\theta)}
		{H_{\mathrm{MFA}}(e_\theta)}
		+
		C\mathscr R(e_\theta),
		\label{equ:relative principal damping}
	\end{align}
	where in the last step we used Young's inequality and \eqref{equ:relative energy equivalence}.
	
	For the lower-order precession term, the part containing
	$\kappa e_\theta$ vanishes. Hence
	\eqref{equ:cross relative estimate}, applied with
	$v=H_{\mathrm{MFA}}(n)$, yields
	\begin{align}
		\left|
		\gamma
		\inpro{
			\kappa e_\theta-H_{\mathrm{MFA}}(e_\theta)}
		{e_\theta\times H_{\mathrm{MFA}}(n)}
		\right|
		\leq
		C
		\bigl(
		1+
		\norm{H_{\mathrm{MFA}}(n)}
		{W^{1,\infty}(\mathscr D)}
		\bigr)
		\mathscr R(e_\theta).
		\label{equ:relative lower precession}
	\end{align}
	
	For the mobility-difference term, the variational definition of
	$H_{\mathrm{MFA}}$ gives
	\[
	\left|
	\inpro{
		\kappa e_\theta-H_{\mathrm{MFA}}(e_\theta)}
	{v}
	\right|
	\leq
	C
	\norm{e_\theta}{H^1(\mathscr D)}
	\norm{v}{H^1(\mathscr D)}.
	\]
	Consequently, using \eqref{equ:mobility difference H1 estimate} and
	\eqref{equ:relative energy equivalence}, we have
	\begin{align}
		\left|
		\inpro{
			\kappa e_\theta-H_{\mathrm{MFA}}(e_\theta)}
		{
			\bigl(
			\mathbb M(m_\theta)-\mathbb M(n)
			\bigr)
			H_{\mathrm{MFA}}(n)}
		\right|
		&\leq
		C \norm{e_\theta}{H^1(\mathscr D)}
		\norm{\bigl(\mathbb M(m_\theta)-\mathbb M(n)
			\bigr)
			H_{\mathrm{MFA}}(n)}{H^1(\mathscr{D})}
		\nonumber\\
		&\leq
		C
		\bigl(
		1+\norm{\nabla n}{L^\infty(\mathscr D)}
		\bigr)
		\norm{H_{\mathrm{MFA}}(n)}
		{W^{1,\infty}(\mathscr D)}
		\norm{e_\theta}{H^1(\mathscr D)}^2
		\nonumber\\
		&\leq
		C
		\bigl(
		1+
		\norm{H_{\mathrm{MFA}}(n)}
		{W^{1,\infty}(\mathscr D)}^2
		\bigr)
		\mathscr R(e_\theta)
		\label{equ:relative mobility difference}
	\end{align}
	
	Finally, we estimate the thermal forcing term. Noting the identity
	\[
	\inpro{
		\kappa e_\theta-H_{\mathrm{MFA}}(e_\theta)}
	{e_\theta}
	=
	\kappa\norm{e_\theta}{L^2}^2
	+
	\alpha\norm{\nabla e_\theta}{L^2}^2
	-
	\inpro{\mathsf Be_\theta}{e_\theta}
	+
	\norm{H_{\mathrm d}[e_\theta]}{L^2(\mathbb R^3)}^2
	=
	2\mathscr R(e_\theta),
	\]
	we obtain, by using $m_\theta=e_\theta+n$ and the self-adjointness of $H_{\mathrm{MFA}}$,
	\begin{align*}
		\inpro{
			\kappa e_\theta-H_{\mathrm{MFA}}(e_\theta)}
		{m_\theta}
		&=
		\inpro{
			\kappa e_\theta-H_{\mathrm{MFA}}(e_\theta)}
		{e_\theta}
		+
		\inpro{
			\kappa e_\theta-H_{\mathrm{MFA}}(e_\theta)}
		{n}
		\\
		&=
		2\mathscr R(e_\theta)
		+
		\inpro{e_\theta}
		{\kappa n-H_{\mathrm{MFA}}(n)}.
	\end{align*}
	Therefore, since $|n|=1$, $\theta\leq1$, and
	$\norm{e_\theta}{L^2}^2\leq C\mathscr R(e_\theta)$ by
	\eqref{equ:relative energy equivalence}, we apply Young's inequality to obtain
	\begin{align}
		2\theta
		\left|
		\inpro{
			\kappa e_\theta-H_{\mathrm{MFA}}(e_\theta)}
		{m_\theta}
		\right|
		&\leq
		4\theta\mathscr R(e_\theta)
		+
		2\theta
		\norm{e_\theta}{L^2}
		\bigl(
		\kappa\norm{n}{L^2}
		+
		\norm{H_{\mathrm{MFA}}(n)}{L^2}
		\bigr)
		\nonumber\\
		&\leq
		C\mathscr R(e_\theta)
		+
		C\theta^2
		\bigl(
		1+
		\norm{H_{\mathrm{MFA}}(n)}{L^2}^2
		\bigr)
		\nonumber\\
		&\leq
		C
		\bigl(
		1+
		\norm{H_{\mathrm{MFA}}(n)}{L^2}^2
		\bigr)
		\mathscr R(e_\theta)
		+
		C\theta^2
		\bigl(
		1+
		\norm{H_{\mathrm{MFA}}(n)}{L^2}^2
		\bigr).
		\label{equ:relative thermal forcing}
	\end{align}
	
	Substituting
	\eqref{equ:relative principal precession}--\eqref{equ:relative thermal forcing}
	into \eqref{equ:relative energy decomposition}, we deduce that, for a.e. $t\in(\delta,T)$,
	\begin{align}
		\frac{\mathrm d}{\dt}\mathscr R(e_\theta)
		+
		\frac12 \inpro{
			\mathbb M(m_\theta)
			H_{\mathrm{MFA}}(e_\theta)}
		{H_{\mathrm{MFA}}(e_\theta)}
		&\leq
		C
		\bigl(
		1+
		\norm{H_{\mathrm{MFA}}(n)}
		{W^{1,\infty}(\mathscr D)}^2
		\bigr)
		\mathscr R(e_\theta)
		\nonumber\\
		&\qquad+
		C\theta^2
		\bigl(
		1+
		\norm{H_{\mathrm{MFA}}(n)}
		{W^{1,\infty}(\mathscr D)}^2
		\bigr)
		\label{equ:relative differential inequality}
	\end{align}
	Integrating \eqref{equ:relative differential inequality} over
	$(\delta,t)$ and applying Gronwall's lemma yields
	\begin{align}
		\mathscr R(e_\theta(t))
		+
		\frac12
		\int_\delta^t
		\inpro{
			\mathbb M(m_\theta)
			H_{\mathrm{MFA}}(e_\theta)}
		{H_{\mathrm{MFA}}(e_\theta)} \ds
		\leq
		C_T
		\left(
		\mathscr R(e_\theta(\delta))
		+
		\theta^2
		\right),
		\label{equ:gronwall R}
	\end{align}
	where $C_T$ is a constant independent of $\delta$.

	Since
	$m_\theta,n\in C([0,T];H^1(\mathscr D))$, we have $e_\theta(\delta) \to m_{\theta,0}-n_0$ strongly in $H^1(\mathscr D)$ as $\delta\downarrow 0$,
	and thus, by the continuity of $\mathscr R$ on $H^1(\mathscr D)$, we get
	\[
	\mathscr R(e_\theta(\delta))
	\to
	\mathscr R(m_{\theta,0}-n_0).
	\]
	Moreover, since $\mathbb M(m_\theta)$ is positive semi-definite,
	\[
	\inpro{
		\mathbb M(m_\theta)
		H_{\mathrm{MFA}}(e_\theta)}
	{H_{\mathrm{MFA}}(e_\theta)}
	\geq0
	\qquad\text{a.e. in }(0,T).
	\]
	Therefore, by the monotone convergence theorem, as $\delta\downarrow 0$
	\[
	\int_\delta^t
	\inpro{
		\mathbb M(m_\theta)
		H_{\mathrm{MFA}}(e_\theta)}
	{H_{\mathrm{MFA}}(e_\theta)}
	\ds
	\to
	\int_0^t
	\inpro{
		\mathbb M(m_\theta)
		H_{\mathrm{MFA}}(e_\theta)}
	{H_{\mathrm{MFA}}(e_\theta)}
	\ds.
	\]
	Passing to the limit in \eqref{equ:gronwall R} then shows \eqref{equ:low temperature relative estimate}. The equivalence
	\eqref{equ:relative energy equivalence} gives
	\eqref{equ:low temperature H1 estimate}.
	
	Finally, since $|n|=1$ and $|m_\theta|\leq1$,
	\[
	0
	\leq
	1-|m_\theta|^2
	=
	(1-|m_\theta|)(1+|m_\theta|)
	\leq
	2|m_\theta-n|.
	\]
	Hence, \eqref{equ:low temperature length defect} also follows from
	\eqref{equ:low temperature H1 estimate}.
\end{proof}

\begin{remark}
	The existence of arbitrarily regular solutions to the Landau--Lifshitz--Gilbert equation for initial data that is sufficiently close to a constant function is proven in~\cite{FeiTra17}. In particular, the regularity hypothesis \eqref{equ:LLG strong regularity} is not vacuous.
\end{remark}

\subsection{Convergence to a stationary state}
\label{subsec:convergence stationary state}

We begin by identifying the stationary states governing the long-time dynamics. Since every global solution is uniformly separated from saturation on intervals bounded away from $t=0$, the mobility $\mathbb M(m)$ is uniformly positive definite there. Consequently, a stationary solution of the thermodynamic formulation satisfies
\[
0
=
\gamma m\times h(m)
+
\mathbb M(m)h(m),
\]
and testing with $h(m)$ yields $h(m)=0$. Thus, the stationary states are precisely the critical points of the free energy in the open unit ball.
For $r\geq 2$, let
\[
H^r_n(\mathscr D)
\coloneq
\left\{
v\in H^r(\mathscr D;\mathbb R^3):
\partial_\nu v=0
\text{ on }\partial\mathscr D
\right\},
\]
and set
\begin{equation}
	\mathcal S_\theta
	\coloneq
	\left\{
	\bar m\in H^2_n(\mathscr D):
	\norm{\bar m}{L^\infty(\mathscr D)}<1,
	\quad
	h(\bar m)=0
	\right\}.
	\label{equ:stationary set}
\end{equation}
The key tool for proving convergence to a single element of $\mathcal S_\theta$ is the following local \L{}ojasiewicz--Simon inequality near each stationary state.

\begin{lemma}[\L{}ojasiewicz--Simon inequality]\label{lem:LLB Lojasiewicz Simon}
	Let $\mathcal E_\theta$ be the energy functional given in \eqref{equ:full nonlinear LLB energy}.
	Let $\bar m\in\mathcal S_\theta$, and choose
	$\rho\in \bigl(\norm{\bar m}{L^\infty(\mathscr D)}, 1\bigr)$. Then, there exist
	constants $\sigma\in\left(0,1/2\right]$, $C>0$, and $\eta>0$, such that
	\begin{equation}
		\left|
		\mathcal E_\theta(v)
		-
		\mathcal E_\theta(\bar m)
		\right|^{1-\sigma}
		\leq
		C\norm{h(v)}{L^2(\mathscr D)}
		\label{equ:LLB Lojasiewicz Simon}
	\end{equation}
	for every $v\in H^2_n(\mathscr D)$ satisfying $\norm{v}{L^\infty(\mathscr D)}\leq\rho$ and $\norm{v-\bar m}{H^1(\mathscr D)}< \eta$.
\end{lemma}

\begin{proof}
	Since $H^2(\mathscr D)\hookrightarrow L^\infty(\mathscr D)$ and
	$H^2(\mathscr D)$ is a Banach algebra in three dimensions, a sufficiently small $H^2$-neighbourhood of $\bar m$ consists of functions taking values in $\overline{B_\rho(0)}$. On this ball, $v\mapsto\Phi(|v|)$ and $v\mapsto G(v)$ are real analytic. Hence, the associated Nemytskii maps are real analytic on this neighbourhood, and so is $\mathcal E_\theta\colon H^2_n(\mathscr D)\to\mathbb R$.
	Its $L^2$-gradient is
	\[
	\nabla_{L^2} \mathcal E_\theta (v)
	=
	-h(v).
	\]
	For $\bar m\in H_n^2(\mathscr D)$, the Hessian is
	\[
	D^2\mathcal E_\theta(\bar m)[v,w]
	=
	\inpro{
		\left(
		-\alpha\Delta
		-\mathsf B
		-H_{\mathrm d}
		+\theta DG(\bar m)
		\right)v
	}{w},
	\]
	for $v\in H_n^2(\mathscr D)$ and $w\in L^2(\mathscr D)$. Equivalently,
	the associated self-adjoint $L^2$-operator is
	\[
	\mathscr L_{\bar m}
	\coloneq
	-\alpha\Delta
	-\mathsf B
	-H_{\mathrm d}
	+\theta DG(\bar m)
	\colon
	H_n^2(\mathscr D)\to L^2(\mathscr D).
	\]
	Moreover, it is a Fredholm operator of index zero: the Neumann
	Laplacian is Fredholm from $H^2_n(\mathscr D)$ to $L^2(\mathscr D)$,
	and all remaining terms are compact perturbations as maps from
	$H^2_n(\mathscr D)$ to $L^2(\mathscr D)$.
	
	The abstract \L{}ojasiewicz--Simon theorem~\cite{FeeMar20} then yields constants
	$\sigma\in(0,1/2]$, $C>0$, and $\delta>0$ such that
	\[
	\left|
	\mathcal E_\theta(v)
	-
	\mathcal E_\theta(\bar m)
	\right|^{1-\sigma}
	\leq
	C\norm{h(v)}{L^2(\mathscr D)}
	\]
	whenever
	$\norm{v-\bar m}{H^2(\mathscr D)}<\delta$.
	
	It remains to replace the $H^2$-neighbourhood by an $H^1$-neighbourhood
	for functions uniformly separated from the saturation boundary. Since
	$h(\bar m)=0$, we have
	\[
	\alpha\Delta(v-\bar m)
	=
	h(v)
	-
	\mathsf B(v-\bar m)
	-
	H_{\mathrm d}[v-\bar m]
	+
	\theta\bigl(G(v)-G(\bar m)\bigr).
	\]
	If $\norm{v}{L^\infty(\mathscr D)}\leq \rho$ and $\norm{\bar m}{L^\infty(\mathscr D)} \leq\rho <1$,
	then $G\colon\overline{B_\rho(0)}\to\mathbb R^3$ is Lipschitz, and thus
	\[
	\norm{G(v)-G(\bar m)}{L^2(\mathscr D)}
	\leq
	C_\rho
	\norm{v-\bar m}{L^2(\mathscr D)}.
	\]
	Therefore, the Neumann elliptic estimate and the $L^2$-boundedness of
	$H_{\mathrm d}$ yield
	\begin{equation}
		\norm{v-\bar m}{H^2(\mathscr D)}
		\leq
		C_\rho
		\left(
		\norm{h(v)}{L^2(\mathscr D)}
		+
		\norm{v-\bar m}{L^2(\mathscr D)}
		\right).
		\label{equ:elliptic LS extension}
	\end{equation}
	Now, choose $\eta>0$ sufficiently small, such that $C_\rho\eta< \delta/2$ and
	$\left|\mathcal E_\theta(v) - \mathcal E_\theta(\bar m) \right| \leq 1$ whenever
	$\norm{v-\bar m}{H^1(\mathscr D)}<\eta$. If $\norm{h(v)}{L^2(\mathscr D)} < \delta/(2C_\rho)$, then \eqref{equ:elliptic LS extension} gives
	$\norm{v-\bar m}{H^2}<\delta$, and the original \L{}ojasiewicz--Simon inequality applies. If instead $\norm{h(v)}{L^2(\mathscr D)} \geq \delta/(2C_\rho)$, then
	\[
	\left|
	\mathcal E_\theta(v)
	-
	\mathcal E_\theta(\bar m)
	\right|^{1-\sigma}
	\leq
	1
	\leq
	\frac{2C_\rho}{\delta}
	\norm{h(v)}{L^2(\mathscr D)}.
	\]
	This proves \eqref{equ:LLB Lojasiewicz Simon}.
\end{proof}

As a consequence of the above \L{}ojasiewicz--Simon inequality, we have the following result on quantitative convergence to an equilibrium.

\begin{theorem}[Convergence to a single equilibrium]
	\label{thm:LLB convergence equilibrium}
	Let $m$ be a global weak solution of the full nonlinear LLB equation given by Theorem~\ref{thm:global existence full nonlinear LLB} and let $\mathcal S_\theta$ be given by~\eqref{equ:stationary set}. Then, there
	exists $m_\infty\in\mathcal S_\theta$ such that
	\begin{equation}
		m(t)
		\to
		m_\infty
		\qquad
		\text{strongly in }H^1(\mathscr D)
		\quad\text{as }t\to\infty.
		\label{equ:LLB H1 convergence equilibrium}
	\end{equation}
	Moreover, there exists
	$\sigma\in(0,1/2]$ such that, for all sufficiently large $t$,
	\begin{equation}
		\mathcal E_\theta(m(t))
		-
		\mathcal E_\theta(m_\infty)
		\leq
		\begin{cases}
			Ce^{-ct},
			&\sigma=\dfrac12,
			\\[1ex]
			C(1+t)^{-\frac{1}{1-2\sigma}},
			&0<\sigma<\dfrac12,
		\end{cases}
		\label{equ:LLB energy convergence rate}
	\end{equation}
	and
	\begin{equation}
		\norm{m(t)-m_\infty}{L^2(\mathscr D)}
		\leq
		C
		\left[
		\mathcal E_\theta(m(t))
		-
		\mathcal E_\theta(m_\infty)
		\right]^\sigma.
		\label{equ:LLB L2 equilibrium rate}
	\end{equation}
\end{theorem}

\begin{proof}
	Apply \eqref{equ:positive time separation final} with $\tau=1$ and
	write $\rho\in(0,1)$ for the corresponding separation constant, so that $|m(t)|\leq\rho$ a.e. in $\mathscr D$ for every
	$t\geq1$. The energy is bounded from below and, by \eqref{equ:positive time energy identity}, $t\mapsto\mathcal E_\theta(m(t))$ is non-increasing on $[1,\infty)$. Hence, there exists
	$\mathcal E_\infty\in\mathbb R$ such that
	$\mathcal E_\theta(m(t)) \to \mathcal E_\infty$ as $t\to \infty$.
	
	Moreover, by \eqref{equ:mobility radial coercivity} and the
	monotonicity of $q$,
	\begin{align}
		2q(\rho)\norm{h(m(t))}{L^2(\mathscr D)}^2
		&\leq
		\inpro{\mathbb M(m(t))h(m(t))}{h(m(t))}.
		\label{equ:q rho h Mh}
	\end{align}
	Using again the energy identity \eqref{equ:positive time energy identity}, and the lower boundedness of $\mathcal E_\theta$, we have
	\begin{align}
		\int_1^\infty
		\norm{h(m(t))}{L^2(\mathscr D)}^2\dt
		&\leq
		\frac{1}{2q(\rho)}
		\int_1^\infty
		\inpro{\mathbb M(m(t))h(m(t))}{h(m(t))}\dt
		\nonumber\\
		&=
		\frac{
			\mathcal E_\theta(m(1))-\mathcal E_\infty
		}{
			2q(\rho)
		}
		<\infty,
		\label{equ:finite h dissipation}
	\end{align}
	where
	$\mathcal E_\infty=\lim_{t\to\infty}\mathcal E_\theta(m(t))$.
	
	Furthermore, noting~\eqref{equ:intro full nonlinear LLB}, we apply 
	\eqref{equ:mobility square estimate} and
	\eqref{equ:precession mobility control} give
	\begin{align*}
		\norm{\partial_t m(t)}{L^2(\mathscr D)}^2
		&=
		\gamma^2\norm{m(t)\times h(m(t))}{L^2(\mathscr D)}^2
		+
		\norm{\mathbb M(m(t))h(m(t))}{L^2(\mathscr D)}^2
		\\
		&\leq
		\left(\frac32\gamma^2+1\right)
		\inpro{\mathbb M(m(t))h(m(t))}{h(m(t))}.
	\end{align*}
	Consequently, we infer that
	\begin{equation}
		\int_1^\infty
		\norm{\partial_t m(t)}{L^2(\mathscr D)}^2\dt
		<\infty.
		\label{equ:finite time derivative dissipation}
	\end{equation}
	
	We next prove relative compactness of the trajectory in $H^1(\mathscr D)$. Let $t_j\to\infty$. From \eqref{equ:finite h dissipation}, we may choose
	$s_j\in[t_j,t_j+1]$ such that $\norm{h(m(s_j))}{L^2(\mathscr D)} \to 0$.
	The sequence $m(s_j)$ is bounded in $H^1(\mathscr D)$, so after
	passing to a subsequence,
	\[
	m(s_j)
	\rightharpoonup
	m_\ast
	\quad\text{weakly in }H^1(\mathscr D), \qquad \text{ and }
	\qquad
	m(s_j)
	\to
	m_\ast
	\quad\text{strongly in }L^2(\mathscr D).
	\]
	Since $G$ is Lipschitz on $\{|v|\leq\rho\}$, we deduce $G(m(s_j)) \to G(m_\ast)$ in $L^2(\mathscr D)$.
	The identity
	\[
	\alpha\Delta m(s_j)
	=
	h(m(s_j))
	-
	\mathsf Bm(s_j)
	-
	H_{\mathrm d}[m(s_j)]
	+
	\theta G(m(s_j))
	\]
	then shows that $\Delta m(s_j)$ converges strongly in
	$L^2(\mathscr D)$. By the Neumann elliptic estimate,
	\[
	m(s_j)
	\to
	m_\ast
	\quad
	\text{strongly in }H^2(\mathscr D).
	\]
	In particular, $h(m_\ast)=0$, thus $m_\ast\in\mathcal S_\theta$.
	
	Furthermore, by \eqref{equ:finite time derivative dissipation}, as $j\to\infty$,
	\begin{align*}
	\norm{m(t_j)-m(s_j)}{L^2(\mathscr D)}
	\leq
	\left(
	\int_{t_j}^{t_j+1}
	\norm{\partial_t m(t)}{L^2(\mathscr D)}^2\dt
	\right)^{1/2}
	\to0.
	\end{align*}
	Hence, $m(t_j)\to m_\ast$ strongly in $L^2(\mathscr D)$ and weakly
	in $H^1(\mathscr D)$. Since $\mathcal E_\theta(m(t_j)) \to \mathcal E_\infty=
	\mathcal E_\theta(m_\ast)$,
	and all terms in the energy except the exchange term converge by
	the strong $L^2$-convergence and the uniform separation, we obtain
	$\norm{\nabla m(t_j)}{L^2(\mathscr D)} \to \norm{\nabla m_\ast}{L^2(\mathscr D)}$.
	Thus
	\[
	m(t_j)
	\to
	m_\ast
	\qquad
	\text{strongly in }H^1(\mathscr D).
	\]
	The trajectory is therefore relatively compact in $H^1(\mathscr D)$,
	and every element of its $\omega$-limit set belongs to
	$\mathcal S_\theta$ and has energy $\mathcal E_\infty$.
	
	The $\omega$-limit set is compact in $H^1(\mathscr D)$. A finite
	collection of the neighbourhoods from
	Lemma~\ref{lem:LLB Lojasiewicz Simon} covers it. Taking the smallest
	\L{}ojasiewicz exponent and enlarging the constant if necessary, we
	obtain numbers $\sigma\in\left(0,1/2\right]$, $C>0$, $t_0>0$, such that
	\begin{equation}
		\left[
		\mathcal E_\theta(m(t))
		-
		\mathcal E_\infty
		\right]^{1-\sigma}
		\leq
		C\norm{h(m(t))}{L^2(\mathscr D)}
		\label{equ:uniform Lojasiewicz trajectory}
	\end{equation}
	for a.e. $t\geq t_0$.
	
	Now, set
	\[
	e(t)
	\coloneq
	\mathcal E_\theta(m(t))
	-
	\mathcal E_\infty.
	\]
	By the energy identity, the separation estimate \eqref{equ:mobility radial coercivity}, and \eqref{equ:uniform Lojasiewicz trajectory},
	\begin{align}
		-e'(t)
		=
		\inpro{\mathbb M(m(t))h(m(t))}{h(m(t))}
		\geq
		2q(\rho)
		\norm{h(m(t))}{L^2(\mathscr D)}^2
		\geq
		c\,e(t)^{2(1-\sigma)}
		\label{equ:energy LS ODE}
	\end{align}
	for a.e. $t\geq t_0$. Solving this scalar differential inequality
	gives \eqref{equ:LLB energy convergence rate}.
	
	Moreover, we obtain
	\begin{align*}
		-\frac{\mathrm d}{\dt}e(t)^\sigma
		&=
		\sigma e(t)^{\sigma-1}
		\inpro{\mathbb M(m(t))h(m(t))}{h(m(t))}
		\\
		&\geq
		c
		e(t)^{\sigma-1}
		\norm{h(m(t))}{L^2(\mathscr D)}^2
		\\
		&\geq
		c\norm{h(m(t))}{L^2(\mathscr D)}
		\geq
		c\norm{\partial_t m(t)}{L^2(\mathscr D)}.
	\end{align*}
	Here, the penultimate inequality follows from
	\eqref{equ:uniform Lojasiewicz trajectory}. Integrating over
	$(t,\infty)$ gives
	\[
	\int_t^\infty
	\norm{\partial_t m(s)}{L^2(\mathscr D)}\ds
	\leq
	Ce(t)^\sigma.
	\]
	Hence, $m(t)$ is Cauchy in $L^2(\mathscr D)$ and converges to a
	single limit $m_\infty\in\mathcal S_\theta$, with
	\eqref{equ:LLB L2 equilibrium rate}. Relative compactness in
	$H^1(\mathscr D)$ then upgrades the convergence to
	\eqref{equ:LLB H1 convergence equilibrium}.
\end{proof}

\subsection{A Curie-type threshold}
\label{subsec:Curie threshold}

In this section, we discuss how temperature affects the stability of the
paramagnetic equilibrium state $m=0$. Since
\begin{equation}\label{equ:Phi s expansion}
\Phi(s)
=
\frac32s^2+O(s^4)
\qquad
\text{as }s\to0,
\end{equation}
the quadratic part of the energy at the origin is
\[
D^2\mathcal E_\theta(0)[v,v]
=
\alpha\norm{\nabla v}{L^2(\mathscr D)}^2
-
\inpro{\mathsf Bv}{v}
+
\norm{H_{\mathrm d}[v]}{L^2(\mathbb R^3)}^2
+
3\theta\norm{v}{L^2(\mathscr D)}^2.
\]
This motivates introducing the spectral quantity
\begin{equation}
	\Lambda_\ast
	\coloneq
	\sup_{0\neq v\in H^1(\mathscr D)}
	\frac{
		\inpro{\mathsf Bv}{v}
		-
		\alpha\norm{\nabla v}{L^2(\mathscr D)}^2
		-
		\norm{H_{\mathrm d}[v]}{L^2(\mathbb R^3)}^2
	}{
		\norm{v}{L^2(\mathscr D)}^2
	}.
	\label{equ:critical spectral value}
\end{equation}
The quantity $\Lambda_\ast$ is finite: $\Lambda_\ast\leq \norm{\mathsf B}{}$. Furthermore, by definition,
\[
D^2\mathcal E_\theta(0)[v,v]
\geq
(3\theta-\Lambda_\ast)
\norm{v}{L^2(\mathscr D)}^2.
\]
Hence, when $\Lambda_\ast>0$, the value $\theta_{\mathrm c}\coloneq \Lambda_\ast/3$
plays the role of a \emph{Curie-type threshold}: above it the paramagnetic state $m=0$ is a strict local energy minimizer with
an $L^2$-coercive Hessian, whereas below it the quadratic energy possesses an
unstable direction. If $\Lambda_\ast\leq0$, then the Hessian at
$m=0$ is coercive throughout all $\theta>0$.

We remark that this threshold should be distinguished from the bulk mean-field value $\theta_{\mathrm c}^{\mathrm{bulk}}=J/3$ introduced in
Section~\ref{subsec:relation reduced LLB}, which is used there only in the asymptotic reduction to the standard LLB model. Indeed, writing
$\mathsf B=JI+\mathsf A$, we have
\[
\Lambda_\ast
=
J+
\sup_{0\neq v\in H^1(\mathscr D)}
\frac{
	\inpro{\mathsf Av}{v}
	-
	\alpha\norm{\nabla v}{L^2(\mathscr D)}^2
	-
	\norm{H_{\mathrm d}[v]}{L^2(\mathbb R^3)}^2
}{
	\norm{v}{L^2(\mathscr D)}^2
}.
\]
Thus, $J/3$ is the Curie threshold of the homogeneous bulk mean-field
reference problem, whereas $\Lambda_\ast/3$ is the corresponding
sample-dependent stability threshold for the full micromagnetic PDE,
including exchange, anisotropy, stray-field effects, and the geometry
of $\mathscr D$. In the homogeneous bulk reduction, where these
micromagnetic corrections are neglected, $\Lambda_\ast$ reduces to
$J$ and the two thresholds coincide.

The following theorem formulates this more precisely. Above the critical temperature $\theta_{\mathrm c}$, the only stationary state is $0$ and every
solution converges to it, whereas below the threshold the energy admits
nonzero states of lower energy.

\begin{theorem}[Characterisation relative to the Curie threshold]
	\label{thm:Curie threshold behaviour}
	Let $m$ be a global weak solution of the full nonlinear LLB equation given by Theorem~\ref{thm:global existence full nonlinear LLB},
	and let $\theta_{\mathrm c}= \Lambda_\ast/3$.
	Then, the following assertions hold.
	\begin{enumerate}[(i)]
		\item If $\theta\geq \theta_{\mathrm c}$, then
		\[
		\mathcal S_\theta=\{0\}.
		\]
		Consequently, every global weak solution converges strongly in
		$H^1(\mathscr D)$ to zero.
		
		\item If $\theta > \theta_{\mathrm c}$, then the convergence is
		exponential: there exist $C,\omega>0$ such that
		\[
		\norm{m(t)}{H^1(\mathscr D)}
		\leq
		Ce^{-\omega t}
		\]
		for all sufficiently large $t$.
		
		\item If $\theta < \theta_{\mathrm c}$, then zero is \emph{not} a minimiser of $\mathcal E_\theta$. More precisely, there exists
		$v\in C^\infty(\overline{\mathscr D})$ satisfying
		$\partial_\nu v=0$ on $\partial\mathscr{D}$ and, for all sufficiently small
		$\delta>0$,
		\[
		\norm{\delta v}{L^\infty(\mathscr D)}<1,
		\qquad
		\mathcal E_\theta(\delta v)<0.
		\]
		In particular, any trajectory satisfying
		$\mathcal E_\theta(m(t_0))<0$ at some $t_0>0$ converges to a
		nonzero stationary state.
	\end{enumerate}
\end{theorem}

\begin{proof}
	Let $\bar m\in\mathcal S_\theta$. Testing $h(\bar m)=0$ by
	$\bar m$ gives
	\begin{align}
		0
		=
		-\alpha
		\norm{\nabla\bar m}{L^2(\mathscr D)}^2
		+
		\inpro{\mathsf B\bar m}{\bar m}
		-
		\norm{H_{\mathrm d}[\bar m]}{L^2(\mathbb R^3)}^2
		-
		\theta
		\int_{\mathscr D}
		\frac{|\bar m|^2}{q(|\bar m|)}\dx.
		\label{equ:stationary test Curie}
	\end{align}
	By the definition of $\Lambda_\ast$ and the monotonicity of $q$,
	\[
	\frac1{q(s)}
	\geq3
	\qquad
	\text{for }s\in[0,1),
	\]
	with strict inequality for $s>0$. Hence
	\begin{align*}
		0
		&\leq
		\Lambda_\ast
		\norm{\bar m}{L^2(\mathscr D)}^2
		-
		\theta
		\int_{\mathscr D}
		\frac{|\bar m|^2}{q(|\bar m|)}\dx
		\leq
		(\Lambda_\ast-3\theta)
		\norm{\bar m}{L^2(\mathscr D)}^2.
	\end{align*}
	If $3\theta>\Lambda_\ast$, this forces $\bar m=0$. If
	$3\theta=\Lambda_\ast$, then equality would require
	$q(|\bar m|)=1/3$ a.e. on $\{|\bar m|>0\}$, which is impossible
	because $q(s)<1/3$ for $s>0$. Thus, $\bar m=0$ also at the critical
	value. Assertion {\rm (i)} follows from Theorem~\ref{thm:LLB convergence equilibrium}.

		Next, assume that $\theta>\theta_{\mathrm c}$ and let $\mu\coloneq 3\theta-\Lambda_\ast>0$.
		We prove exponential convergence in three steps. First, the strict
		inequality $3\theta>\Lambda_\ast$ gives a gradient--energy estimate
		which controls the energy by $h(m)$. Second,
		the same condition implies that the energy controls the full
		$H^1(\mathscr D)$-norm. Finally, positive-time separation makes the
		mobility uniformly coercive, so that the energy identity yields a
		closed differential inequality for $\mathcal E_\theta(m(t))$.
		
		We begin with the first point. Recall that a differentiable
		functional $\mathcal F$ on a convex subset of a Hilbert space is
		called $\mu$-strongly convex with respect to $\norm{\cdot}{}$ if
		\[
		\mathcal F(v)
		\geq
		\mathcal F(u)
		+
		D\mathcal F(u)[v-u]
		+
		\frac{\mu}{2}\norm{v-u}{}^2
		\]
		for all admissible $u$ and $v$; see, e.g.,
		\cite[Chapter~10]{BauCom17}. For a twice-differentiable functional,
		a sufficient condition is the uniform Hessian bound
		$D^2\mathcal F(u)[z,z] \geq \mu\norm{z}{}^2$.
		For completeness, we verify this bound directly for
		$\mathcal E_\theta$.
		
		Set $\phi(v)\coloneq\Phi(|v|)$ for $|v|<1$.
		For $s=|v|\in(0,1)$, the Hessian of $\phi$ has tangential eigenvalue
		\[
		\frac{\mathcal L^{-1}(s)}{s}
		=
		\frac1{q(s)}
		\geq3,
		\]
		since $q(s)\leq q(0)=1/3$. Its radial eigenvalue is
		$(\mathcal L^{-1})'(s)$. Since $q$ is decreasing,
		\[
		0
		\geq
		q'(s)
		=
		\frac{
			\mathcal L^{-1}(s)
			-
			s(\mathcal L^{-1})'(s)
		}{
			(\mathcal L^{-1}(s))^2
		},
		\]
		and hence
		\[
		(\mathcal L^{-1})'(s)
		\geq
		\frac{\mathcal L^{-1}(s)}{s}
		\geq3.
		\]
		At $s=0$, both eigenvalues converge to $3$ by the expansion of
		$\mathcal L^{-1}$. Thus,
		\[
		D^2\phi(v)[z,z]
		\geq
		3|z|^2
		\qquad
		\text{for }|v|<1.
		\]
		Consequently,
		\begin{align*}
			D^2\mathcal E_\theta(v)[z,z]
			&=
			\alpha\norm{\nabla z}{L^2(\mathscr D)}^2
			-
			\inpro{\mathsf Bz}{z}
			+
			\norm{H_{\mathrm d}[z]}{L^2(\mathbb R^3)}^2
			+
			\theta
			\int_{\mathscr D}
			D^2\phi(v)[z,z]\dx
			\\
			&\geq
			\alpha\norm{\nabla z}{L^2(\mathscr D)}^2
			-
			\inpro{\mathsf Bz}{z}
			+
			\norm{H_{\mathrm d}[z]}{L^2(\mathbb R^3)}^2
			+
			3\theta\norm{z}{L^2(\mathscr D)}^2
			\\
			&\geq
			(3\theta-\Lambda_\ast)
			\norm{z}{L^2(\mathscr D)}^2
			=
			\mu\norm{z}{L^2(\mathscr D)}^2,
		\end{align*}
		where the last inequality follows from
		\eqref{equ:critical spectral value}. Hence
		$\mathcal E_\theta$ is $\mu$-strongly convex with respect to the
		$L^2(\mathscr D)$-norm.
		
		We now use this property only to relate the energy to its
		variational derivative. To make the argument explicit, for fixed
		$m$ consider
		\[
		f(r)
		\coloneq
		\mathcal E_\theta\bigl((1-r)m\bigr),
		\qquad 0\leq r\leq1.
		\]
		The segment $\{(1-r)m:0\leq r\leq1\}$ remains in the open unit ball
		whenever $m$ is uniformly separated from saturation, and
		\[
		f''(r)
		=
		D^2\mathcal E_\theta((1-r)m)[m,m]
		\geq
		\mu\norm{m}{L^2(\mathscr D)}^2.
		\]
		Taylor's formula with integral remainder gives
		\[
		f(1)
		\geq
		f(0)+f'(0)
		+
		\frac{\mu}{2}
		\norm{m}{L^2(\mathscr D)}^2.
		\]
		Since $\mathcal E_\theta(0)=0$ and $D\mathcal E_\theta(m)=-h(m)$,
		we have $f'(0)= D\mathcal E_\theta(m)[-m]= \inpro{h(m)}{m}$.
		Therefore,
		\begin{align*}
		0
		&\geq
		\mathcal E_\theta(m)
		+
		\inpro{h(m)}{m}
		+
		\frac{\mu}{2}
		\norm{m}{L^2(\mathscr D)}^2
		\\
		&\geq
		\mathcal E_\theta(m)
		-
		\frac{1}{2\mu}
		\norm{h(m)}{L^2(\mathscr D)}^2
		-
		\frac{\mu}{2}
		\norm{m}{L^2(\mathscr D)}^2
		+
		\frac{\mu}{2}
		\norm{m}{L^2(\mathscr D)}^2,
		\end{align*}
		where we used Young's inequality. We deduce that
		\begin{equation}
			\mathcal E_\theta(m)
			\leq
			\frac{1}{2\mu}
			\norm{h(m)}{L^2(\mathscr D)}^2.
			\label{equ:high temperature PL inequality}
		\end{equation}
		This is the estimate that will turn the energy dissipation into an
		exponential decay inequality.
		
		We next show that $\mathcal E_\theta$ controls the full
		$H^1(\mathscr D)$-norm. Since $\mathcal L^{-1}(s)\geq3s$,
		we have
		\[
		\Phi(s)
		=
		\int_0^s\mathcal L^{-1}(\sigma) \mathrm{d}\sigma
		\geq
		\frac32s^2,
		\qquad 0\leq s<1.
		\]
		Therefore,
		\begin{align*}
			2\mathcal E_\theta(m)
			&\geq
			\alpha\norm{\nabla m}{L^2(\mathscr D)}^2
			-
			\inpro{\mathsf Bm}{m}
			+
			\norm{H_{\mathrm d}[m]}{L^2(\mathbb R^3)}^2
			+
			3\theta\norm{m}{L^2(\mathscr D)}^2
			\\
			&\geq
			(3\theta-\Lambda_\ast)
			\norm{m}{L^2(\mathscr D)}^2
			=
			\mu\norm{m}{L^2(\mathscr D)}^2,
		\end{align*}
		and consequently,
		\begin{equation}
			\norm{m}{L^2(\mathscr D)}^2
			\leq
			\frac{2}{\mu}\mathcal E_\theta(m).
			\label{equ:high temperature L2 coercivity}
		\end{equation}
		On the other hand, after discarding the nonnegative stray-field and
		potential terms,
		\[
		2\mathcal E_\theta(m)
		\geq
		\alpha\norm{\nabla m}{L^2(\mathscr D)}^2
		-
		\norm{\mathsf B}{}
		\norm{m}{L^2(\mathscr D)}^2.
		\]
		Using \eqref{equ:high temperature L2 coercivity} and rearranging the terms, we obtain
		\[
		\mathcal E_\theta(m)
		\geq
		\frac{\alpha\mu}
		{2(\mu+\norm{\mathsf B}{})}
		\norm{\nabla m}{L^2(\mathscr D)}^2.
		\]
		Together with \eqref{equ:high temperature L2 coercivity}, this
		implies that there exists $c_\theta>0$ such that
		\begin{equation}
			\mathcal E_\theta(m)
			\geq
			c_\theta
			\norm{m}{H^1(\mathscr D)}^2.
			\label{equ:high temperature H1 coercivity}
		\end{equation}
		
		It remains to combine these two estimates with the dissipation law.
		Using the energy identity, estimate \eqref{equ:q rho h Mh}, and
		\eqref{equ:high temperature PL inequality}, we infer that, for
		a.e. $t\geq t_0$,
		\begin{align*}
			\frac{\mathrm d}{\dt}
			\mathcal E_\theta(m(t))
			&=
			-
			\inpro{
				\mathbb M(m(t))h(m(t))
			}{
				h(m(t))
			}
			\\
			&\leq
			-2q(\rho)
			\norm{h(m(t))}{L^2(\mathscr D)}^2
			\\
			&\leq
			-4\mu q(\rho)
			\mathcal E_\theta(m(t)).
		\end{align*}
		Gronwall's lemma then yields
		\[
		\mathcal E_\theta(m(t))
		\leq
		\mathcal E_\theta(m(t_0))
		\exp\bigl(
		-4\mu q(\rho)(t-t_0)
		\bigr),
		\qquad t\geq t_0.
		\]
		Finally, \eqref{equ:high temperature H1 coercivity} gives
		\[
		\norm{m(t)}{H^1(\mathscr D)}^2
		\leq
		C
		\exp\bigl(
		-4\mu q(\rho)t
		\bigr),
		\]
		after absorbing the fixed time $t_0$ and
		$\mathcal E_\theta(m(t_0))$ into the constant.
		This proves assertion {\rm (ii)}.
	
	Finally, suppose $3\theta<\Lambda_\ast$, i.e. $\theta< \theta_{\mathrm c}$. By definition of
	$\Lambda_\ast$ and density, there exists
	$v\in C^\infty(\overline{\mathscr D})$ with
	$\partial_\nu v=0$ on $\partial\mathscr{D}$ such that
	\[
	\inpro{\mathsf Bv}{v}
	-
	\alpha\norm{\nabla v}{L^2(\mathscr D)}^2
	-
	\norm{H_{\mathrm d}[v]}{L^2(\mathbb R^3)}^2
	>
	3\theta\norm{v}{L^2(\mathscr D)}^2.
	\]
	Using \eqref{equ:Phi s expansion}, we obtain
	\begin{align*}
		\mathcal E_\theta(\delta v)
		={}&
		\frac{\delta^2}{2}
		\left[
		\alpha\norm{\nabla v}{L^2(\mathscr D)}^2
		-
		\inpro{\mathsf Bv}{v}
		+
		\norm{H_{\mathrm d}[v]}{L^2(\mathbb R^3)}^2
		+
		3\theta\norm{v}{L^2(\mathscr D)}^2
		\right]
		+
		O(\delta^4).
	\end{align*}
	The coefficient of $\delta^2$ is negative. Hence
	$\mathcal E_\theta(\delta v)<0$ for all sufficiently small
	$\delta>0$, and $\delta$ may also be chosen so that
	$\norm{\delta v}{L^\infty}<1$.
	
	If a trajectory satisfies
	$\mathcal E_\theta(m(t_0))<0$, then its limiting equilibrium from
	Theorem~\ref{thm:LLB convergence equilibrium} satisfies
	\[
	\mathcal E_\theta(m_\infty)
	\leq
	\mathcal E_\theta(m(t_0))
	<0.
	\]
	Since $\mathcal E_\theta(0)=0$, one has $m_\infty\neq0$. This shows
	assertion {\rm (iii)}, thus completing the proof of the theorem.
\end{proof}

\section{Non-saturated case: strong solutions and uniqueness}
\label{subsec:strong separated solutions}

Although the preceding weak-solution theory was formulated for the physical three-dimensional problem, the arguments extend, with only minor dimensional modifications, to $d\in\{1,2\}$. In this section, we explicitly allow $\mathscr D\subset\mathbb R^d$, $d\in\{1,2,3\}$, while the magnetisation remains $\mathbb R^3$-valued. For $d=3$, $H_{\mathrm d}$ denotes the physical whole-space stray-field operator introduced before. For $d\in \{1,2\}$, we consider a corresponding lower-dimensional operator, still denoted by $H_{\mathrm d}$, which is linear, self-adjoint and non-positive on $L^2(\mathscr D;\mathbb R^3)$, and satisfies \eqref{equ:stray field Lp}--\eqref{equ:stray field W2r}.
Under these assumptions, the global weak existence, separation results, and long-time behaviour established in Section~\ref{sec:global existence} remain valid for $d\in \{1,2\}$; in fact, the relevant estimates are
more favourable in lower dimensions.

Here, we consider initial data uniformly separated from saturation.
In this case, the equation is initially uniformly parabolic, and the
separation estimates established in Section~\ref{sec:global existence} can be combined with quasilinear parabolic theory to obtain stronger well-posedness results.

For ease of presentation, define
\[
\mathbb A(v)z
\coloneq
\mathbb M(v)z+\gamma v\times z.
\]
The closed-ball equation \eqref{equ:closed ball LLB} in strong form can be written as
\begin{equation}
	\partial_t m
	-
	\alpha\mathbb A(m)\Delta m
	=
	\mathbb A(m)
	\bigl(\mathsf Bm+H_{\mathrm d}[m]\bigr)
	-
	2\theta m.
	\label{equ:strong closed ball equation}
\end{equation}
For every $\rho<1$, we have
\begin{equation}
	\operatorname{Re}
	\bigl(\overline z\cdot\mathbb A(v)z\bigr)
	=
	\overline z\cdot\mathbb M(v)z
	\geq
	2q(\rho)|z|^2,
	\qquad
	|v|\leq\rho,\quad z\in\mathbb C^3.
	\label{equ:separated strong ellipticity}
\end{equation}
Thus, the skew-symmetric precession term does not affect ellipticity.

\subsection{Existence of a unique maximal strong solution}
 
The following theorem concerns the existence of a unique maximal strong solution, which is global when $d\leq 2$. To state this theorem, we define the interpolation space 
\[
\mathcal X_s
\coloneq
\bigl(
H^{s-1}(\mathscr D),
H_n^{s+1}(\mathscr D)
\bigr)_{\frac12,\,2}.
\]
Away from the borderline exponent $s=3/2$, this space is concretely given by
\[
\mathcal X_s
=
\begin{cases}
	H^s(\mathscr D), & 1<s<3/2,\\
	H_n^s(\mathscr D), & 3/2<s\leq2,
\end{cases}
\]
while for $s=3/2$ we retain the interpolation-space notation.

\begin{theorem}[Strong solutions]
	\label{thm:strong separated wellposedness}
	Let $s\in (1,2]$ satisfy $s>d/2$.
	Suppose $m_0\in\mathcal X_s$, and there exists a constant $\rho_0>0$ such that $\norm{m_0}{L^\infty(\mathscr D)}\leq\rho_0<1$.
	Then, there exists a unique maximal strong solution of
	\eqref{equ:strong closed ball equation}, with maximal existence time
	$T_{\max}\in(0,\infty]$, such that, for every $T<T_{\max}$,
	\begin{equation}
		m\in
		C\bigl([0,T];\mathcal X_s\bigr)
		\cap
		H^1\bigl(0,T;H^{s-1}(\mathscr D)\bigr)
		\cap
		L^2\bigl(0,T;H_n^{s+1}(\mathscr D)\bigr).
		\label{equ:strong separated regularity}
	\end{equation}
	For each such $T$, there exists $\rho_T<1$ satisfying
	\[
	\norm{m(t)}{L^\infty(\mathscr D)}
	\leq\rho_T
	\qquad
	\text{for every }t\in[0,T].
	\]
	The solution depends locally Lipschitz continuously on its initial
	datum in the function spaces appearing in
	\eqref{equ:strong separated regularity}.
	
	If $d\leq2$, then $T_{\max}=\infty$. In particular, for every $T>0$,
	\begin{equation}
		\sup_{t\in[0,T]}\norm{m(t)}{H^s}^2
		+
		\int_0^T
		\left(
		\norm{\partial_t m(t)}{H^{s-1}}^2
		+
		\norm{m(t)}{H^{s+1}}^2
		\right)\dt
		\leq C_T.
		\label{equ:global strong separated estimate}
	\end{equation}
	Moreover, for every $0<\tau<T$,
	\begin{equation}
		m\in
		C\bigl([\tau,T];H_n^2(\mathscr D)\bigr)
		\cap
		H^1\bigl(\tau,T;H^1(\mathscr D)\bigr)
		\cap
		L^2\bigl(\tau,T;H_n^3(\mathscr D)\bigr),
		\label{equ:low dimensional strong improvement}
	\end{equation}
	and
	\begin{equation}
		\operatorname*{ess\,sup}_{t\in[\tau,T]}
		\norm{\partial_t m(t)}{L^2}^2
		+
		\int_\tau^T
		\norm{\nabla\partial_t m(t)}{L^2}^2\dt
		\leq C_{\tau,T}.
		\label{equ:low dimensional time derivative improvement}
	\end{equation}
	The constants may depend on the initial datum, the fixed
	coefficients, the operator bounds, and the indicated time intervals.
\end{theorem}

\begin{proof}
	We first apply local maximal regularity and then prove that, in
	dimensions one and two, a finite maximal time is impossible.
	The latter step uses an estimate for $\partial_t m$ to recover
	a continuous $H^2$-trace at the maximal time.
	
	For local existence proof, set $X_0=H^{s-1}(\mathscr D)$ and $X_1=H_n^{s+1}(\mathscr D)$.
	The Neumann interpolation identity gives $(X_0,X_1)_{\frac12,\,2}= \mathcal X_s$.
	Since $\mathcal X_s\hookrightarrow H^s(\mathscr D)$ and
	$s>d/2$, the Sobolev multiplication theorem yields
	\[
	\norm{fg}{H^{s-1}}
	\leq
	C\norm{f}{H^s}\norm{g}{H^{s-1}}.
	\]
	Therefore, the maps
	\[
	v\longmapsto
	-\alpha\mathbb A(v)\Delta
	\in\mathcal L(X_1,X_0)
	\quad
	\text{ and }
	\quad
	v\longmapsto
	\mathbb A(v)(\mathsf Bv+H_{\mathrm d}[v])-2\theta v
	\in X_0
	\]
	are locally smooth on the open set
	$\mathcal{U} \coloneq \left\{
	v\in\mathcal X_s:
	\norm{v}{L^\infty}<1
	\right\}$.
	
	We also verify maximal regularity for the frozen principal operator.
	For each fixed $v\in\mathcal U$, the principal symbol of
	$-\alpha\mathbb A(v)\Delta$ is $\alpha |\xi|^2\mathbb A(v(x))$.
	By \eqref{equ:separated strong ellipticity}, this symbol is uniformly
	normally elliptic on every subset of $\mathcal U$ that is uniformly
	separated from saturation. Together with the homogeneous Neumann
	boundary condition, standard maximal-regularity theory for normally
	elliptic parabolic systems thus shows that
	\[
	-\alpha\mathbb A(v)\Delta
	\colon
	H_n^{s+1}(\mathscr D)
	\to
	H^{s-1}(\mathscr D)
	\]
	has maximal $L^2$-regularity; see, e.g.,~\cite{Ama95,Ama05}.
	The maximal regularity theorem for quasilinear parabolic problems then gives the
	unique maximal solution and continuous dependence stated above;
	see~\cite[Theorems~2.1 and~3.1]{Ama05}.
	Since $H^s(\mathscr D)\hookrightarrow C(\overline{\mathscr D})$ and
	$m\in C([0,T];\mathcal X_s)\hookrightarrow
		C([0,T];H^s(\mathscr D))$, we have
	\[
	\norm{m(t)-m_0}{L^\infty(\mathscr D)}
	\to0
	\qquad\text{as }t\downarrow0.
	\]
	Hence, the strict initial separation
	$\norm{m_0}{L^\infty}\leq\rho_0<1$ persists on a sufficiently short
	time interval.
	Moreover, if $T_{\max}<\infty$ and
	$m(t)\to m_\ast$ in $\mathcal X_s$ as
	$t\uparrow T_{\max}$ with
	$\norm{m_\ast}{L^\infty}<1$, then $m_\ast$ belongs to the open phase
	space $\mathcal U$. Restarting the local problem from $m_\ast$ at
	$t=T_{\max}$ therefore extends the solution beyond $T_{\max}$,
	contradicting maximality.

	We now assume $d\leq2$ and suppose, for contradiction, that
	$T_{\max}<\infty$. Since $s>1$, we have
		$m_0\in\mathcal X_s\hookrightarrow H^1(\mathscr D)$.
	The finite-time $H^1$ and
	time-derivative \emph{a priori} estimates, as well as the
	positive-time separation established in Theorem~\ref{thm:global existence full nonlinear LLB},
	apply on every interval $[0,S]$ with $S<T_{\max}$, with constants that
	may be chosen uniformly for $S\leq T_{\max}$.
	Consequently, for every fixed $a\in(0,T_{\max})$, there exist
	$\rho_a\in(0,1)$ and $C_a>0$ such that
	\begin{equation}
		|m|\leq\rho_a<1
		\quad\text{on }[a,T_{\max})\times\mathscr D,
		\label{equ:strong separation near maximal time}
	\end{equation}
	and
	\begin{equation}
		\sup_{t<T_{\max}}\norm{m(t)}{H^1}^2
		+
		\int_a^{T_{\max}}
		\left(
		\norm{m(t)}{H^2}^2
		+
		\norm{\partial_t m(t)}{L^2}^2
		\right)\dt
		\leq C_a.
		\label{equ:strong continuation basic bounds}
	\end{equation}
	Here, the constants are independent of $S<T_{\max}$.
	The $H^2$-bound follows from the weighted dissipation and
	$q(|m|)\geq q(\rho_a)>0$ as in \eqref{equ:Lap bound by weight}.
	
	Now, set $w\coloneq \partial_t m$. Time differentiation of
	\eqref{equ:strong closed ball equation} gives
	\[
	\partial_t w
	=
	\alpha\mathbb A(m)\Delta w
	+
	D\mathbb A(m)[w]H_{\mathrm{MFA}}(m)
	+
	\mathbb A(m)(\mathsf Bw+H_{\mathrm d}[w])
	-
	2\theta w.
	\]
	Testing by $w$, integrating the principal term by parts, and using
	the skew-symmetry of the precessional coefficient, we obtain
	\begin{align}
		&\frac12\frac{\mathrm d}{\dt}\norm{w}{L^2}^2
		+
		\alpha\sum_{i=1}^d
		\inpro{\mathbb M(m)\partial_iw}{\partial_iw}
		+
		2\theta\norm{w}{L^2}^2
		\nonumber\\
		&\quad
		=
		-\alpha\sum_{i=1}^d
		\inpro{
			D\mathbb A(m)[\partial_i m]\partial_iw
		}{w}
		+
		\inpro{
			D\mathbb A(m)[w]H_{\mathrm{MFA}}(m)
		}{w}
		+
		\inpro{
			\mathbb A(m)(\mathsf Bw+H_{\mathrm d}[w])
		}{w}
		\nonumber\\
		&\quad\leq
		C_a
		\left(
		\norm{\nabla m}{L^4}
		\norm{w}{L^4}
		\norm{\nabla w}{L^2}
		+
		\norm{H_{\mathrm{MFA}}(m)}{L^2}
		\norm{w}{L^4}^2
		+
		\norm{w}{L^2}^2
		\right).
		\label{equ:strong differentiated estimate}
	\end{align}
	For $d\leq2$, the Gagliardo--Nirenberg inequality gives
	\[
	\norm{w}{L^4}^2
	\leq
	C\norm{w}{L^2}
	\bigl(\norm{\nabla w}{L^2}+\norm{w}{L^2}\bigr).
	\]
	Therefore, applying Young's inequality in
	\eqref{equ:strong differentiated estimate}, and using the ellipticity
	\eqref{equ:separated strong ellipticity}, we obtain
	\begin{equation}
		\frac{\mathrm d}{\dt}\norm{w}{L^2}^2
		+
		c_a\norm{\nabla w}{L^2}^2
		\leq
		C_a
		\left(
		1+\norm{\nabla m}{L^4}^4
		+\norm{H_{\mathrm{MFA}}(m)}{L^2}^2
		\right)
		\norm{w}{L^2}^2,
		\label{equ:strong time derivative Gronwall}
	\end{equation}
	where the coefficient on the right-hand side is integrable up to $T_{\max}$.
	Therefore, Gronwall's inequality yields, for every
	$b\in (a,T_{\max})$,
	\begin{equation}
		\operatorname*{ess\,sup}_{b<t<T_{\max}}
		\norm{w(t)}{L^2}^2
		+
		\int_b^{T_{\max}}\norm{\nabla w(t)}{L^2}^2\dt
		<\infty.
		\label{equ:strong time derivative bounded}
	\end{equation}
	
	We now recover two spatial bounds from this estimate.
	Since $\mathbb A(m)$ is uniformly invertible on
	$[b,T_{\max})$, equation
	\eqref{equ:strong closed ball equation} gives
	\begin{equation}
		\alpha\Delta m
		=
		\mathbb A(m)^{-1}(w+2\theta m)
		-
		\mathsf Bm-H_{\mathrm d}[m].
		\label{equ:strong elliptic recovery}
	\end{equation}
	Hence, we have
	\[
	\norm{m(t)}{H^2}
	\leq C_a\bigl(1+\norm{w(t)}{L^2}\bigr),
	\]
	so $m$ is bounded in $L^\infty(b,T_{\max};H_n^2)$.
	
	We next recover one additional spatial derivative. By the uniform
	$H^2$-bound for $m$ and the separation
	\eqref{equ:strong separation near maximal time}, the standard Sobolev
	composition estimate gives
	\[
	\sup_{b<t<T_{\max}}
	\norm{\mathbb A(m(t))^{-1}}{H^2}
	\leq C_a.
	\]
	Noting that the product of $H^2$ and $H^1$ functions lies in $H^1$,
	by the elliptic recovery identity
	\eqref{equ:strong elliptic recovery} and the $H^1$-boundedness
	of $H_{\mathrm d}$, we obtain
	\[
	\norm{\Delta m}{H^1}
	\leq
	C_a\bigl(1+\norm{w}{H^1}\bigr).
	\]
	Therefore, by \eqref{equ:strong time derivative bounded} and Neumann
	elliptic regularity,
	\[
	m\in L^2(b,T_{\max};H_n^3(\mathscr D)),
	\qquad
	\partial_t m=w\in L^2(b,T_{\max};H^1(\mathscr D)).
	\]
	The standard interpolation result then shows that $m\in C([b,T_{\max}];H_n^2(\mathscr D))$.
	Since $H_n^2(\mathscr D)\hookrightarrow\mathcal X_s$ for $s\in (1,2]$,
	in particular, $m(t)$ has a limit in $\mathcal X_s$ as
	$t\uparrow T_{\max}$. By
	\eqref{equ:strong separation near maximal time}, this limit remains
	strictly separated from saturation.
	Thus, we deduce $T_{\max}=\infty$.
	The maximal-regularity class gives
	\eqref{equ:global strong separated estimate} on every finite interval.
	Repeating the preceding estimates on $[\tau/2,T]$ proves
	\eqref{equ:low dimensional strong improvement} and
	\eqref{equ:low dimensional time derivative improvement}.
\end{proof}

\subsection{Exponential stability}

We next obtain an exponential stability result. This implies the existence of a unique global strong solution for small initial data when $d=3$.

\begin{theorem}[Exponential stability of non-degenerate equilibria]
	\label{thm:strong stable equilibrium}
	Let $d\leq3$ and let $s\in (1,2]$ satisfy $s>d/2$. Let $\bar m\in\mathcal S_\theta$, and
	define the operator 
	\[
	\mathscr L_{\bar m} \coloneq D^2\mathcal E_\theta(\bar m)= -\alpha\Delta-\mathsf B-H_{\mathrm d}
	+\theta DG(\bar m),
	\]
	which is the
	Hessian of the free energy at the equilibrium $\bar m$.
	Assume that there exists $\mu>0$ such that
	\begin{equation}
		\inpro{\mathscr L_{\bar m}z}{z}
		\geq
		\mu\norm{z}{L^2}^2
		\qquad
		\text{for every }z\in H_n^2(\mathscr D),
		\label{equ:strong equilibrium stability}
	\end{equation}
	Then, there exist $\eta,\omega,C>0$ such that, whenever $m_0\in\mathcal X_s$ and $\norm{m_0-\bar m}{\mathcal X_s}\leq\eta$,
	the corresponding solution is global, remains uniformly separated
	from saturation, and satisfies
	\begin{align}
		&\sup_{t\geq0}
		e^{2\omega t}\norm{m(t)-\bar m}{H^s}^2
		\nonumber\\
		&\qquad+
		\int_0^\infty e^{2\omega t}
		\left(
		\norm{\partial_t m(t)}{H^{s-1}}^2
		+
		\norm{m(t)-\bar m}{H^{s+1}}^2
		\right)\dt
		\leq
		C\norm{m_0-\bar m}{\mathcal X_s}^2.
		\label{equ:strong exponential stability}
	\end{align}
	In particular,
	\[
	\norm{m(t)-\bar m}{H^s}
	\leq
	Ce^{-\omega t}
	\norm{m_0-\bar m}{\mathcal X_s}.
	\]
	The constants depend on the fixed coefficients, the domain,
	$\bar m$, and the stability bound, but not on $m_0$ within the
	stated neighbourhood.
\end{theorem}

\begin{proof}
	We first verify stability of the linearised operator, then estimate
	the nonlinear remainder in the same maximal-regularity spaces used
	in Theorem~\ref{thm:strong separated wellposedness}.
	
	Since $h(\bar m)=0$ and $\bar m$ is strictly separated, the
	$H^2$-mapping property of $H_{\mathrm d}$ and elliptic regularity
	give $\bar m\in H^4(\mathscr D)$.
	In particular, all coefficients below have the regularity required
	for their realizations on $X_0=H^{s-1}(\mathscr D)$ and $X_1=H_n^{s+1}(\mathscr D)$.
	Set
	$\mathscr A_{\bar m}\coloneq \mathbb A(\bar m)\mathscr L_{\bar m}$.
	Since $h(\bar m)=0$ and
	$Dh(\bar m)=-\mathscr L_{\bar m}$, the linearisation of the LLB
	equation about $\bar m$ is $\partial_t u+\mathscr A_{\bar m}u=0$.
	
	Choose $\rho\in \bigl(\norm{\bar m}{L^\infty},1\bigr)$.
	The operator $\mathscr L_{\bar m}$ is self-adjoint and satisfies
	$\mathscr L_{\bar m}\geq\mu I$ on $L^2$.
	If
	\[
	\mathscr A_{\bar m}z=\lambda z,
	\qquad z\neq0,
	\]
	then, using the complexified $L^2$ inner product,
	\begin{align*}
		\operatorname{Re}\lambda\,
		\inpro{\mathscr L_{\bar m}z}{z}
		&=
		\operatorname{Re}
		\inpro{
			\mathbb A(\bar m)\mathscr L_{\bar m}z
		}{\mathscr L_{\bar m}z}
		\\
		&=
		\inpro{
			\mathbb M(\bar m)\mathscr L_{\bar m}z
		}{\mathscr L_{\bar m}z}
		\\
		&\geq
		2q(\rho)\norm{\mathscr L_{\bar m}z}{L^2}^2
		\\
		&\geq
		2\mu q(\rho)
		\inpro{\mathscr L_{\bar m}z}{z}.
	\end{align*}
	Therefore, $\operatorname{Re}\lambda\geq2\mu q(\rho)>0$.
	Moreover, $\mathscr A_{\bar m}$ is normally elliptic with compact
	resolvent, so its spectrum consists only of isolated eigenvalues of
	finite multiplicity. Elliptic regularity identifies the spectra of its
	$L^2$- and $X_0$-realizations. The preceding estimate then gives
	a positive spectral gap,
	\[
	\inf\{
	\operatorname{Re}\lambda:
	\lambda\in\sigma(\mathscr A_{\bar m})
	\}
	\geq
	2\mu q(\rho)>0.
	\]
	Consequently, $-\mathscr A_{\bar m}$ generates an exponentially stable
	analytic semigroup on $X_0$.
	For some $\omega>0$, exponentially weighted maximal regularity then yields
	\begin{align}
		&\sup_{t\geq0}e^{\omega t}\norm{u(t)}{H^s}
		+
		\norm{e^{\omega t}\partial_t u}
		{L^2(0,\infty;H^{s-1})}
		+
		\norm{e^{\omega t}u}
		{L^2(0,\infty;H^{s+1})}
		\nonumber\\
		&\qquad\leq
		C\left(
		\norm{u(0)}{\mathcal{X}_s}
		+
		\norm{
			e^{\omega t}
			(\partial_t u+\mathscr A_{\bar m}u)
		}{L^2(0,\infty;H^{s-1})}
		\right);
		\label{equ:stable maximal regularity}
	\end{align}
	see, e.g.~\cite[Proposition~3.3]{PruSimZac09}.
	
	Writing $u=m-\bar m$ and using $h(\bar m)=0$, we rewrite the full
	nonlinear LLB equation \eqref{equ:closed ball LLB} as a perturbation of its linearisation,
	\[
	\partial_t u+\mathscr A_{\bar m}u
	=
	\mathcal N_{\bar m}(u),
	\]
	where
	\begin{align*}
		\mathcal N_{\bar m}(u)
		=
		-\bigl(
		\mathbb A(\bar m+u)-\mathbb A(\bar m)
		\bigr)\mathscr L_{\bar m}u
		-
		\theta\mathbb A(\bar m+u)
		\bigl(
		G(\bar m+u)-G(\bar m)-DG(\bar m)u
		\bigr).
	\end{align*}
	For $u$ in a sufficiently small $H^s$-neighbourhood of zero,
	$\bar m+u$ remains in a fixed compact subset of the open unit ball.
	The Sobolev multiplier estimate and Taylor's theorem then give
	\begin{align}
		\norm{\mathcal N_{\bar m}(u)}{H^{s-1}}
		&\leq
		C\norm{u}{H^s}\norm{u}{H^{s+1}}
		+
		C\norm{u}{H^s}^2
		\leq
		C\norm{u}{H^s}\norm{u}{H^{s+1}}.
		\label{equ:stable quadratic remainder}
	\end{align}
	The same estimates for differences yield
	\begin{align}
		\norm{
			\mathcal N_{\bar m}(u)-\mathcal N_{\bar m}(v)
		}{H^{s-1}}
		\leq
		C\bigl(\norm{u}{H^s}+\norm{v}{H^s}\bigr)
		\norm{u-v}{H^{s+1}}
		+
		C\bigl(\norm{u}{H^{s+1}}+\norm{v}{H^{s+1}}\bigr)
		\norm{u-v}{H^s}.
		\label{equ:stable quadratic difference}
	\end{align}
	Introduce the exponentially weighted maximal-regularity norm
	\[
	\norm{u}{\mathbb E_\omega}
	\coloneq
	\sup_{t\geq0}
	e^{\omega t}\norm{u(t)}{H^s}
	+
	\norm{e^{\omega t}\partial_tu}
	{L^2(0,\infty;H^{s-1})}
	+
	\norm{e^{\omega t}u}
	{L^2(0,\infty;H^{s+1})}.
	\]
	The preceding estimates \eqref{equ:stable quadratic remainder} and \eqref{equ:stable quadratic difference} imply, for $u,v$ in a sufficiently small
	$\mathbb E_\omega$-ball,
	\begin{align*}
		\norm{e^{\omega t}\mathcal N_{\bar m}(u)}
		{L^2(0,\infty;H^{s-1})}
		&\leq
		C\norm{u}{\mathbb E_\omega}^2,
		\\
		\norm{
			e^{\omega t}
			\bigl(
			\mathcal N_{\bar m}(u)-\mathcal N_{\bar m}(v)
			\bigr)
		}
		{L^2(0,\infty;H^{s-1})}
		&\leq
		C
		\bigl(
		\norm{u}{\mathbb E_\omega}
		+
		\norm{v}{\mathbb E_\omega}
		\bigr)
		\norm{u-v}{\mathbb E_\omega}.
	\end{align*}
	
	Now, let $\mathcal T(u)$ denote the solution of
	\[
	\partial_t\mathcal T(u)
	+
	\mathscr A_{\bar m}\mathcal T(u)
	=
	\mathcal N_{\bar m}(u),
	\qquad
	\mathcal T(u)(0)=m_0-\bar m.
	\]
	Applying \eqref{equ:stable maximal regularity} to
	$\mathcal T(u)$, with forcing $\mathcal N_{\bar m}(u)$ and initial
	datum $m_0-\bar m$, then subsequently to
	$\mathcal T(u)-\mathcal T(v)$, whose initial datum vanishes, gives
	\begin{align*}
		\norm{\mathcal T(u)}{\mathbb E_\omega}
		&\leq
		C_0\norm{m_0-\bar m}{\mathcal X_s}
		+
		C_1\norm{u}{\mathbb E_\omega}^2,
		\\
		\norm{\mathcal T(u)-\mathcal T(v)}{\mathbb E_\omega}
		&\leq
		C_1
		\bigl(
		\norm{u}{\mathbb E_\omega}
		+
		\norm{v}{\mathbb E_\omega}
		\bigr)
		\norm{u-v}{\mathbb E_\omega}.
	\end{align*}
	We now choose the fixed-point ball. First, select $R>0$ sufficiently
	small that
	\[
	C_1R\leq\frac14
	\qquad\text{and}\qquad
	C_{\rm emb}R
	\leq
	\frac12
	\bigl(
	1-\norm{\bar m}{L^\infty(\mathscr D)}
	\bigr),
	\]
	where $C_{\rm emb}$ is the norm of the embedding
	$H^s(\mathscr D)\hookrightarrow L^\infty(\mathscr D)$.
	Then, for $u,v$ in the closed ball $B_R \coloneq
	\left\{u\in\mathbb E_\omega: \norm{u}{\mathbb E_\omega}\leq R \right\}$,
	we have
	\[
	\norm{\mathcal T(u)-\mathcal T(v)}{\mathbb E_\omega}
	\leq
	2C_1R
	\norm{u-v}{\mathbb E_\omega}
	\leq
	\frac12
	\norm{u-v}{\mathbb E_\omega}.
	\]
	Hence, $\mathcal T$ is a contraction on $B_R$.
	
	Having fixed $R$, choose the initial perturbation sufficiently small
	that $\norm{m_0-\bar m}{H^s} \leq R/(2C_0)$.
	Then, for every $u\in B_R$,
	\[
	\norm{\mathcal T(u)}{\mathbb E_\omega}
	\leq
	C_0\norm{m_0-\bar m}{\mathcal X_s}
	+
	C_1R^2
	\leq
	\frac R2+\frac R4
	<R,
	\]
	so $\mathcal T$ maps $B_R$ into itself. Moreover, by the choice of
	$R$,
	\[
	\norm{\bar m+u(t)}{L^\infty(\mathscr D)}
	\leq
	\norm{\bar m}{L^\infty(\mathscr D)}
	+
	C_{\rm emb}R
	<1
	\qquad
	\text{for all }t\geq 0,
	\]
	so the whole fixed-point ball remains uniformly separated from saturation.
	Hence, it has a unique fixed point $u$ on $[0,\infty)$. Moreover,
	\[
	\norm{u}{\mathbb E_\omega}
	\leq
	C\norm{m_0-\bar m}{\mathcal X_s},
	\]
	which is precisely \eqref{equ:strong exponential stability}.
	The choice of $R$ also ensures
	\[
	\norm{\bar m+u(t)}{L^\infty(\mathscr D)}<1
	\qquad
	\text{for all }t\geq0,
	\]
	so the solution remains uniformly separated from saturation. This completes the proof.
\end{proof}

As discussed in Section~\ref{subsec:relation reduced LLB}, near the
paramagnetic state $m=0$ the full model formally reduces above the bulk
Curie temperature to the standard reduced LLB equation \eqref{equ:standard above Curie LLB}, for which existence of global strong solutions for $d=3$ is typically proved under a small-data assumption; see, e.g.,~\cite{PuYan22, Soe26c}. For our problem, the condition $3\theta>\Lambda_\ast$ provides the corresponding stability condition for the zero equilibrium. Motivated by this connection, the following corollary gives global strong well-posedness for sufficiently small perturbations of $m=0$, in the Sobolev scale $s\in(3/2,2]$.

\begin{corollary}[Small initial magnetization]
	\label{cor:strong small magnetization}
	Let $d=3$, $s\in (3/2, 2]$, and suppose $3\theta>\Lambda_\ast$,
	where $\Lambda_\ast$ was defined in~\eqref{equ:critical spectral value}, or equivalently $\theta>\theta_{\mathrm c}$.
	Then, sufficiently small initial data in $H_n^s(\mathscr D)$
	generate a unique global strong solution satisfying
	\eqref{equ:strong exponential stability} with $\bar m=0$.
\end{corollary}

\begin{proof}
	Since $G(0)=0$ and $DG(0)=3I$, zero is an equilibrium and
	\begin{align*}
		\inpro{\mathscr L_0z}{z}
		&=
		\alpha\norm{\nabla z}{L^2}^2
		-
		\inpro{\mathsf Bz}{z}
		-
		\inpro{H_{\mathrm d}[z]}{z}
		+
		3\theta\norm{z}{L^2}^2
		\geq
		(3\theta-\Lambda_\ast)\norm{z}{L^2}^2.
	\end{align*}
	Thus, \eqref{equ:strong equilibrium stability} holds.
	The conclusion then follows from Theorem~\ref{thm:strong stable equilibrium}.
\end{proof}




\appendix

\section{Properties of the nonlinear coefficients}

We record several elementary properties of $q$ and $\mathbb M$ which are used in our analysis.

\begin{lemma}\label{lem:properties q and mobility}
	Let $\mathcal{L}$ be the Langevin function defined in \eqref{equ:L xi} and let $\mathcal{L}^{-1}$ be its inverse. Then, the function
	\[
	q(s)=\frac{s}{\mathcal L^{-1}(s)},
	\qquad 0<s<1,
	\]
	admits a smooth extension to $[0,1]$ satisfying $q(0)=\frac13$ and $q(1)=0$, and it is strictly decreasing on $[0,1]$. Moreover, there exist constants
	$c_q,C_q>0$ such that
	\begin{equation}
		c_q(1-s^2)
		\leq
		q(s)
		\leq
		C_q(1-s^2)
		\qquad
		\text{for all }s\in[0,1].
		\label{equ:q saturation equivalence}
	\end{equation}
\end{lemma}

\begin{proof}
	The Maclaurin expansion of the Langevin function is
	\[
	\mathcal L(\xi)
	=
	\frac{\xi}{3}
	-
	\frac{\xi^3}{45}
	+
	O(\xi^5)
	\qquad\text{as }\xi\to0.
	\]
	Since $\mathcal L'(0)=1/3$, its local inverse is analytic at the
	origin and satisfies
	\[
	\mathcal L^{-1}(s)
	=
	3s+\frac95s^3+O(s^5).
	\]
	Consequently,
	\[
	q(s)
	=
	\frac13-\frac15s^2+O(s^4)
	\qquad\text{as }s\to0,
	\]
	and the analyticity of $\mathcal L^{-1}$ shows that $q$ extends
	smoothly to $s=0$ with $q(0)=1/3$.
	To verify monotonicity, we use the Mittag--Leffler expansion
	\[
	\coth\xi
	=
	\frac1\xi
	+
	2\xi
	\sum_{k=1}^\infty
	\frac{1}{\xi^2+\pi^2k^2}.
	\]
	It follows that
	\[
	\frac{\mathcal L(\xi)}{\xi}
	=
	2
	\sum_{k=1}^\infty
	\frac{1}{\xi^2+\pi^2k^2},
	\]
	which is strictly decreasing for $\xi>0$. Since
	$q(\mathcal L(\xi))
	=
	\mathcal L(\xi)/\xi$ and $\mathcal L$ is strictly increasing, $q$ is strictly decreasing.
	Furthermore, by noting that
	$\coth \xi=1+2e^{-2\xi}(1-e^{-2\xi})^{-1}$, we obtain
	\[
	\mathcal L(\xi)
	=
	1-\frac1\xi+O(e^{-2\xi})
	\qquad\text{as }\xi\to\infty.
	\]
	To see that the extension at $s=1$ is smooth, set $u\coloneq1/\xi$ and define
	\[
	F(u)
	\coloneq
	\mathcal L(1/u)
	=
	1-u+r(u),
	\qquad
	r(u)
	=
	\frac{2}{e^{2/u}-1}.
	\]
	The function $r$ extends smoothly to $u=0$ by $r(0)=0$, since
	$r(u)=O(e^{-2/u})$ and all derivatives of $r$ vanish as $u\downarrow0$.
	Since $F(0)=1$ and $F'(0)=-1$, the inverse function theorem yields a
	smooth inverse $u=u(s)$ near $s=1$, with $u(1)=0$. For $s<1$, this
	inverse satisfies $u(s)=1/\mathcal L^{-1}(s)$. Consequently, $q(s)= su(s)$,
	so $q$ extends smoothly to $s=1$ with $q(1)=0$.
	
	In particular, as $s\uparrow1$, we have
	\[
	\mathcal L^{-1}(s)
	\sim
	\frac1{1-s},
	\qquad
	q(s)
	\sim
	s(1-s).
	\]
	Hence,
	\[
	\frac{q(s)}{1-s^2}
	\to
	\frac12
	\qquad\text{as }s\uparrow1.
	\]
	The quotient $q(s)/(1-s^2)$ thus extends to a positive continuous
	function on $[0,1]$, with values $1/3$ at $s=0$ and $1/2$ at $s=1$.
	Its positive minimum and finite maximum yield
	\eqref{equ:q saturation equivalence}.
\end{proof}

\begin{lemma}
	\label{lem:properties M}
	The matrix $\mathbb M$ defined by~\eqref{equ:nonlinear LLB mobility} admits an extension $\mathbb M \in C^{1,1} \bigl(\overline{B_1(0)}; \mathbb R^{3\times3}\bigr)$, with
	\[
	\mathbb M(0)=\frac23\mathsf I,
	\qquad
	\mathbb M(v)=\mathsf I-v\otimes v
	\quad\text{if }|v|=1.
	\]
	Furthermore, for every $v\in\overline{B_1(0)}$ and
	$z\in\mathbb R^3$,
	\begin{align}
		\label{equ:mobility bounded}
		z\cdot\mathbb M(v)z
		&\leq
		|z|^2,
		\\
		\label{equ:mobility radial coercivity}
		2q(|v|)|z|^2
		&\leq
		z\cdot\mathbb M(v)z,
		\\
		\label{equ:mobility square estimate}
		|\mathbb M(v)z|^2
		&\leq
		z\cdot\mathbb M(v)z,
		\\
		\label{equ:precession mobility control}
		|v\times z|^2
		&\leq
		\frac32\,
		z\cdot\mathbb M(v)z.
	\end{align}
\end{lemma}

\begin{proof}
	For $v\neq0$, write
	\[
	z=z_\parallel+z_\perp,
	\qquad
	z_\parallel
	=
	\frac{v\otimes v}{|v|^2}z,
	\qquad
	z_\perp=z-z_\parallel.
	\]
	Then, we have
	\begin{equation}
		\label{equ:mobility decomposition proof}
		z\cdot\mathbb M(v)z
		=
		2q(|v|)|z_\parallel|^2
		+
		\bigl(1-q(|v|)\bigr)|z_\perp|^2.
	\end{equation}
	Since $0\leq q\leq1/3$, the eigenvalues $2q(|v|)$ and
	$1-q(|v|)$ belong to $[0,1]$, with
	$2q(|v|)\leq1-q(|v|)$. This proves
	\eqref{equ:mobility bounded} and
	\eqref{equ:mobility radial coercivity}. Since
	$\lambda^2\leq\lambda$ for $\lambda\in[0,1]$, the same
	decomposition yields \eqref{equ:mobility square estimate}.
	Moreover, \eqref{equ:precession mobility control} follows since
	\[
	|v\times z|^2
	=
	|v|^2|z_\perp|^2
	\leq
	|z_\perp|^2
	\leq
	\frac32\,
	z\cdot\mathbb M(v)z,
	\]
	where we used $1-q(|v|)\geq2/3$.
	
	It remains to verify the regularity at the apparent singularity
	$v=0$. Rewriting
	\[
	\mathbb M(v)
	=
	\bigl(1-q(|v|)\bigr)\mathsf I
	+
	\frac{3q(|v|)-1}{|v|^2}v\otimes v,
	\]
	and using
	\[
	q(s)
	=
	\frac13-\frac15s^2+O(s^4)
	\qquad\text{as }s\to0,
	\]
	we obtain
	\[
	\mathbb M(v)
	=
	\frac23\mathsf I
	+
	\frac15|v|^2\mathsf I
	-
	\frac35v\otimes v
	+
	O(|v|^4).
	\]
	In particular, the singularity is removable and
	$\mathbb M(0)=\frac23\mathsf I$. More generally, the expansion of
	$q$ at the origin contains only even powers of $s$, so the functions
	\[
	s\mapsto1-q(s),
	\qquad
	s\mapsto\frac{3q(s)-1}{s^2}
	\]
	extend smoothly across $s=0$ as even functions. Together with the
	regular extension of $q$ at $s=1$, this shows that $\mathbb M$ is
	$C^2$ on $\overline{B_1(0)}$, and hence
	$\mathbb M\in C^{1,1}(\overline{B_1(0)};\mathbb R^{3\times3})$.
\end{proof}

We also recall the continuity properties of the physical stray field.

\begin{lemma}
	\label{lem:stray field mapping}
	Let $p\in (1,\infty)$ and $s\in \{0,1,2\}$. The stray-field linear operator $H_{\mathrm d}$ satisfies
	\begin{align}\label{equ:stray field Lp}
		\norm{H_{\mathrm d}[v]}{L^p(\mathbb R^3)}
		&\leq
		C_p\norm{v}{L^p(\mathscr D)},
		\\
		\label{equ:stray field W2r}
		\norm{H_{\mathrm d}[v]}{W^{s,p}(\mathscr D)}
		&\leq
		C_{s,p} \norm{v}{W^{s,p}(\mathscr D)}.
	\end{align}
\end{lemma}

\begin{proof}
	See~\cite[Lemma~2.3]{CarFab01} and \cite[Section~3.1]{DumLab12}.
\end{proof}

\makeatletter
\begingroup
\let\@tocwrite\@gobbletwo

\section*{Statements and declarations}

\subsection*{Conflict of interest}
The author has no competing interests to declare that are relevant to the content of this article.

\subsection*{Funding statement}
This research was initiated while the author was supported by the Commonwealth through an Australian Government Research Training Program (RTP) Scholarship (\href{https://doi.org/10.82133/C42F-K220}{DOI: 10.82133/C42F-K220}) at UNSW Sydney. The author has also been supported by an Australian Mathematical Society Lift-off Fellowship. The work was completed during the author's postdoctoral appointment at the Institute of Analysis and Scientific Computing of TU Wien, supported by the Austrian Science Fund (FWF) through the international project I6802 ``Functional error estimates for PDEs on unbounded domains'' (\href{https://doi.org/10.55776/I6802}{DOI: 10.55776/I6802}), led by Prof. Dirk Praetorius.

\endgroup
\makeatother

\end{document}